\documentclass{article}
\usepackage{convexity}

\begin{document}

\title{Semiconcave functions in Alexandrov's geometry}

\author{Anton~Petrunin\footnote{\it Supported in part by the National Science Foundation under grant \# DMS-0406482.
}}
\date{}

\maketitle

\begin{abstract}
The following is a compilation of some techniques in Alexandrov's geometry which are directly connected to convexity.
\end{abstract}

\renewcommand{\thefootnote}{}

\renewcommand{\thefootnote}{\arabic{footnote}}

\section*{Introduction}

This paper is not about results, it is about available techniques in Alexandrov's geometry which are linked to semiconcave functions.  
We consider only spaces with lower curvature bound, but
most techniques described here also work for upper curvature bound and even in more general settings. 

Many proofs are omitted, I include only those which necessary for a continuous story and some easy ones.
The proof of the existence of quasigeodesics is included in appendix~\ref{constr-qg} (otherwise it would never be published).

I did not bother with rewriting basics of Alexandrov's geometry but I did change notation, so it does not fit exactly in any introduction.
  
I tried to make it possible to read starting from any place. 
As a result the dependence of statements is not linear, some results in the very beginning depend on those in the very end and the other way around (but there should not be any cycle).

Here is a list of available introductions to Alexandrov's geometry: 
\begin{enumerate}[$\diamond$]
\item The original paper \cite{BGP}
by Burago, Gromov and Perelman 
and its extension
\cite{perelman:spaces2} by Perelman 
is the first introduction to Alexandrov's geometry. 
I use it as the main reference.

\item  A reader friendly introduction to Alexandrov's geometry \cite{shiohama} by Shiogama. 

\item A survey in Alexandrov's geometry \cite{plaut:survey}
written for topologists by Plaut. 
The first 8 sections can be used as an introduction. 
The material covered in my paper is closely related to sections 7--10 of the   Plaut's survey.

\item Chapter 10 in the book \cite{BBI} by Burago, Burago and Ivanov is yet an other reader friendly introduction.
\end{enumerate}

I want to thank Karsten Grove for making me write this paper, Stephanie Alexander, Richard Bishop, Sergei Buyalo, Vitali Kapovitch, Alexander Lytchak and Conrad Plaut for many useful discussions during its preparation and correction of mistakes, Irina Pugach for correcting my English.

\tableofcontents

\subsection{Notation and conventions}
\begin{enumerate}[$\diamond$]
\item By $\Alex^m(\kappa)$ we will denote the class of $m$-dimensional Alexandrov's spaces with curvature $\ge \kappa$. 
In this notation we may omit $\kappa$ and $m$, but if not stated otherwise we assume that dimension is finite.

\item Gromov--Hausdorff convergence is understood with fixed sequence of approximations. 
That is, 
once we write $X_n\GHto X$ that means that we fixed a sequence of Hausdorff approximations $f_n\:X_n\to X$ (or equivalently $g_n\:X\to X_n$).

This makes possible to talk about limit points in $X$ for a sequence $x_n\in X_n$, limit of functions $f_n\:X_n\to \RR$, Hausdorff limit of subsets $S_n\i X_n$ as well as weak limit of measures $\mu_n$ on $X_n$.

\item regular fiber --- see page~\pageref{reg-fib}

\item $\angle x y z$ --- angle at $y$ in a geodesic triangle $\triangle x y z\i A$

\item $\angle(\xi,\eta)$ --- an angle between two directions $\xi,\eta\in \Sigma_p$

\item $\tilde\angle_\kappa x y z$ --- a comparison angle; 
that is the angle of the model
triangle $\tilde\triangle x y z$ in $\Lob_\kappa$ at $y$.

\item $\tilde\angle_\kappa(a,b,c)$ --- an angle opposite $b$ of a triangle in $\Lob_\kappa$ with sides $a,b$ and $c$. 
In case $a+b<c$ or $b+c<a$ we assume $\tilde\angle_\kappa(a,b,c)=0$.

\item $\uparrow_p^q$ --- a direction at $p$ of a minimizing geodesic from $p$ to $q$ 

\item $\Uparrow_p^q$ --- the set of all directions at $p$ of minimizing geodesics from $p$ to $q$ 

\item $A$ --- usually an Alexandrov's space

\item $\argmax$ --- see page~\pageref{argmax}

\item $\partial A$ --- boundary of $A$

\item $\dist_x(y)=|x y|$  --- distance between $x$ and $y$

\item $\d_p f$ --- differential of $f$ at $p$, see page~\pageref{def:df}

\item $\gexp_p$ --- see section~\ref{gexp} 

\item $\gexp_p(\kappa;v)$ --- see section~\ref{sph-hyp-exp}

\item $\gamma^\pm$ --- right/left tangent vector, see \ref{def-grad-curv}

\item $\Lob_\kappa$ --- model plane see page~\pageref{Lob_k}

\item $\Lob_\kappa^+$ --- model halfplane see page~\pageref{Lob_k^+}

\item $\Lob_\kappa^m$ --- model $m$-space, see page~\pageref{lob-k-m}

\item $\log_p$ --- see page~\pageref{log}

\item $\nabla_p f$ --- gradient of $f$ at $p$, see definition~\ref{def:grad}

\item $\rho_\kappa$ --- see page~\pageref{rho_k}.

\item $\Sigma(X)$ --- the spherical suspension over $X$ see \cite[4.3.1]{BGP}, in \cite[89]{plaut:survey} and \cite{berestovskii} it is called \emph{spherical cone}.

\item $\sigma_\kappa$ --- see footnote~\ref{sigma_k} on page~\pageref{sigma_k}.

\item $T_p=T_p A$ --- tangent cone at $p\in A$, see page~\pageref{T_p}.

\item $T_p E$ --- see page~\pageref{T_pE}

\item $\Sigma_p=\Sigma_p A$ --- see footnote~\ref{U_p} on page~\pageref{U_p}.

\item $\Sigma_p E$ --- see page~\pageref{U_pX}

\item $f^\pm$ --- see page~\pageref{def-grad-curv}

\end{enumerate}

\newpage


\section{Semiconcave functions.}\label{CF}

\subsection{Definitions}

\begin{thm}{\bf Definition for a space without boundary.}\label{def:with-no-bry} 
Let $A\in\Alex$, 
$\partial A=\emptyset$ 
and $\Omega\i A$ be an open subset.

A locally Lipschitz function $f\colon\Omega \rightarrow \RR$ is called
$\lambda$-\emph{concave} (briefly $f''\le \lambda$) if for any unit-speed geodesic $\gamma$ 
in $\Omega$, the function
$$f\circ\gamma(t)-\tfrac\lambda2{\cdot}t^2$$
is concave.
\end{thm}

If $A$ is an Alexandrov's space with non-empty boundary\footnote{Boundary of Alexandrov's space is defined in \cite[7.19]{BGP}.}, 
then its \emph{doubling} $\tilde A$
(that is, two copies of $A$ glued along their boundaries) 
is also an Alexandrov's space (see \cite[5.2]{perelman:spaces2}) and $\partial \tilde A=\emptyset$.

Set $\mathtt{p}\:\tilde A\to A$ to be the canonical map. 

\begin{thm}{\bf Definition for a space with boundary.}\label{def:with-bry}  Let $A\in\Alex$, $\partial A\not=\emptyset$ and $\Omega\i A$ be an open subset.

A locally  Lipschitz function $f\colon \Omega \rightarrow \RR$ is called $\lambda$-\emph{concave} (briefly $f''\le \lambda$) if $ f\circ \mathtt{p}$ is $\lambda$-{concave} in $\mathtt{p}^{-1}(\Omega)\i \tilde A$.
\end{thm}

\noi{\bf Remark.} Note that the restriction of a linear function on $\RR^n$ to a ball is not $0$-concave in this sense.

\subsection{Variations of definition.}\label{variation-CF}

A function $f\:A\to \RR$ is called \textit{semiconcave} if for any point $x\in A$ there is a neighborhood $\Omega_x\ni x$ and $\lambda\in \RR$ such that the restriction $f|_{\Omega_x}$ is $\lambda$-concave.

Let $\phi\:\RR\to \RR$ be a continuous function. 
A function $f\:A\to \RR$ is called \emph{$\phi(f)$-concave} (briefly $f''\le \phi(f)$) if for any point $x\in A$ and any $\eps>0$ there is a neighborhood $\Omega_x\ni x$ such that $f|_{\Omega_x}$ is $(\phi\circ f(x)+\eps)$-concave

For the Alexandrov's spaces with curvature $\ge \kappa$, it is natural to consider the functions satisfying $f''+\kappa{\cdot}f\le 1$.
The advantage of such functions comes from the fact that on the model space%
\footnote{\label{Lob_k} that is, 
the simply connected 2-manifold of constant curvature $\kappa$ 
(the Russian $\Lobm$ is for Lobachevsky)}  
$\Lob_\kappa$,  one can construct model $(1-\kappa{\cdot} f)$-concave functions which are equally concave in all directions at any fixed point.
The most important example of $(1-\kappa{\cdot} f)$-concave function is $\rho_\kappa\circ\dist_x$, where $\dist_x(y)=|x y|$ denotes distance function from $x$ to $y$ and
$$\label{rho_k}\rho_\kappa(x) =
\l[\begin{matrix} 
\tfrac1\kappa{\cdot} (1-\cos(x{\cdot}\sqrt{\kappa} ))  & {\text{if }} & \kappa > 0 \\
 x^2/2                          & {\text{if }} & \kappa=0 \\
\tfrac1\kappa {\cdot}(\cosh(x{\cdot}\sqrt{-\kappa})-1) & {\text{if }} & k < 0
\end{matrix}\r.$$

In the above definition of $\lambda$-concave function one can exchange Lipschitz continuity for usual continuity. 
Then it will define the same set of functions, see corollary~\ref{cor:cont-conv}.

\subsection{Differential} 

Given a point $p$ in an Alexandrov's space  $A$, we denote by \label{T_p} $T_p=T_p A$ the tangent cone at $p$.

For an Alexandrov's space, the tangent cone can be defined in two equivalent ways (see \cite[7.8.1]{BGP}):
\begin{enumerate}[$\diamond$]
\item As a cone over space of directions at a point
\item \label{tangent-def2}
As a limit of rescalings of the Alexandrov's space, that is:

Given $\lam>0$, we denote the space by $\lam{\cdot} A$, the rescaling of $A$ by factor $\lam$.
Let $i_\lam\: \lam{\cdot} A\to A$ be the canonical map. 
The limit of $(\lam{\cdot} A,p)$ for $\lam\to\infty$ is  the tangent cone $(T_p,o_p)$ at $p$ with marked origin $o_p$. 

\end{enumerate}

\begin{thm}{\bf Definition.}\label{def:df} Let $A\in\Alex$ and $\Omega\i A$ be an open subset.

For any function $f\:\Omega \rightarrow \RR$ the function $\d_p f\colon T_p \rightarrow \RR$, $p\in\Omega$ defined by
$$\d_p f=\lim_{\lam\to\infty} \lam{\cdot}(f\circ i_\lam-f(p)), \ \ f\circ i_\lam\:\lam{\cdot} A\to\RR$$
is called the \emph{differential} of $f$ at $p$.
\end{thm}

It is easy to see that the differential $\d_p f$ is well defined for any semiconcave function $f$. 
Moreover, $\d_p f$
is a concave function on the tangent cone $T_p$ which is positively homogeneous; 
that is
$d_p f(r\cdot v)=r\cdot d_p f(v)$
for $r\ge0$.

\paragraph*{Gradient.} With a slight abuse of notation, we will call elements of the tangent cone $T_p$ the ``tangent vectors'' at $p$.
The origin $o=o_p$ of $T_p$ plays the role of a ``zero 
vector''.
For a tangent vector $v$ at $p$ we define its absolute value $|v|$ 
as the distance $|o v|$ in $T_p$.
For two tangent vectors $u$ and $v$ at $p$ we can define 
their ``scalar product'' 
$$\langle u, v\rangle\df
=
(|u|^2+|v|^2-|u v|^2)/2
=
|u|\cdot |v|\cdot\cos\alpha,$$ 
where $\alpha=\angle u o v=\tilde\angle_0 u o v$ in $T_p$.

 It is easy to see that for any $u\in T_p$, the function $x\mapsto -\langle u, x\rangle$ on $T_p$ is concave.

\begin{thm}{\bf Definition.}\label{def:grad} 
Let $A\in\Alex$ and $\Omega\i A$ be an open subset.
Given a $\lambda$-concave function
$f\:\Omega\to\RR$, 
a vector $g\in T_p$ is called a \emph{gradient} of $f$ at $p\in \Omega$  
(in short:  $g=\nabla_p f$) if

(i) $\d_p f(x)\le \langle g , x\rangle\ \hbox{for any}\ x\in T_p$, and

(ii) $\d_p f(g) = \langle g,g \rangle .$
\end{thm}

It is easy to see that any $\lambda$-concave function $f\:\Omega\to \RR$
has a uniquely defined gradient vector field. 
Moreover, if $\d_p f(x) \le 0$ for all $x\in T_p$,
then $\nabla_p f=o_p$; otherwise,
$$\nabla_p f=\d_p f(\xi_{\max})\cdot\xi_{\max} $$
where $\xi_{\max}\in \Sigma_p$%
\footnote{\label{U_p}By $\Sigma_p\i T_p$ we denote the set of unit
vectors, which we also call directions at $p$. 
The space $(\Sigma_p,\angle)$ with angle metric is an Alexandrov's space with curvature $\ge 1$. $(\Sigma_p,\angle)$ it is also path-isometric to the subset $\Sigma_p\i T_p$.} 
is the (necessarily
unique) 
unit vector for which the function $\d_p f$ attains its maximum.

For two points $p,q\in A$ we denote by $\uparrow_p^q\,\in \Sigma_p$ a direction of a minimizing geodesic from $p$ to $q$. 
Set \label{log} $\log_p q=|p q|\cdot\!\!\uparrow_p^q\in T_p$.
In general, $\uparrow_p^q$ and $\log_p q$ are not uniquely
defined.

The following inequalities describe an important property of the ``gradient
vector field'' which will be used throughout this paper.

\begin{wrapfigure}{r}{30mm}
\begin{lpic}[t(-2mm),b(0mm),r(0mm),l(0mm)]{pics/grad-inq(0.4)}
\lbl[br]{0,55;$p$}
\lbl[tr]{47,-2;$q$}
\lbl[tr]{11,44;$\uparrow_p^q$}
\lbl[r]{39,8;$\uparrow_q^p$}
\lbl[lb]{26,27;$\ell$}
\lbl[bl]{14,53;$\nabla_p f$}
\lbl[tl]{57,2;$\nabla_q f$}
\end{lpic}
\end{wrapfigure}

\begin{thm}{\bf Lemma.} 
\label{lem:grad}
Let $A\in\Alex$ and $\Omega\i A$ be an open subset,
$f\:\Omega\to\RR$ be a $\lambda$-concave function. 
Assume all minimizing
geodesics between $p$ and $q$ belong to $\Omega$, set $\ell=|p q|$. 
Then
\end{thm}

\vspace{-5mm}
$$\<\uparrow_p^q,\nabla_p f\>\ge
{\{f(q)-f(p)-\tfrac\lambda2{\cdot}\ell^2\}}/{\ell},$$
{\it and in particular}
$$
\<\uparrow_p^q,\nabla_p f\>
+
\<\uparrow_q^p,\nabla_q f\>
\ge 
-\lambda{\cdot}\ell.
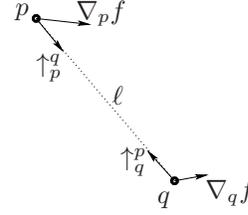
$$

\Proof. 
Let $\gamma\:[0,\ell]\to \Omega$ be a unit-speed minimizing geodesic from $p$ to $q$, so 
$$\gamma(0)=p,\ \ \gamma(\ell)=q,\ \ \gamma^+(0)=\uparrow_p^q.$$
From definition~\ref{def:grad} and the $\lambda$-concavity of $f$ we get
\begin{align*}
\<\uparrow_p^q,\nabla_p f\>
&=
\< \gamma^+(0),\nabla_p f\>
\ge
\\
&\ge
d_p f(\gamma^+(0))=
\\
&=(f\circ\gamma)^+(0)
\ge
\\
&\ge
\frac{f\circ\gamma(\ell)-f\circ\gamma(0)-\tfrac{\lambda}2{\cdot}\ell^2}{\ell}.
\end{align*}
Hence the first inequality follows. 
(For definition of $\gamma^+$ and $(f\circ\gamma)^+$ 
see \ref{def-grad-curv}.)

The second inequality is a sum of two inequalities of the first type. 
\qeds

\begin{thm}{\bf Lemma.} \label{lem:gradcon}
Let $A_n \GHto A$, $A_n\in \Alex^m(\kappa)$.

Let  $f_n\:A_n\to \RR$ be a sequence of $\lambda$-concave 
functions and $f_n\to f\:A\to \RR$. 

Let $x_n\in A_n$ and $x_n\to x\in A$.

Then 
$$|\nabla_x f|\le \liminf_{n\to \infty} |\nabla_{x_n} f_n|.$$
\end{thm}

The corollary below states that the function $x\mapsto |\nabla_x f|$ is lower-semicontinuous.

\begin{thm}{\bf Corollary.} \label{cor:gradlim} Let $A\in\Alex$ and $\Omega\i A$
be an open subset. 

If $f\:\Omega\to\RR$ is a semiconcave function then the function $$x\mapsto|\nabla_x f|$$ 
is lower-semicontinus;
that is,
for any sequence $x_n\to x\in \Omega$, we
have 
$$|\nabla_x f|\le \liminf_{n\to \infty} |\nabla_{x_n} f|.$$

\end{thm}

\Proof\  {\it of lemma~\ref{lem:gradcon}}. Fix an $\eps>0$ and choose $q$ near
$p$ such
that 
$$\frac{f(q)-f(p)}{|p q|}> |\nabla_p f|-\eps.$$
Now choose $q_n\in A_n$ such that $q_n\to q$. 
If $|p q|$ is sufficiently small and $n$ is
sufficiently large,  the $\lambda$-concavity of $f_n$ then implies that
$$\liminf_{n\to \infty}d_{p_n}f_n(\uparrow_{p_n}^{q_n})\ge |\nabla_p f|-2{\cdot}\eps.$$
Therefore
$$\liminf_{n\to \infty} |\nabla_{p_n} f_n|\ge
|\nabla_p f|-2{\cdot}\eps \ \ \text{for any}\ \ \eps>0.$$
Whence the lemma follows.
\qeds

\paragraph*{Supporting and polar vectors.}\label{supp-polar}

\begin{thm}{\bf Definition.}\label{def-support} 
Assume $A\in\Alex$ and $\Omega\i A$ is an open subset, $p\in
\Omega$. 
Let $f\:\Omega\to \RR$ be a semiconcave function.

A vector $s\in T_p$ is called a \emph{supporting vector} of $f$ at $p$ if 
$$\d_p f(x)\le -\langle s , x\rangle\ \ \hbox{for any}\ \ x\in T_p$$
\end{thm}

The following lemma sates that the set of supporting vectors is not empty.

\begin{thm}{\bf Lemma.} 
Assume $A\in\Alex$ and $\Omega\i A$ is an open subset and $f\:\Omega\to \RR$ is a semiconcave function, $p\in \Omega$. 
Then the set of supporting vectors of $f$ at $p$ form a non-empty convex subset of $T_p$.
\end{thm}

\Proof. Convexity of the set of supporting vectors follows from concavity of the function $x\to -\<u,x\>$ on $T_p$.
To show existence, consider a minimum point  $\xi_{\min}\in \Sigma_p$ of the function $\d_p f|_{\Sigma_p}$. 
We will show that the vector
$$s=\l[-\d_p f(\xi_{\min})\r]\cdot\xi_{\min}$$
is a supporting vector for $f$ at $p$.
Assume that we know the existence of supporting vectors for dimensions $<m$.
Applying it to $d_p f|_{\Sigma_p}$ at $\xi_{\min}$, we get $d_{\xi_{\min}}(d_p f|_{\Sigma_p})\equiv 0$.
Therefore, since $d_p f|_{\Sigma_p}$ is $(-d_p f)$-concave 
(see section~\ref{variation-CF}) 
for any $\eta\in \Sigma_p$ we have 
$$d_p f(\eta)\le d_p f(\xi_{\min})\cdot\cos\angle(\xi_{\min},\eta).$$
Hence the result follows. \qeds

In particular, if the space of directions $\Sigma_p$ has
a diameter%
\footnote{We always consider $\Sigma_p$ with angle metric.} 
$\le \tfrac\pi2$ then
$\nabla_p f=o$ for any $\lambda$-concave function $f$.

Clearly, for any vector $s$, supporting  $f$ at $p$ we have 
$$|s|\ge|\nabla_p f|.$$

\begin{thm}{\bf Definition.}\label{defn:polar}
Two vectors $u,v\in T_p$ are called \emph{polar} if for any vector $x\in T_p$ we
have 
$$\<u,x\>+\<v,x\>\ge0.$$
More generally, a vector $u\in T_p$ is called polar to a set of vectors
$\mathcal{V}\i T_p$ if 
$$\<u,x\>+\sup_{v\in \mathcal{V}}\<v,x\>\ge0.$$
\end{thm}

Note that if $u,v\in T_p$
are polar to each other then \label{*-polar-inq}
$$\d_p f(u)+\d_p f(v)\le 0\eqno(*)$$
for any semiconcave function $f$ defined at $p$.
Indeed, if $s$ is a supporting vector then
$$\d_p f(u)+\d_p f(v)\le-\<s,u\>-\<s,v\>\le 0.$$
Similarly, if $u$ is polar to a set $\mathcal{V}$ then \label{**-polar-inq}
$$\d_p f(u)+\inf_{v\in \mathcal{V}}\{\d_p f(v)\}\le 0\eqno(**)$$
for any semiconcave function $f$ defined at $p$.

\noi{\it Examples of pairs of polar vectors.}
\begin{enumerate}[(i)]
\item If two vectors $u,v\in T_p$ are \emph{antipodal};
that is, $|u|=|v|$ and $\angle u o_p v=\pi$ then they are polar to each other. 

In general, if $|u|=|v|$ then they are polar if and only if for any $x\in T_p$
we have $\angle u o_p x+\angle x o_p v\le\pi$.

\item \label{polar} If $\uparrow_q^p$ is uniquely defined then $\uparrow_q^p$ is
polar to $\nabla_q\dist_p$. 

More generally, if $\Uparrow_p^q\i \Sigma_p$ denotes the set of all directions from $p$
to $q$ then $\nabla_q\dist_p$ is polar to the set $\Uparrow_q^p$. 
\end{enumerate}
Both statement follow from the identity 
$$d_q(v)=\min_{\xi\in \Uparrow_q^p}\{-\<\xi,v\>\}$$
and the definition of gradient (see \ref{def:grad}).

Fix a vector $v\in T_p$.
Let us apply above property~(\ref{polar}) to the function
$f_v=\dist_v\:T_p\z\to\RR$. 
We get that $\nabla_o f_v$ is polar to $\uparrow_o^v$.
Since there is a natural isometry $T_o T_p\to T_p$ we get the following.

\begin{thm}{\bf Lemma.}\label{lem:polar} Given any vector $v\in T_p$ there is a
polar vector $v^*\in T_p$. 
Moreover, one can assume that $|v^*|\le |v|$
\end{thm}

In \ref{lem:milka} using quasigeodesics we will show that in fact one can assume
$|v^*|=|v|$

\section{Gradient curves.}

The technique of gradient curves was influenced by 
Sharafutdinov's retraction introduced in \cite{sharafutdinov}. 
These curves were designed to simplify Perelman's proof of
existence of quasigeodesics. 
However, it turned out that gradient curves themselves
provide a superior tool, which is in fact almost universal in Alexandrov's
geometry.
Unlike most of Alexandrov's techniques, gradient curves work equally well for
infinitely dimensional Alexandrov's spaces (the proof requires some quasifications). 
As it was shown by Lytchak in \cite{lytchak:open-map},
this technique also works for spaces with curvature bounded above and for locally compact
spaces with well defined tangent cone at each point. 
Some traces of these
properties can be found even in general metric spaces, 
see \cite{grad-flow-book}.

\subsection{Definition and main properties}\label{def-grad-curv}

Given a curve $\gamma(t)$ in an Alexandrov's space $A$, we denote by $\gamma^+(t)$ 
the right, and by 
$\gamma^-(t)$ the left, tangent vectors to  $\gamma(t)$, where, respectively,
$$\gamma^\pm(t)\in T_{\gamma(t)},\ \ \ 
\gamma^\pm(t)
=
\lim_{\eps\to 0+}
\frac{\log_{\gamma(t)}\gamma(t\pm\eps)}{\eps}.$$ 
This sign convention is not quite standard; 
in particular, for a function
$f\:\RR\z\to\RR$,  its right derivative is equal to $f^+$ and  its left derivative is equal to $-f^-(t)$.
For example, if $f(t)=t$ then $f^+(0)= 1$ and $f^-(0)= -1$.

\begin{thm}{\bf Definition.} Let $A\in\Alex$ and $f\:A\to \RR$ be a semiconcave
function. 

A curve $\alpha(t)$ is called $f$-\emph{gradient curve} if 
$$\alpha^+(t)=\nabla_{\alpha(t)}f$$
for any $t$.
\end{thm}

\begin{thm}{\bf Proposition.}\label{prop:gradlip} Given a $\lambda$-concave
function $f$ on an Alexandrov's space $A$ and a point $p\in A$ there is a unique
gradient curve $\alpha\colon [0,\infty) \rightarrow A$ such that $\alpha(0)=p$. 
\end{thm} 

\label{grad-constr} The gradient curve  can be constructed as a limit of broken
geodesics,
made up of short segments with directions close to the gradient. 
Convergence, uniqueness,
follow from lemma \ref{lem:grad}, while corollary~\ref{cor:gradlim} guarantees
that the limit is indeed
a gradient curve.

\paragraph*{Distance estimates.}

\begin{thm}{\bf Lemma.} 
\label{lem:concave}
Let $A\in\Alex$ and $f\:A\to \RR$ be a
$\lambda$-concave function and $\alpha(t)$ be an $f$-gradient curve.

Assume $\bar\alpha(s)$ is the reparametrization of $\alpha(t)$ by arclength. 
Then $f\circ\bar\alpha$ is  $\lambda$-concave.
\end{thm} 

\Proof. For $s>s_0$,
\begin{align*}
(f\circ\bar\alpha)^+(s_0)
&=
|\nabla_{\bar\alpha(s_0)}f|
\ge
\\
&\ge
d_{\bar\alpha(s_0)}f\l(\uparrow_{\bar\alpha(s_0)}^{\bar\alpha(s)}\r)
\ge
\\
&\ge
\frac{f(\bar\alpha(s))-f(\bar\alpha(s_0))-\tfrac\lambda2{\cdot}
|\bar\alpha(s)\,\bar\alpha(s_0)|^2 }{|\bar\alpha(s)\,\bar\alpha(s_0)|}.
\end{align*}
Therefore, since $s-s_0\ge|\bar\alpha(s)\,\bar\alpha(s_0)|=s-s_0-o(s-s_0)$, we
have 
$$(f\circ\bar\alpha)^+(s_0)\ge
\frac{f(\bar\alpha(s))-f(\bar\alpha(s_0))-\tfrac\lambda2{\cdot}(s-s_0)^2}{s-s_0} +o(s-s_0);$$
that is $(f\circ\bar\alpha)''\le\lambda$.
\qeds  

The following lemma gives a nice parametrization of a gradient
curve (by $\theta_\lambda$) 
so that they start to behave like a geodesic in some comparison inequalities.

\begin{thm}{\bf Lemma.} \label{lem:dist-est}
Let $A\in\Alex$, $f\:A\to \RR$ 
be a $\lambda$-concave function 
and
$\alpha,\beta\:[0,\infty)\to A$ be two $f$-gradient curves with $\alpha(0)=p$,
$\beta(0)=q$. 

Then 
\begin{enumerate}[(i)]
\item\label{two-equal-ends} for any $t\ge0$,
$$|\alpha(t)\beta(t)|\le e^{\lambda{\cdot}t}|p q|$$
\item \label{one-end} for any $t\ge0$,
$$|\alpha(t)q|^2\le|p q|^2+ 
\l\{2{\cdot}f(p)-2{\cdot}f(q)+\lambda{\cdot}
|p q|^2\r\}\cdot\theta_\lambda(t)
+|\nabla_p f|^2\cdot\theta^2_\lambda(t),$$
where 
$$\theta_\lambda(t)=\int_0^t e^{\lambda{\cdot}t}\cdot d t=
\l[
\begin{matrix}
t &\text{if}&\lambda=0\\
\frac{e^{\lambda{\cdot}t}-1}\lambda&\text{if}&\lambda\not=0
\end{matrix}
\r.
$$
 
\item \label{two-ends} if $t_p\ge t_q\ge 0$ then

$|\alpha(t_p)\beta(t_q)|^2
\le e^{2\lambda{\cdot}t_q}
\bigl[|p q|^2+ 
 $

$\ \ \ \ \ \ \ \ \ \ \ \ \ \ \ \ \ \ \ +\l\{2{\cdot}f(p)-2{\cdot}f(q)+\lambda{\cdot}
|p
q|^2\r\}\cdot\theta_\lambda(t_p-t_q)+$

$\ \ \ \ \ \ \ \ \ \ \ \ \ \ \ \ \ \ \ +|\nabla_p f|^2
\cdot\theta^2_\lambda(t_p-t_q)\bigr].$
\end{enumerate}
\end{thm}

In case $\lambda>0$, this lemma can also be reformulated in a geometer-friendly way:

\bigskip

\noi{\bf \ref{lem:dist-est}$'\!\!\!$. Lemma.}
{\it 
Let $\alpha$, $\beta$, $p$ and $q$ be as in lemma~\ref{lem:dist-est} and $\lambda>0$. 
Consider points $\tilde o, \tilde p, \tilde q\i \RR^2$ defined by the
following:
$$|\tilde p\tilde q|=|p q|,\ \ \lambda{\cdot}
|\tilde o\tilde p|=|\nabla_p f|,$$
$$\tfrac\lambda2{\cdot}\l(|\tilde o\tilde q|^2-|\tilde o\tilde p|^2\r)=f(q)-f(p)$$
Let $\tilde\alpha(t)$ and $\tilde\beta(t)$ be 
$\l({\frac\lambda2{\cdot}\dist^2_{\tilde o}}\r)$-gradient curves in $\RR^2$ with
$\tilde \alpha(0)=\tilde p$, $\tilde \beta(0)=\tilde q$.
Then,
\begin{enumerate}[(i)]
\item $|\alpha(t)q|\le |\tilde\alpha(t)\tilde q|$ for any $t>0$
\item $|\alpha(t)\beta(t)|\le|\tilde\alpha(t)\tilde\beta(t)|$
\item if $t_p\ge t_q$ then $|\alpha(t_p)\beta(t_q)|\le
|\tilde\alpha(t_p)\tilde\beta(t_q)|$
\end{enumerate}
}

\begin{wrapfigure}{r}{27mm}
\begin{lpic}[t(-10mm),b(0mm),r(0mm),l(0mm)]{pics/alpha-curve(0.4)}
\lbl[tr]{2,73;$p$}
\lbl[tr]{50,3;$q$}
\lbl[br]{35,90;$\alpha(t)$}
\lbl[br]{68,16;$\beta$}
\end{lpic}
\end{wrapfigure}

\Proof. (\ref{one-end}). 
If $\lambda=0$, then from lemma~\ref{lem:concave} it follows
that%
\footnote{For $\lambda\ne0$ it will be $f\circ\alpha(t)-f\circ\alpha(0)
\le \l|\nabla_{\bar\alpha(0)}f\r|^2\cdot [\theta_\lambda(t)+\tfrac\lambda2{\cdot}\theta_\lambda^2(t)]$.}
$$f\circ\alpha(t)-f\circ\alpha(0)
\le \l|\nabla_{\bar\alpha(0)}f\r|^2\cdot t.$$
Therefore from lemma~\ref{lem:grad}, setting $\ell=\ell(t)=|q\alpha(t)|$, we
get%
\footnote{For $\lambda\not=0$ it will be
$\l({\ell^2}/2\r)'-\tfrac{\lambda}2{\cdot}\ell^2\le f(p)-f(q)+
\l|\nabla_{p}f\r|^2\cdot[\theta_\lambda(t)+\tfrac\lambda2{\cdot}\theta_\lambda^2(t)]$.}
$$\l({\ell^2}/2\r)'\le f(p)-f(q)+ \l|\nabla_{p}f\r|^2\cdot t,$$
hence the result.

\noi(\ref{two-equal-ends}) follows from the second inequality in
lemma~\ref{lem:grad}; 

\noi(\ref{two-ends}) follows from (\ref{two-equal-ends}) and (\ref{one-end}).
\qeds

\paragraph*{Passage to the limit.} 
The next lemma states that gradient curves behave
nicely with Gromov--Hausdorff convergence;
that is, a limit of gradient curves is a
gradient curve for the limit function.

\begin{thm}{\bf Lemma.} \label{lem:stable-grad-curves}
Let $A_n \GHto A$, $A_n\in \Alex^m(\kappa)$, $A_n\ni p_n\to p\in A$.

Let  $f_n\:A_n\to \RR$ be a sequence of $\lambda$-concave functions and $f_n\to
f\:A\to \RR$.
 
Let $\alpha_n\: [0,\infty) \to A_n$ be the sequence of 
$f_n$-gradient curves with $\alpha_n(0)=p_n$ and 
let $\alpha\: [0,\infty) \to A$ be the $f$-gradient curve with $\alpha(0)=p$.

Then $\alpha_n\to\alpha$ as $n\to\infty$.
\end{thm}

\Proof.
Let $\bar\alpha_n(s)$ denote the reparametrization of $\alpha_n(t)$ 
by arc length.
Since all $\bar\alpha_n$ are $1$-Lipschitz, 
 we can choose a partial limit, say $\bar\alpha(s)$ in $A$.
Note
that we may assume that $f$ has no critical points and so
$d(f\circ\bar\alpha)\not=0$. Otherwise consider instead the  sequence
$A'_n=A_n\times\RR$ with $f'_n(a\times x)=f_n(a)+x$.

Clearly, $\bar\alpha$ is also 1-Lipschitz and hence, by Lemma \ref{lem:gradcon},
\begin{align*}
\lim_{n\to\infty}f_n\circ\bar\alpha_n|_a^b
&=
\lim_{n\to\infty}
\int_a^b|\nabla_{\bar\alpha_n(s)} f_n|\cdot\d s\ge
\\
&\ge\int_a^b|\nabla_{\bar\alpha(s)} f|\cdot\d s
\ge 
\\
&\ge
\int_a^b d_{\bar\alpha(s)} f(\bar\alpha^+(s))\cdot\d s
=
\\
&=
f\circ\bar\alpha|_a^b, 
\end{align*}
where $\bar\alpha^+(s)$ denotes any partial limit of
$\log_{\bar\alpha(s)}\bar\alpha(s+\eps)/\eps$, $\eps\to0+$.

On the other hand, 
since $\bar\alpha_n\to\bar\alpha$ 
and $f_n\to f$ 
we have $f_n\circ\bar\alpha_n|_a^b \to f\circ\bar\alpha|_a^b$;
that is, equality holds in both of these inequalities. 
Hence 
$$|\nabla_{\bar\alpha(s)} f|= \lim_{n\to\infty} |\nabla_{\bar\alpha_n(s)} f_n|,\
\ \ 
|\bar\alpha^+(s)|= 1\ \ \ \text{a.e.}$$
and the directions of $\bar\alpha^+(s)$ and  $\nabla_{\bar\alpha(s)} f$ 
coincide almost everywhere.

This implies that $\bar\alpha(s)$ is a gradient curve reparametrized by 
arc length. 
It only remains to show that the original 
parameter $t_n(s)$ of $\alpha_n$ converges to the original 
parameter $t(s)$ of $\alpha$.

Notice that $|\nabla_{\bar\alpha_n(s)} f_n|\cdot d t_n=d s$ or 
$d t_n/d s=d s/d(f_n\circ\bar\alpha_n)$. 
Likewise, $d t/d s=d s/d(f\circ\bar\alpha)$. 
Then the convergence $t_n\to t$ follows from the $\lambda$-concavity of 
$f_n\circ\bar\alpha_n$ (see Lemma~\ref{lem:concave}) 
and the convergence $f_n\circ\bar\alpha_n\to f\circ\bar\alpha.$\qeds

\subsection{Gradient flow}\label{grad-flow}

Let $f$ be a semi-concave function on an Alexandrov's space $A$. 
We define the $f$-\emph{gradient flow} to be the one parameter family of maps 
$$\Phi^t_f\:A\to A, \ \ \Phi^t_f(p)=\alpha_p(t),$$
where $t\ge 0$ and $\alpha_p\:[0,\infty)\to A$ is the $f$-gradient curve which
starts at $p$ (that is, $\alpha_p(0)=p$).
\footnote{In general the domain of definition of $\Phi^t_f$ can be smaller than $A$, 
but it is defined on all $A$ for a reasonable type of function, say for $\lambda$-concave and for $(1-\kappa{\cdot}f)$-concave functions.}
Obviously
$$ \Phi^{t+\tau}_f=\Phi^t_f\circ\Phi^\tau_f.$$
This map has the following main properties:
\begin{enumerate}
\item $\Phi^t_f$ is locally Lipschitz (in the domain of definition). 
Moreover, if $f$ is $\lambda$-concave then it is $e^{\lambda{\cdot}t}$-Lipschitz.

This follows from lemma~\ref{lem:dist-est}(\ref{two-equal-ends}).

\item Gradient flow is stable under Gromov--Hausdorff convergence, namely:

If $A_n\in \Alex^m(\kappa)$, $A_n\GHto A$, $f_n\:A_n\to\RR$ is a sequence of
$\lambda$-concave functions which converges to $f\:A\to \RR$ then
$\Phi_{f_n}^t\:A_n\to A_n$ converges pointwise to $\Phi_f^t\:A\to A$.

This follows  from lemma~\ref{lem:stable-grad-curves}.

\item\label{grad-onto} For any $x\in A$ and all sufficiently small $t\ge 0$, there is $y\in A$ so that 
$\Phi_f^t(y)=x$.

For spaces without boundary this follows from Lemma 1 proved by Grove and Petersen in \cite{grove-petersen:rad-sphere}.
For spaces with boundary one should pass to its doubling.
\end{enumerate}

Gradient flow can be used to deform a mapping with target in $A$. 
For example, if $X$ is a metric space, then given a Lipschitz map $F\:X\to A$ and
a positive Lipschitz function $\tau\:X\to \RR_+$ one can consider the map $F'$ called
\emph{gradient deformation} of $F$ which is defined by
$$F'(x)=\Phi_f^{\tau(x)}\circ F(x),\ \ \ F'\:X\to A.$$

From lemma~\ref{lem:dist-est} it is easy to see that the \emph{dilation}\footnote{That is the optimal Lipschitz constant.}
of $F'$
can be estimated in terms of $\lambda$, $\sup_x\tau(x)$, dilation of $F$ and the
Lipschitz constants of $f$ and $\tau$.

Here is an optimal estimate for the length element of a curve which follows from
lemma~\ref{lem:dist-est}:

\begin{thm}{\bf Lemma.} \label{lem:grad-variation} Let $A\in\Alex$.
Let $\gamma_0(s)$ be a curve in $A$ parametrized by arc-length, 
$f\:A\to\RR$ be a $\lambda$-concave function, 
and $\tau(s)$ be a non-negative Lipschitz function. 
Consider the curve 
$$\gamma_1(s)=\Phi^{\tau(s)}_f \circ\gamma_0(s).$$ 
If $\sigma=\sigma(s)$ is its
arc-length parameter then
$$\d\sigma^2\le e^{2\lambda\tau}\l[\d
s^2+2\cdot d(f\circ\gamma_0)\d\tau+|\nabla_{\gamma_0(s)}f|^2\cdot\d\tau^2\r]$$
\end{thm}

\subsection{Applications} 

\begin{thm}{\bf Toponogov's splitting Theorem.}\label{thm:splitting} Let $A\in\Alex(0)$, and $\gamma:\RR\to A$ be a line 
(that is, a unit-speed geodesic which is minimizing on each segment). 
Then there is an isometry $h:A\to \RR\times A'$ where $A'\in\Alex(0)$.
Moreover, $i$ can be chosen on such a way that if $\pi_\RR$ denotes the projection $\pi:\RR\times A'\to \RR$, $\pi((t,x))=t$ then $\pi\circ h\circ\gamma(t)=t$.
\end{thm}

For smooth $2$-dimensional surfaces, 
this theorem was proved by Cohn-Vossen in \cite{cohn-vossen_line}.
For the Riemannian manifolds of higher dimensions 
it was proved by Toponogov in \cite{toponogov-globalization+splitting}.
Then it was generalized by Milka in  \cite{milka-line}
to Alexandrov's spaces, almost the same prove is given in \cite[10.3]{BBI}.
The Milka's proof is a beautiful Euclid-like argument, but it is not short;
here we present a short proof which use gradient flow for Busemann's functions.

\Proof. Consider two Buseman's functions $b_+$ and $b_-$ associated with rays%
\footnote{Let $X$ be a metric space, a unit-speed geodesic $\gamma:[0,\infty)\to X$ is called a ray if it is minimizing on each bounded segment.} $\gamma:[0,\infty)\to A$ and $\gamma:(-\infty,0]\to A$;
that is,
$$b_\pm(x)=\lim_{t\to\infty}|\gamma(\pm t)\,x|- t.$$
(It converges since for any $x$, $|\gamma(\pm t)\,x|- t$ decreases in $t$).
For any $R>0$ and $\eps>0$ there is $T$ such that for $t>T$, the function $\dist_{\gamma(\pm t)}$ is $(-\eps)$-concave in $B_R(\gamma(0))$.
Therefore, both functions $b_\pm$ are concave.

Note that since $\gamma$ is a line, we have $b_+(x)+b_-(x)\ge0$ for any $x\in A$.
On the other hand $f(t)={\dist^2_x(\gamma(t))}$ is $2$-concave; in particular, $f(t)\le t^2+at+b$ for some constants $a,b\in\RR$. 
Passing to $t\to\pm\infty$, we get $b_+(x)+b_-(x)\le0$ for any $x\in A$.
Hence $$b_+(x)+b_-(x)\equiv 0.$$

Set $A'=b_+^{-1}(0)\i A$, it forms an Alexandrov space since it is a closed convex set.
Note that $|\nabla b_\pm|\equiv 1$, therefore $1$-Lipschitz curve $\alpha$, such that $b_\pm(\alpha(t))=t+\Const$ is a $b_\pm$-gradient curve. 
In this case curve $\alpha(-t)$ is a $b_\mp$-gradient curve.
It follows that for any $t>0$, $\Phi_{b_+}^t\circ\Phi_{b_-}^t=\id_A$.
Set
$$\Phi^t=\left[\begin{matrix}
\Phi_{b_+}^t&\hbox{if}\ t\ge0\\
\Phi_{b_-}^t&\hbox{if}\ t<0
               \end{matrix}\right.$$
Consider map $h:A'\times \RR\to A$ defined by $h:(x,t)\mapsto \Phi^t(x)$.
It is easy to see that $h$ is onto.
Applying lemma \ref{lem:dist-est}(\ref{two-ends}) for $\Phi_{b_+}^t$ and $\Phi_{b_-}^t$, we get that $h$ is a short and non-contracting at the same time; therefore $h$ is an isometry.\qeds

Gradient flow gives a simple proof to the following result which generalizes a Liberman's lemma in \cite{liberman}. 
This generalization was first obtained by Perelman and me in
\cite[5.3]{perelman-petrunin:extremal}, 
a simplified proof was given in
\cite[1.1]{petrunin:extremal}.
See sections~\ref{extremal} and~\ref{QG}  
for the definitions of extremal subset and
quasigeodesic. 

\begin{thm}{\bf Generalized Lieberman's Lemma.}\label{lib-lem} Any unit-speed geodesic for the
induced intrinsic metric on an extremal subset is a quasigeodesic in the ambient
Alexandrov's space.
\end{thm}

\Proof. Let $\gamma\:[a,b]\to E$ be a unit-speed minimizing geodesic in an extremal subset
$E\i A$ and $f$ be a $\lambda$-concave function defined in a neighborhood of
$\gamma$.
Assume $f\circ\gamma$ is not $\lambda$-concave, then there is a non-negative
Lipschitz function $\tau$ with support in $(a,b)$ such that
$$\int\limits_a^b\l[(f\circ\gamma)'\tau'+\lambda\tau\r]\cdot d s< 0$$
Then as follows from lemma~\ref{lem:grad-variation}, for small $t\ge 0$
$$\gamma_t(s)=\Phi^{t\cdot\tau(s)}_f \circ\gamma_0(s)$$
gives a length-contracting homotopy of curves relative to ends and according to
definition~\ref{def:extrim}, it stays in $E$ --- this is a contradiction.\qeds

The fact that gradient flow is stable with respect to collapsing has the
following useful consequence: Let $M_n$ be a collapsing sequence of Riemannian
manifolds with curvature $\ge\kappa$ and $M_n\GHto A$.
For a regular point $p$ let us denote by $F_n(p)$ the \emph{regular fiber}%
\footnote{see footnote~\ref{reg-fib} on page~\pageref{reg-fib}} 
over $p$, 
it is well defined for all large $n$.
Let $f\:A\to\RR$ be a $\lambda$-concave function.
If $\alpha(t)$ is an $f$-gradient curve in $A$  which passes only through
regular points, then for any $t_0<t_1$
there is a homotopy equivalence $F_n(\alpha(t_0))\to F_n(\alpha(t_1))$ with
dilation $\approx e^{\lambda{\cdot}(t_1-t_0)}$.

This observation was used in \cite{KPT} by Kapovitch Tuschmann and me
to prove some properties of almost
nonnegatively curved manifolds. 
In particular, it gave simplified proofs of
the results of Fukaya and Yamaguchi in \cite{FY}):

\begin{thm}{\bf Nilpotency theorem.} Let $M$ be a closed almost nonnegatively curved manifold.
Then a finite cover of $M$ is a \emph{nilpotent space};
that is, 
its fundamental group is nilpotent and acts nilpotently on higher homotopy groups.
\end{thm}

\begin{thm}{\bf Theorem.} Let $M$ be an almost nonnegatively curved
$m$-manifold. Then $\pi_1(M)$ is $\Const(m)$-nilpotent;
that is, 
$\pi_1(M)$ contains a nilpotent subgroup of
index at most $\Const(m)$.
\end{thm}

Gradient flow also gives an alternative proof of the homotopy lifting theorem~\ref{thm:per-ser}. 

Let us start with the definition.
Given a topological space $X$, a map $F\:X\to A$, a finite sequence of $\lambda$-concave functions $\{f_i\}$ on $A$ and continuous functions $\tau_i\:X\to\RR_+$ one can consider a composition of gradient deformations (see \ref{grad-flow})
$$F'(x)=\Phi_{f_N}^{\tau_N(x)}\circ\cdots\circ\Phi_{f_2}^{\tau_2(x)}\circ\Phi_{f_1}^{\tau_1(x)}\circ F(x),\ \ \ F'\:X\to A,$$
which we also call \textit{gradient deformation} of $F$.

Let us define {\it gradient homotopy} to be a gradient deformation of trivial homotopy 
$$F\:[0,1]\times X\to A,\ \ \ F_t(x)=F_0(x)$$
with the functions 
$$\tau_i\:[0,1]\times X\to\RR_+\ \ \ \text{such that}\ \ \ \tau_i(0,x)\equiv 0.$$
If $Y\i X$, then to define \emph{gradient homotopy relative to} $Y$ we assume in addition $$\tau_i(t,y)= 0 \ \ \ \text{for any }\ \ \ y\in Y,\ \ t\in[0,1].$$

Then theorem~\ref{thm:per-ser} follows from lemma~\ref{lem:stable-grad-curves} and the following lemma:

\begin{thm}{\bf Lemma \cite{petrunin:grad-hom}.}\label{lem:hom-approx} Let $A$ be an Alexandrov's space without proper extremal subsets and $K$ be a finite simplicial complex. 
Then, given $\eps>0$, for any  homotopy 
$$F_t\:K\to A,\ \ t\in [0,1]$$ one can construct an $\eps$-close gradient homotopy 
$$G_t\:K\to A$$ such that $G_0\equiv F_0$.
\end{thm}

\section{Gradient exponent}\label{gexp}

\addtocounter{subsection}{1}
\setcounter{thm}{0}

One of the technical difficulties in Alexandrov's geometry comes from
nonextendability of geodesics. 
In particular, the exponential map, $\exp_p\:T_p\to A$, if defined the usual way, can
be undefined in an arbitrary small neighborhood of the origin. 
Here we construct its analog, the \emph{gradient exponential map} $\gexp_p\:T_p\to
A$, which practically solves this problem. 
It has many important properties of the ordinary exponential map, and is even
``better'' in certain respects,
even in the Riemannian universe.

\ 

Let $A$ be an Alexandrov's space and $p\in A$, consider the function $f\z=\dist_p^2/2$.
Recall that $i_{\lam}\:\lam{\cdot} A\to A$ denotes canonical maps (see
page~\pageref{tangent-def2}). 
Consider the one parameter family of maps
$$\Phi^t_{f}\circ i_{e^t}\:e^t{\cdot} A\to A\ \ \ \text{as}\ \ \ t\to\infty \ \ \
\text{so}\ \ \ (e^t{\cdot} A,p)\GHto (T_p,o_p)$$
where $\Phi^t_{f}$ denotes gradient flow (see section~\ref{grad-flow}). 
Let us define the \textit{gradient exponential map} as the limit
$$\gexp_p\:T_p A\to A,\ \ \  \gexp_p=\lim_{t\to\infty}\Phi^t_{f}\circ i_{e^t}.$$

\noi \textit{Existence and uniqueness of gradient exponential.} 
If $A$ is
an Alexandrov's space with curvature $\ge 0$, then $f''\le 1$  and from
lemma~\ref{lem:dist-est}, $\Phi^t_{f}$ is an $e^t$-Lipschitz and therefore
compositions $\Phi^t_{f}\circ i_{e^t}\:e^t{\cdot} A\to A$ are \emph{short}%
\footnote{that is, the maps with Lipschitz constant 1.}. 
Hence a partial limit $\gexp_p\:T_p
A\to A$ exists, and it is a short map.
\footnote{For general lower curvature
bound, $f$ is only $(1+O(r^2))$-concave in the ball $B_r(p)$.
Therefore $\Phi^1_f\:B_{r/e}(p)\to B_{r}(p)$ is $e(1+O(r^2))$-Lipschitz.
By taking compositions of these maps for different $r$ we get that
$\Phi^N_f\:B_{r/e^N}(p)\to B_{r}(p)$ is $e^N(1+O(r^2))$-Lipschitz. 
Obviously, the same is true for any $t\ge 0$;
that is, 
$\Phi^t_f\:B_{r/e^t}(p)\to B_{r}(p)$ 
is $e^t(1+O(r^2))$-Lipschitz, or 
$$\Phi^t_{f}\circ i_{e^t}\:e^t{\cdot} A\to A$$
is $(1+O(r^2))$-Lipschitz on $B_r(p)\i e^t{\cdot} A$. 
This is sufficient for existence of partial limit $\gexp_p\:T_p A\to A$, which
turns out to be $(1+O(r^2))$-Lipschitz on a central ball of radius $r$ in $T_p$.}

Clearly for any partial limit we have
$$\Phi^t_f\circ\gexp_p(v)=\gexp_p(e^t\cdot v)\eqno(*)$$
and since $\Phi^t$ is $e^t$-Lipschitz, it follows that $\gexp_p$ is uniquely
defined.

From above and the definition of extremal subset (\ref{def:extrim}),
we get the following.

\begin{thm}{\bf Property.}
\label{grad-in-extr} If $E\in A$ is an extremal subset, $p\in E$ and
$\xi\in \Sigma_p E$ then $\gexp_p(t\cdot\xi)\in E$ for any $t\ge0$.
\end{thm}

\noi{\bf Radial curves.} 
From identity $(*)$, it follows that for any $\xi\in \Sigma_p$, curve
$$\alpha_\xi\:t\mapsto\gexp_p(t\cdot\xi)$$ 
satisfies the following differential equation
$$\alpha_\xi^+(t)=\frac{|p\,\alpha_\xi(t)|}{t}{\cdot}\nabla_{\alpha_\xi(t)}\dist_p\ \ \
\text{for all}\ \ t>0 \ \ \ \text{and}\ \ \ \alpha_\xi^+(0)=\xi\eqno(\diamond)$$
We will call such a curve {\it radial curve} from $p$ in the direction $\xi$.
From above, such radial curve exists and is unique in any direction.

Clearly, for any radial curve from $p$, $|p\alpha_\xi(t)|\le t$; and if this inequality is
exact for some $t_0$ then $\alpha_\xi\:[0,t_0]\to A$ is a unit-speed minimizing
geodesic starting at $p$ in the direction $\xi\in \Sigma_p$.
In other words,
$$\gexp_p\circ\log_p=\id_A.%
\footnote{In proposition~\ref{prop:unique-gexp-inverse} we will show that
$\alpha_\xi((0,t_0))$ does not meet any other radial curve from $p$.}$$

Next lemma gives a comparison inequality for radial curves.

\begin{thm}{\bf Lemma.} \label{lem:monotonic}
Let $A\in \Alex$,  $f\:A\to\RR$ be a $\lambda$-concave function $\lambda\ge 0$ then for any $p\in A$ and $\xi\in \Sigma_p$
$$f\circ\gexp_p(t\cdot\xi)\le f(p)+t\cdot\d_p f(\xi)+t^2{\cdot}\tfrac\lambda2.$$
Moreover, the function
$$\vartheta(t)=\{f\circ\gexp_p(t\cdot\xi)-f(p)-t^2\cdot\tfrac\lambda2\}/t$$
is non-increasing.
\end{thm}

In particular, applying this lemma for $f=\dist_q^2/2$ we get

\begin{thm}{\bf Corollary.} \label{cor:angle--}
If $A\in\Alex(0)$ then for any $p,q,\in A$ and
$\xi\in \Sigma_p$,
$$\tilde\angle_0(t,|\gexp_p(t{\cdot}\xi)q|,|p q|)$$ 
is non-increasing in $t$.%
\footnote{$\tilde\angle_\kappa(a,b,c)$ denotes angle
opposite to $b$ in a triangle with sides $a,b,c$ in $\Lob_\kappa$.} In
particular,
$$\tilde\angle_{0}(t,|\gexp_p(t{\cdot}\xi)\,q|,|p q|) \le\angle(\xi,\uparrow_p^q).$$
\end{thm}

In \ref{sph-hyp-exp} you will find a version of this corollary for arbitrary lower curvature bound.

\Proof\ {\it of lemma~\ref{lem:monotonic}.} 
Recall that $\nabla_q\dist_p$ is polar to the set $\Uparrow_q^p\i T_q$ (see
example~(\ref{polar}) on page~\pageref{polar}). 
In particular, from inequality $(**)$ on page~\pageref{**-polar-inq},
$$\d_q f(\nabla_q\dist_p) 
+\inf_{\zeta\in\Uparrow_q^p}\{\d_q f(\zeta)\}
\le 0$$

On the other hand, since $f''\le \lambda$,
$$\d_q f(\zeta)\ge \frac{f(p)-f(q)-\lambda{\cdot}
|p q|^2/2}{|p q|}\ \ \text{for any}\
\ \zeta\in \Uparrow_q^p,$$
therefore
$$\d_q f(\nabla_q\dist_p)
\le
\frac{f(q)-f(p)+\tfrac\lambda2{\cdot}
|p q|^2}{|p q|}.$$

Set $\alpha_\xi(t)=\gexp(t\cdot\xi)$, $q=\alpha_\xi(t_0)$, then 
$\alpha^+_\xi(t_0)=\tfrac{|p q|}{t}{\cdot}\nabla_{q}\dist_p$ as in $(\diamond)$. 
Therefore,
\begin{align*}
(f\circ\alpha_\xi)^+(t_0)
&=
\d_q f(\alpha^+_\xi(t_0))\le
\\
&\le
\frac{|p q|}{t_0}
\cdot 
\l[\frac{f(q)-f(p)+\tfrac\lambda2{\cdot}
|p q|^2}{|p q|}\r]
=
\\
&=
\frac{f(q)-f(p)
+\tfrac\lambda2{\cdot}
|p q|^2}{t_0}
\le
\intertext{since $|p q|\le t_0$ and $\lambda\ge 0$,}
&\le \frac{f(q)-f(p)+\tfrac\lambda2{\cdot}t^2_0}{t_0}
=
\\
&=
\frac{f(\alpha_\xi(t_0))-f(p)+\tfrac\lambda2{\cdot}t^2_0}{t_0}.
\end{align*}
Substituting this inequality in the expression for derivative of $\vartheta$,
$$\vartheta^+(t_0)=\frac{(f\circ\alpha_\xi)^+(t)}{t_0}
-\frac{f\circ\gexp_p(t_0\cdot\xi)- f(p)}{t_0^2}-\tfrac\lambda2,$$ we get
$\vartheta^+\le 0$; 
that is, $\vartheta$ is non-increasing. 

Clearly, $\vartheta(0)=\d_p f(\xi)$ and so the first statement follows.\qeds

\subsection{Spherical and hyperbolic gradient exponents}\label{sph-hyp-exp}

The gradient exponent described above is sufficient for most applications. 
It works perfectly for non-negatively curved Alexandrov's spaces and 
where one does not care for the actual lower curvature bound.
However, for fine analysis on spaces with curvature $\ge \kappa$, there is a better
analog of this map, which we denote
$\gexp_p(\kappa;v)$; $\gexp_p(0;v)=\gexp_p(v)$.

In addition to case $\kappa=0$, it is enough to consider only two cases:
$\kappa=\pm1$, the rest can be obtained by rescalings. 
We will define two maps: $\gexp_p(-1,*)$ and  $\gexp_p(1,*)$, and list their
properties, leaving calculations to the reader. 
These properties are analogous to the following properties of the ordinary gradient
exponent:

\begin{enumerate}[$\diamond$]
\item if $A\in\Alex({0})$, then $\gexp_p\:T_p\to A$ is distance non-increasing. 

Moreover, for any $q\in A$, the angle $$\tilde\angle_{0}(t,|\gexp_p(t\cdot\xi)\,q|,|p q|)$$ is
non-increasing in $t$ (see corollary~\ref{cor:angle--}). In particular
$$\tilde\angle_{0}(t,|\gexp_p(t\cdot\xi)\,q|,|p q|) \le\angle(\xi,\uparrow_p^q).$$
\end{enumerate}

\begin{thm}
{\bf Case $\kappa=-1$.}
\end{thm}
The hyperbolic radial curves are defined by the following differential equation
$$\alpha^+_\xi(t)=
\frac{\tanh|p\alpha_\xi(t)|}{\tanh t}\cdot\nabla_{\alpha_\xi(t)}\dist_p\ \ \
\text{and}\ \ \ \alpha^+_\xi(0)=\xi.$$
These radial curves are defined for all $t\in [0,\infty)$.
Let us define
$$\gexp_p(-1;t\cdot \xi)=\alpha_\xi(t).$$
This map is defined on tangent cone $T_p$.
Let us equip the tangent cone with a hyperbolic metric $\mathfrak h(u,v)$ defined by
the hyperbolic rule of cosines
$$\cosh(\mathfrak h(u,v)) 
=\cosh|u|\cdot\cosh|v|-\sinh|u|\cdot\sinh|v|\cdot\cos\alpha,$$
where $u,v\in T_p$ and $\alpha=\angle u o_p v$.
$(T_p,\mathfrak h)\in\Alex(-1)$, this is a so called \emph{elliptic cone} over $\Sigma_p$; 
see \cite[4.3.2]{BGP} and \cite{alexander-bishop:worps}.
Here are the main properties of $\gexp(-1;*)$:
\begin{enumerate}[$\diamond$]
\item if $A\in\Alex(-1)$, then $\gexp(-1;*)\:(T_p,\mathfrak h)\to A$ is distance
non-increasing. 

Moreover, the function 
$$t\mapsto \tilde\angle_{-1}(t,|\gexp(-1;t\cdot\xi )\,q|,|p q|)$$ 
is non-increasing in $t$. 
In particular for any $t>0$, $$\tilde\angle_{-1}(t,|\gexp(-1;t\cdot\xi )\,q|,|p
q|)\le \angle(\xi,\uparrow_p^q).$$
\end{enumerate}

\begin{thm}\label{kappa=1}
{\bf Case $\kappa=1$.}
\end{thm}
For unit tangent vector $\xi\in \Sigma_p$, the spherical radial curve is defined to satisfy the
following identity:
$$\alpha^+_\xi(t)=\frac{\tan|p\alpha_\xi(t)|}{\tan t}\cdot\nabla_{\alpha_\xi(t)}\dist_p\ \ \ \text{and}\ \ \
\alpha^+_\xi(0)=\xi.$$
These radial curves are defined for all $t\in [0,\tfrac\pi2]$.
Let us define the spherical gradient exponential map by
$$\gexp_p(1;t\cdot \xi)=\alpha_\xi(t).$$
This map is well defined on $\bar B_{\pi/2}(o_p)\i T_p$.
Let us equip $\bar B_{\pi/2}(o_p)$ with a spherical distance $\mathfrak{s}(u,v)$ defined by
the spherical rule of cosines
$$\cos(\mathfrak{s}(u,v)) 
=
\cos|u||\cdot\cos|v|+\sin|u||\cdot\sin|v||\cdot\cos\alpha,$$
where $u,v\in B_\pi(o_p)\i T_p$ and $\alpha=\angle u o_p v$.
$(\bar B_{\pi}(o_p),\mathfrak{s})\in\Alex(1)$, this is isometric to \emph{spherical suspension} $\Sigma(\Sigma_p)$, see \cite[4.3.1]{BGP} and \cite{alexander-bishop:worps}.
Here are the main properties of $\gexp(1;*)$:
\begin{enumerate}[$\diamond$]
\item If $A\in\Alex(1)$ then $\gexp_p(1,*)\:(\bar B_{\pi/2}(o_p),\mathfrak{s})\to
A$ is distance non-increasing.

Moreover, if $|p q|\le\tfrac\pi2$, then function
$$t\mapsto\tilde\angle_{1}(t,|\gexp_p(1;t\cdot \xi)\,q|,|p q|)$$ 
is non-increasing in $t$. 
In particular, for any $t>0$
$$\tilde\angle_{1}(t,|\gexp_p(1;t\cdot \xi)\,q|,|p
q|)\le\angle(\xi,\uparrow_p^q).$$
\end{enumerate}

\subsection{Applications}

One of the main applications of gradient exponent and radial curves is the proof
of existence of quasigeodesics; 
see property~\ref{exist-qg} page~\pageref{exist-qg} and appendix~\ref{constr-qg} for the proof.

An infinite-dimensional generalization of gradient exponent was introduced by Perelman to make the last step in the proof of equality of Hausdorff and topological dimension for Alexandrov's spaces, 
see \cite[A.4]{perelman-petrunin:qg}.
As it was shown by Plaut (see \cite{plaut:dimension} or
\cite[151]{plaut:survey}), 
if $\dim_H A\ge m$, then there is a point $p\in A$, the tangent cone of which contains a subcone $W\i T_p$ isometric to Euclidean $m$-space.
Then infinite-dimensional analogs of properties  in section~\ref{sph-hyp-exp} ensure that image $\gexp_p (W)$ has topological dimension $\ge m$ and therefore $\dim A\ge m$.

The following statement has been proven by Perelman in \cite{perelman:spaces2}, 
then its
formulation was made more exact by Alexander and Bishop in \cite{alexander-bishop:fk}. 
Here we give a simplified proof with the use of a gradient exponent.

\begin{thm}{\bf Theorem.} \label{thm:dist-to-bry} 
Let $A\in \Alex(\kappa)$ and
$\partial A\not=\emptyset$; 
then the function $f\z=\sigma_\kappa\circ\dist_{\partial
A}$%
\footnote{\label{sigma_k}$\sigma_\kappa\:\RR\to\RR$ is defined by
$$\sigma_\kappa(x)
=
\sum_{n=0}^\infty\frac{(-\kappa)^{n}}{ (2n+1)!}{\cdot}x^{2n+1}=
\left[\begin{matrix}
{ \frac{1}{\sqrt \kappa}{\cdot}\sin({x{\cdot}\sqrt \kappa})}&{\hbox{if}\ \kappa>0}\\
            {x}&{\hbox{if}\ \kappa=0}\\
            { \frac{1}{\sqrt{-\kappa}}{\cdot}\sinh({x{\cdot}\sqrt {-\kappa}})}&{\hbox{if}\
\kappa<0}\\
\end{matrix}\right. .$$} is $(-\kappa{\cdot}f)$-concave in $\Omega=A\backslash\partial A$.%
\footnote{Note that by definition~\ref{def:with-bry}, 
$f$ is not semiconcave in $A$.} 

In particular,
\begin{enumerate}[(i)]
\item if $\kappa=0$, $\dist_{\partial A}$ is concave in $\Omega$;
\item if $\kappa>0$, the level sets $L_x=\dist^{-1}_{\partial A}(x)\i A$, $x>0$
are strictly concave hypersurfaces.
\end{enumerate}
\end{thm}

\begin{wrapfigure}[14]{r}{37mm}
\begin{lpic}[t(-15mm),b(0mm),r(0mm),l(0mm)]{pics/dist-to-bry(0.4)}
\lbl[b]{14,144;$\tilde\gamma(0)$}
\lbl[b]{63,158;$\tilde\gamma(\tau)$}
\lbl[lt]{22,134;$\alpha$}
\lbl[b]{18,39;$\tilde\beta$}
\lbl[t]{15,15;$\tilde p$}
\lbl[t]{63,15;$\tilde q$}
\lbl[lb]{75,17;$\partial\Lob^+_\kappa$}
\end{lpic}
\end{wrapfigure}

\Proof. 
We have to show that for any unit-speed geodesic $\gamma$, the function 
$f\circ\gamma$ is $(-\kappa{\cdot}f\circ\gamma)$-concave; 
that is, for any $t_0$,
$$(f\circ\gamma)''(t_0)\le -\kappa{\cdot}f\circ\gamma(t_0)$$
\emph{in a barrier sense}%
\footnote{For a continuous function $f$, $f''(t_0)\le
c$ \emph{in a barrier sense} means that there is a smooth function $\bar f$ such
that $f\le \bar f$, $f(t_0)=\bar f(t_0)$ and $\bar f''(t_0)\le c$}. 
Without loss
of generality we can assume $t_0=0$.

Direct calculations show that the statement is true for \label{Lob_k^+} $A=\Lob_\kappa^+$, the
halfspace  of the model space $\Lob_\kappa$.

Let $p\in \partial A$ be a closest point to $\gamma(0)$ and
$\alpha=\angle(\gamma^+(0),\uparrow_{\gamma(0)}^p)$.

Consider the following configuration in the model halfspace $\Lob_\kappa^+$: Take a
point $\tilde p\in\partial\Lob_\kappa^+$ and consider the geodesic $\tilde\gamma$ in
$\Lob_\kappa^+$ such that 
$$|\gamma(0) p|
=|\tilde\gamma(0)\tilde
p|
=|\tilde\gamma(0)\,\partial\Lob^+_\kappa|,$$ 
so $\tilde p$ is the closest point to $\tilde\gamma(0)$ on the
boundary%
\footnote{in case $\kappa>0$ it is possible only if 
$|\gamma(0) p|\le
\frac{\pi}{2{\cdot}\sqrt\kappa}$, but this is always the case since otherwise any small
variation of $p$ in $\partial A$ decreases distance $|\gamma(0) p|$.} and
$$\angle(\tilde\gamma^+(0),\uparrow_{\tilde\gamma(0)}^{\tilde p})=\alpha.$$
Then it is enough to show that 
$$\dist_{\partial A} \gamma(\tau) \le
\dist_{\partial\Lobs_\kappa^+}\tilde\gamma(\tau)+o(\tau^2).$$
Set 
$$\beta(\tau)=\angle \gamma(0)\, p\, \gamma(\tau)$$ 
and
$$\tilde\beta(\tau)=\angle \tilde\gamma(0)\,\tilde p\, \tilde\gamma(\tau).$$
From the comparison inequalities
$$|p\gamma(\tau)|\le|\tilde p\tilde \gamma(\tau)|$$
and
$$\theta(\tau)=\max\l\{0,\,\tilde\beta(\tau)-\beta(\tau)\r\}=o(\tau).\eqno(*)$$
Note that the tangent cone at $p$ splits: 
$T_p A=\RR_+\times T_p\partial A$.%
\footnote{This follows from the fact that $p$ lies on a shortest path between
two inverse images of $\gamma(0)$ in the doubling $\tilde A$ of $A$, see
\cite[7.15]{BGP}.}
Therefore we can represent $v=\log_p\gamma(\tau)\in T_p A$ as $v=(s,w)\in
\RR_+\times T_p\partial A$.
Let $\tilde q=\tilde q(\tau)\in\partial\Lob_\kappa$ be the closest point to
$\tilde\gamma(\tau)$, so 
\begin{align*}
\angle(\uparrow_p^{\gamma(\tau)},w)&=
\tfrac\pi2-\beta(\tau)\le 
\\
&\le \tfrac\pi2-\tilde\beta(\tau) -\theta(\tau)=
\\
&=
\angle\tilde\gamma(\tau)\tilde p\tilde q+o(\tau).
\end{align*}
Set $q=\gexp_p\l(\kappa;|\tilde p\tilde
q|\frac{w}{|w|}\r)$.%
\footnote{\label{qg-grad} 
Alternatively, one can set $q=\gamma(|\tilde p\tilde q|)$, where $\gamma$ is a
quasigeodesic in $\partial A$ starting at $p$ in direction $\frac{w}{|w|}\in
\Sigma_p$ (it exists by second part of property~\ref{exist-qg} on
page~\pageref{exist-qg}).} 
Since gradient curves preserve extremal subsets $q\in \partial A$ (see
property~\ref{grad-in-extr} on page~\pageref{grad-in-extr}).
Clearly $|\tilde p\tilde q|=O(\tau)$, therefore applying the comparison from
section~\ref{sph-hyp-exp} (or Corollary~\ref{cor:angle--} if $\kappa=0$)
together with $(*)$, we get
\begin{align*}
\dist_{\partial A} \gamma(\tau) 
&\le|q \gamma(\tau)|
\le
\\
&\le|\tilde q\tilde\gamma(\tau)|+O\l(|\tilde p\tilde q|\cdot\theta(\tau)\r)
=
\\
&=\dist_{\partial\Lobs_\kappa^+}\tilde\gamma(\tau)+o(\tau^2).
\end{align*}
\qedsf

The following corollary implies that the Lipschitz condition in the definition of
convex function~\ref{def:with-bry}-- \ref{def:with-no-bry} can be relaxed to usual continuity.

\begin{thm}{\bf Corollary.}\label{cor:cont-conv} Let $A\in\Alex$, $\partial A=\emptyset$,
$\lambda\in\RR$ and $\Omega\i A$ be open.

Assume $f\:\Omega\to\RR$ is a continuous function such that for any unit-speed geodesic $\gamma$ in $\Omega$ we have that the function
$$t\mapsto f\circ\gamma-\tfrac\lambda2{\cdot}t^2$$
is concave; then $f$ is locally Lipschitz.
 
In particular, $f''\le \lambda$ in the sense of
definition~\ref{def:with-bry}.
\end{thm}

\Proof. Assume $f$ is not Lipschitz at $p\in \Omega$.
Without loss of generality we can assume that $\Omega$ is
convex%
\footnote{Otherwise, pass to a small convex neighborhood of $p$ which
exists by by corollary~\ref{cor:convex-nbhd}.} 
and $\lambda<0$%
\footnote{Otherwise, add a very concave (Lipschitz) function which
exists by theorem~\ref{thm:strictly-concave}}.
Then, since $f$ is continuous, sub-graph
$$X_f=\{(x,y)\in \bar\Omega\times\RR|y\le f(x)\}$$
is closed convex subset of $A\times\RR$, therefore it forms an Alexandrov's
space. 

Since $f$ is not Lipschitz at $p$, there is a sequence of pairs of points
$(p_n,q_n)$ in $A$, such that 
$$p_n,q_n\to p\ \ \ \text{and}\ \ \  \frac{f(p_n)-f(q_n)}{|p_n
q_n|}\to+\infty.$$ 
Consider a sequence of radial curves $\alpha_n$ in $X_f$ which extend  shortest
paths from $(p_n,f(p_n))$ to $(q_n,f(q_n))$.
Since the boundary $\partial X_f\i X_f$ is an extremal subset, we have $\alpha_n(t)\in
\partial X_f$ for all 
\begin{align*}
t
&\ge \ell_n
=
\\
&=|(p_n,f(p_n))(q_n,f(q_n))|
=
\\
&=\sqrt{|p_n
q_n|^2+(f(p_n)-f(q_n))^2}.
\end{align*}
Clearly, the function $h\:X_f\to\RR$, $h\:(x,y)\mapsto y$ is concave.
Therefore, from \ref{lem:monotonic}, there is a sequence $t_n>\ell_n$, so
$\alpha_n(t_n)\to (p,f(p)-1)$.
Therefore, $(p,f(p)-1)\in\partial X_f$ thus $p\in\partial A$;
that is, $\partial A\not=\emptyset$, a
contradiction.
\qeds

\begin{thm}{\bf Corollary.} \label{cor:eq-qg}
Let $A\in\Alex^m(\kappa)$, $m\ge2$ and $\gamma$ be a unit-speed curve in $A$
which has a convex $\kappa$-developing with respect to any point.
Then $\gamma$ is a quasigeodesic;
that is, for any $\lambda$-concave function $f$,
function $f\circ\gamma$ is $\lambda$-concave.
\end{thm}

\Proof.  Let us first note that in the proof of theorem~\ref{thm:dist-to-bry} we
used only two properties of curve $\gamma$:
$|\gamma^\pm|=1$ and the convexity of the $\kappa$-development of $\gamma$ with respect
to $p$.

Assume $\kappa=\lambda=0$ then sub-graph of $f$
$$X_f=\{(x,y)\in A\times\RR\ |\ y\le f(x)\}$$ 
is a closed convex subset, therefore it forms an Alexandrov's space.

Applying the above remark, we get that if $\gamma$ is a unit-speed curve in
$X_f\backslash \partial X_f$ with convex $0$-developing with respect to any
point then $\dist_{\partial X_f}\circ\gamma$ is concave.
 Hence, for any $\eps>0$, the function $f_\eps$, which has the level set $\dist_{\partial X_f}^{-1}(\eps)\i \RR\times A$ like the graph, has a concave restriction to any curve $\gamma$ in
$A$ with a convex $0$-developing with respect to any
point in $A\backslash\gamma$. 
Clearly, $f_\eps\to f$ as $\eps\to 0$, hence $f\circ\gamma$ is concave. 

For $\lambda$-concave function the set $X_f$ is no longer convex, but it becomes convex if one changes metric on $A\times\RR$ to \emph{parabolic cone}%
\footnote{\label{par-cone} that is, a warped-product
$\RR\times_{\exp(\Const{\cdot}t)} A$, which is an Alexandrov's space, see \cite[4.3.3]{BGP} and \cite{alexander-bishop:worps}} 
and then one can repeat the same arguments.\qeds

\noi{\bf Remark} One can also get this corollary from the following lemma:

\begin{thm}{\bf Lemma.}\label{f_eps}
Let $A\in\Alex^m(\kappa)$, $\Omega$ be an open subset of $A$ and
 $f\:\Omega\to\RR$ be a $\lambda$-concave $L$-Lipschitz function.
Then function
$$f_\eps(y)=\min_{x\in\Omega}\{f(x)+\tfrac1\eps{\cdot}|x y|^2\}$$
is $(\lambda+\delta)$-concave in the domain of definition%
\footnote{that is the set where the minimum is defined.} 
for some%
\footnote{this function $\delta(L,\lambda,\kappa,\eps)$ is achieved for the model space $\Lambda_\kappa$} 
$\delta=\delta(L,\lambda,\kappa,\eps)$, $\delta\to0$ as $\eps\to 0$.

Moreover, if $m\ge 2$ and $\gamma$ is a unit-speed curve in $A$ with $\kappa$-convex developing with respect to any point then 
$f_\eps\circ\gamma$ is also $(\lambda+\delta)$-concave.
\end{thm}

\Proof. 
It is analogous to theorem~\ref{thm:dist-to-bry}. 
We only indicate it in the simplest case, $\kappa=\lambda=0$. 
In this case $\delta$ can be taken to be $0$.

Let $\gamma$ be a unit-speed geodesic (or it satisfies the last condition in the lemma). 
It is enough to show that for any $t_0$
$$(f_\eps\circ\gamma)''(t_0)\le0$$
in a barrier sense. 

Let $y=\gamma(t_0)$ and $x\in \Omega$ be a point for which 
$f_\eps(y)=f(x)+\tfrac1\eps{\cdot}|x y|^2$.
The tangent cone $T_x$ splits in direction $\uparrow_y^x$;
that is, there is an isometry $T_x\to \RR\times \Cone$ such that $\uparrow_x^y\mapsto(1,o)$, 
where $o\in\Cone$ is its origin.
Let 
$$\log_x\gamma(t)=(a(t),v(t))\in\RR\times \Cone=T_x.$$
Consider vector 
$$w(t)=(a(t)-|x y|,v(t))\in \RR\times \Cone=T_x.$$ 
Clearly $|w(t)|\ge|x\gamma(t)|$.
Set $x(t)=\gexp_y(w(t))$ then lemma~\ref{lem:monotonic} gives an estimate for $f\circ x(t))$ while corollary~\ref{cor:angle--} gives an estimate for $|\gamma(t)x(t)|^2$. Hence the result.
\qeds

Here is yet another illustration for the use of gradient exponents.
At first sight it seems very simple, but the proof is not quite obvious.
In fact, I did not find any proof of this without applying the gradient exponent. 

\begin{thm} {\bf Lytchak's problem.}\label{lyt-prob} Let $A\in\Alex^m(1)$. 
Show
that
$$\vol_{m-1}\partial A\le \vol_{m-1}S^{m-1}$$
where $\partial A$ denotes the boundary of $A$ and $S^{m-1}$ the unit $(m-1)$-sphere.
\end{thm}

The problem would have followed from conjecture~\ref{conj:bry} (that boundary of
an Alexandrov's space is an Alexandrov's space), but before this conjecture has been proven, any partial
result is of some interest.
Among other corollaries of
conjecture~\ref{conj:bry}, it is expected that if $A\in \Alex(1)$ then $\partial
A$, equipped with induced intrinsic metric, admits a noncontracting map to
$S^{m-1}$. 
In particular, its intrinsic diameter is at most $\pi$, and perimeter
of any triangle in $\partial A$ is at most $2\pi$. 
This does not follow from the proof below, 
since in general 
$\gexp_z(1;\partial B_{\pi/2}(o_z))\not\i\partial A$;
that is, the image $\gexp_z(1;\partial B_{\pi/2}(o_z))$
might have creases inside of $A$, 
which could be used as a shortcut for
curves with ends in $\partial A$.

Let us first prepare a proposition:

\begin{thm}{\bf Proposition.}\label{prop:unique-gexp-inverse}
The inverse of the gradient exponential map $\gexp^{-1}_p(\kappa;*)$ is uniquely
defined inside any minimizing geodesic starting at $p$.
\end{thm}

\Proof. Let $\gamma\:[0,t_0]\to A$ be a unit-speed minimizing geodesic,
$\gamma(0)=p$, $\gamma(t_0)=q$.
 From the angle comparison we get that $|\nabla_x\dist_p|\ge-\cos\tilde\angle_\kappa p
x q$. Therefore, for any $\zeta$ we have
$$|p\alpha_\zeta(t)|^+_t
\ge 
-|\alpha^+_\zeta(t)|{\cdot}\cos\tilde\angle_\kappa
p\,\alpha_\zeta(t)\,q\ \ \text{and}\ \ |\alpha_\zeta(t)q|^+_t\ge-|\alpha^+_\zeta(t)|.$$
Therefore, $\tilde\angle_\kappa p\, q\, \alpha_\zeta(t)$ is nondecreasing in $t$, hence the result.
\qeds

\Proof\ \textit{of \ref{lyt-prob}}. 
Let $z\in A$ be the point at maximal distance from $\partial A$, in particular
it realizes maximum of $f=\sigma_1\circ\dist_{\partial
A}=\sin\circ\dist_{\partial A}$.
From theorem~\ref{thm:dist-to-bry}, $f''+f\le 0$ and $f(z)\le 1$.

Note that $A\i \bar B_{\pi/2}(z)$, otherwise if $y\in A$ with $|y z|>\tfrac\pi2$, then
since $f''+f\le 0$ and $f(y)\ge 0$, 
we have $\d f(\uparrow_z^y)>0$; 
that is, $z$ is not a maximum of $f$.

From this it follows that gradient exponent
$$\gexp_z(1;*)\:(\bar B_{\pi/2}(o_z),\mathfrak s)\to A$$
is a short onto map. 

Moreover,
$$\partial A\i\gexp_z(\partial B_{\pi/2}(o_z)).$$ 
Indeed, $\gexp$ gives a homotopy equivalence $\partial B_{\pi/2}(o_z)\to
A\backslash \{z\}$. 
Clearly, $\Sigma_z\z=\partial (B_{\pi/2}(o_z),\mathfrak s)$ has no boundary, therefore 
$H_{m-1}(\partial A,\ZZ_2)\not=0$, 
see Lemma 1 in \cite{grove-petersen:rad-sphere} by Grove and Petersen. 
Hence for any point $x\in\partial A$, any minimizing geodesic $z x$ must have
a point of the image $\gexp(1;\partial B_{\pi/2}(o))$ but, as it is shown in
proposition~\ref{prop:unique-gexp-inverse}, it can only be its end $x$. 

Now since $$\gexp_z(1;*)\:(\bar B_{\pi/2}(o_z),\mathfrak s)\to A$$ is short and
$(\partial B_{\pi/2}(o),\mathfrak s)$ is isometric to $\Sigma_z A$ we get
$\vol\partial A\le\vol \Sigma_z A$ and clearly, $\vol \Sigma_z A\le \vol S^{m-1}$.\qeds

\section{Extremal subsets}\label{extremal}

Imagine that you want to move a heavy box inside an empty room by pushing it around. 
If the box is located in the middle of the room, you  can push it in any direction. 
But once it is pushed against a wall you can not push it back to the center;
and once it is pushed into a corner you cannot push it anywhere anymore. 
The same is true if one tries to move a point in an Alexandrov's space by pushing
it along a gradient flow, but the role of walls and corners is played by  extremal
subsets.

Extremal subsets first appeared in the study of their special case --- the boundary of
an Alexandrov's space; 
introduced by Perelman and me in \cite{perelman-petrunin:extremal}, 
and
were studied further in \cite{petrunin:extremal}
and \cite{perelman:collapsing}.

An Alexandrov's space without extremal subsets resembles a very non-smooth Riemannian manifold.
The presence of extremal subsets makes it behave as something new and maybe
interesting; it gives an additional combinatorial structure which
reflects geometry and topology of the space itself, as well as of nearby spaces.

\subsection{Definition and properties.}

We define extremal subsets as ``ideals'' of the gradient flow. 

\begin{thm}{\bf Definition.}\label{def:extrim} Let $A\in \Alex$.

$E\i A$ is an \emph{extremal subset}, if for any semiconcave function
$f$ on $A$, $t\ge 0$ and $x\in E$, we have $\Phi_f^t(x)\in E$.
\end{thm}

Recall that $\Phi^t_f$ denotes the $f$-gradient flow for time $t$, see
\ref{grad-flow}.
Here is a quick corollary of this definition:
\newcounter{extr-prop}
\begin{enumerate}
\item Extremal subsets are closed. Moreover:
\begin{enumerate}[(i)]
\item For any point $p\in A$, there is an $\eps>0$, such that if an extremal subset
intersects $\eps$-neighborhood of $p$ then it contains $p$.

\item On each extremal subset the intrinsic metric is locally finite.
\end{enumerate}
\setcounter{extr-prop}{\value{enumi}}
These properties follow from the fact that the gradient flow for a $\lambda$-concave
function with $d_p f|_{\Sigma_p}<0$ pushes a small ball $B_\eps(p)$ to $p$ in time
proportionate to~$\eps$.
\end{enumerate}

\noi\textit{Examples.}
\begin{enumerate}[(i)]
\item An Alexandrov's space itself, as well as the empty set, forms an extremal subsets.
\item A point $p\in A$ forms a one-point extremal subset if its
space of directions $\Sigma_p$ has a diameter $\le\tfrac\pi2$
\item \label{ex:t-cone}
If one takes a subset $X$ of points of an Alexandrov's space with tangent cones
homeomorphic to each other then its closure of $X$ forms an extremal subset.
The same holds for the closure of a connected component of of $X$.

Equivalently, the set $X$ can be described as a set of points with homeomorphic small spherical neighborhoods. The equivalence follows from Perelman's stability theorem.

In particular, if in this construction we take points with tangent cone homeomorphic
to $\RR_{+}\times \RR^{m-1}$ then we get the boundary of an Alexandrov's space.

This follows from theorem~\ref{thm:dist-extr} and the Morse lemma (property~\ref{morse} page~\pageref{morse}). 

\item  Let $G$ be a closed subgroups in the group of isometries of an Alexandrov space $A$.
Denote by $A^G$ the fixed point set of $G$.
Then the projection of $A^G$ in $A/G$ forms an extremal subset.

\item If $\iota\:A\to A'$ is a submetry and $E\subset A$ is an extremal subset then $\iota(E)$ is an extremal subset in $A'$.
\end{enumerate}

\noi The following theorem gives an equivalence of our definition of extremal subset and the definition given in
\cite{perelman-petrunin:extremal}:

\begin{thm}{\bf Theorem.}\label{thm:dist-extr}  A closed subset $E$ in
an Alexandrov's space $A$ is extremal if and only if for any $q\in A\backslash E$,
the following condition is fulfilled:

If $\dist_q$ has a local minimum on $E$ at a point $p$, then $p$ is a critical
point of $\dist_q$ on $A$;
that is, $\nabla_p\dist_q=o_p$.
\end{thm}

\Proof. For the ``only if'' part, note that if $p\in E$ is not a critical point of
$\dist_q$, then one can find a point $x$ close to $p$ so that $\uparrow_p^x$ is
uniquely defined and close to the direction of $\nabla_p\dist_q$,
so $d_p\dist_q(\uparrow_p^x)>0$. 
Since $\nabla_p\dist_x$ is polar to $\uparrow_p^x$ (see
page~\pageref{supp-polar}) we get 
$$(\d_p\dist_q)(\nabla_p\dist_x)<0,$$ 
see inequality~\ref{*-polar-inq} on page~\pageref{*-polar-inq}.
Hence, the gradient flow $\Phi_{\dist_x}^t$ pushes the point $p$ closer to $q$, which
contradicts the fact that $p$ is a minimum point $\dist_q$ on $E$.

To prove the ``if'' part, it is enough to show that if $F\i A$ satisfies the condition of the theorem,
then for any $p\in F$, and any semiconcave function $f$, either $\nabla_p f=o_p$ or 
$\tfrac{\nabla_p f}{|\nabla_p f|}\in \Sigma_p F$.
If so, an $f$-gradient curve
can be obtained as a limit of broken lines with vertexes on $F$, and from
uniqueness, any gradient curve which starts at $F$ lives in $F$.

Let us use induction on $\dim A$. 
Note that if $F\i A$ satisfies the condition,
then the same is true for $\Sigma_p F\i\Sigma_p$, for any $p\in F$.
Then using the inductive hypothesis we get that $\Sigma_p F\i \Sigma_p$ is an extremal subset.

If $p$ is isolated, then clearly $\diam \Sigma_p\le \tfrac\pi2$ and therefore $\nabla_p
f=o$, so we can assume $\Sigma_p F\not=\emptyset$. 

Note that $(\d_p f)''+\d_p f\le 0$ on $\Sigma_p$ 
(see~\ref{variation-CF}, page~\pageref{variation-CF}).
Take $\xi=\tfrac{\nabla_p f}{|\nabla_p f|}$, so 
$\xi\in \Sigma_p$ is the maximal point of $\d_p f$.
Let $\eta\in \Sigma_p F$ be a direction closest  to $\xi$, then
$\angle(\xi,\eta)\le \tfrac\pi2$; otherwise $F$ would not satisfy the condition in the theorem for a point $q$ with $\uparrow_p^q\,\approx\xi$.
Hence, since $\Sigma_p F\i \Sigma_p$ is an extremal subset, $\nabla_\eta(\d_p f)\in \Sigma_\eta \Sigma_p F$ and
therefore 
$$(\d_\eta\d_p f)(\uparrow_\eta^\xi)\le\<\nabla_\eta\d_p
f,\,\uparrow_\eta^\xi\>\le0.$$
Hence, $\d_p f(\eta)\ge\d_p f(\xi)$, and therefore $\xi=\eta$;
that is,  $\tfrac{\nabla_p f}{|\nabla_p f|}\in \Sigma_p F$. 
\qeds

From this theorem it follows that in the definition of extremal subset (\ref{def:extrim}),
one has to check only squares of distance functions. 
Namely:
\textit{Let $A\in \Alex$, then $E\i A$ is an extremal subset, if for any point $p\in A$, and any $x\in E$, we
have $\Phi_{\dist_p^2}^t(x)\in E$ for any $t\ge 0$}.

In particular, applying lemma~\ref{lem:stable-grad-curves} we get

\begin{thm}{\bf Lemma.}
\label{lem:limit-extr} The limit of extremal subsets is an extremal subset. 

Namely, if $A_n\in\Alex^m(\kappa)$, $A_n\GHto A$  and $E_n\i A_n$ is a sequence
of extremal subsets such that $E_n\to E\i A$ then $E$ is an extremal subset of
$A$. 
\end{thm}

The following is yet another important technical lemma:

\begin{thm} {\bf Lemma.} {\rm \cite[3.1(2)]{perelman-petrunin:extremal}}
\label{lem:dist-to-extr}
Let $A\in\Alex$ be compact, then there is $\eps>0$
such that $\dist_E$ has no critical values in $(0,\eps)$. 
Moreover,
$$|\nabla_x\dist_E|>\eps\ \ \text{if}\ \ 0<\dist_E(x)<\eps.$$

For a non-compact $A$,  the same is true for the restriction
 $\dist_E|_\Omega$ to any bounded open $\Omega\i A$.
\end{thm}

\Proof. Follows from lemma~\ref{lem:tuda-suda} and
theorem~\ref{thm:dist-extr}.\qeds

\begin{thm}{\bf Lemma about an obtuse angle.} \label{lem:tuda-suda} Given $v>0$, $r>0$,
$\kappa\in\RR$ and $m\in\NN$, there is $\eps=\eps(v,r,\kappa,m)>0$ such that if
$A\in\Alex^m(\kappa)$, $p\in A$, $\vol_m B_r(p)>v$, then for any two points
$x,y\in B_r(p)$, $|x y|<\eps$ there is point $z\in B_r(p)$ such that $\angle z x
y>\tfrac\pi2+\eps$ or $\angle z y x>\tfrac\pi2+\eps$.
\end{thm}
The proof is based on a volume comparison for $\log_x\:A\to T_x$ similar to \cite[lemma 1.3]{grove-petersen:finiteness} by Grove and Petersen.

Note that the tangent cone \label{T_pE}$T_p E$ of an extremal subset $E\i A$ is well defined; 
that is,
for any $p\in E$, the subsets $\lam\cdot E$ in $(\lam{\cdot} A,p)$ converge to a subcone of $T_p E\i T_p A$ as
$\lam\to\infty$.
Indeed, assume $E\i A$ is an extremal subset and $p\in E$.
For any $\xi\in \Sigma_p E$%
\footnote{\label{U_pX}For a closed subset $X\i A$, and $p\in X$, $\Sigma_p X\i \Sigma_p$ denotes the set of tangent directions to $X$ at $p$;
that is the
set of limits of~$\uparrow_p^{q_n}$ for~$q_n\to p$, $q_n\in X$.}%
, the radial curve $\gexp(t\cdot\xi)$ lies in $E$.%
The later follows from the fact that the curves 
$t\mapsto\gexp(t\,\cdot\!\uparrow_p^{q_n})$ starting from $q_n$ belong to $E$
and their converge to $\gexp(t\cdot\xi)$.
In particular, there is a curve which goes in any tangent direction of $E$.
Therefore, as $\lam\to\infty$,  $(\lam\cdot E\i \lam{\cdot} A,p)$ converges to
a subcone $T_p E\i T_p A$, 
which is the cone over $\Sigma_p E$ (see also
\cite[3.3]{perelman-petrunin:extremal})

Let us list some properties of tangent cones of extremal subsets:

\begin{enumerate}\setcounter{enumi}{\value{extr-prop}}

\item\label{ext-tangent} A closed subset $E\i A$ is  extremal if and only if 
the following condition is fulfilled:
\begin{enumerate}[$\diamond$]
\item At any point $p\in E$, its tangent cone $T_p E\i T_p A$ is well defined,
and it is an extremal subset of the tangent cone $T_p A$; compare
\cite[1.4]{perelman-petrunin:extremal}.
\end{enumerate}

(Here is an
equivalent formulation in terms of the space of directions: For any $p\in E$,
either 
(a) $\Sigma_p E=\emptyset$ and $\diam \Sigma_p\le\tfrac\pi2$ or 
(b) $\Sigma_p E=\{\xi\}$ is one point extremal subset and $\bar B_{\pi/2}(\xi)=\Sigma_p$
or
(c) $\Sigma_p E$ is extremal subset of $\Sigma_p$ with at least two points.)

$T_p E$ is extremal as a limit of extremal subsets, see
lemma~\ref{lem:limit-extr}. 
On the other hand for any semiconcave function $f$ and $p\in E$, the differential $d_p f\:T_p\to\RR$ is concave and since $T_p E\i T_p$ is extremal we have $\nabla_p f\z\in T_p E$. 
That is, 
gradient curves can be approximated by broken
geodesics with vertices on $E$, see page~\pageref{grad-constr}.

\item \cite[3.4--5]{perelman-petrunin:extremal} If $E$ and $F$ are extremal subsets then so are
\begin{enumerate}[(i)]
\item $E\cap F$ and for any $p\in E\cap F$ we have $T_p(E\cup F)=T_p E\cup \Sigma_p F$
\item $E\cup F$ and  for any $p\in E\cup F$ we have $T_p(E\cap F)=T_p E\cap \Sigma_p F$
\item \label{E-F}$\overline{E\backslash F}$\ \,\, and for any $p\in \overline{E\backslash F}$\ \,\, we have $T_p(\overline{E\backslash
F})=\overline{T_p E\backslash T_p F}$
\end{enumerate}
In particular, if $T_p E=T_p F$ then $E$ and $F$ coincide in a neighborhood of $p$.

The properties (i) and (ii) are obvious. 
The property (iii) follows from
property~\ref{ext-tangent} and lemma~\ref{lem:dist-to-extr}.
\setcounter{extr-prop}{\value{enumi}}
\end{enumerate}

\noi We continue with properties of the intrinsic metric of extremal subsets:
\begin{enumerate}
\setcounter{enumi}{\value{extr-prop}}

\item
\cite[3.2(3)]{perelman-petrunin:extremal} Let $A\in\Alex^m(\kappa)$ and $E\i A$ be an extremal subset. Then the induced 
metric of $E$ is locally bi-Lipschitz equivalent to its induced intrinsic
metric.
Moreover, the local 
Lipschitz constant at point $p\in E$ can be expressed in terms of $m$, $\kappa$
and volume of a ball $v=\vol B_r(p)$ for some (and therefore any) $r>0$.

From lemma~\ref{lem:tuda-suda}, it follows that for two sufficiently close
points $x,y\in E$ near $p$ there is a point $z$ so that
$\<\nabla_x\dist_z,\uparrow_x^y\> >\eps$ or $\<\nabla_y\dist_z,\uparrow_y^x\>
>\eps$. 
Then, for the corresponding point, say $x$, the gradient curve
$t\to\Phi^t_{\dist_z}(x)$ lies in $E$, it is 1-Lipschitz and the distance
$|\Phi^t_{\dist_z}(x)\,y|$ is decreasing with the speed of at least $\eps$. 
Hence the result.

\item 
Let $A_n\in\Alex^m(\kappa)$, $A_n\GHto A$ without collapse (that is $\dim A=m$) and $E_n\subset A_n$ be extremal subsets. 
Assume $E_n\to E\subset A$ as subsets. 
Then 
\begin{enumerate}[(i)]
\item\cite[9.1]{kapovitch:stability} \label{lim-dim-extr} For all large $n$, there is a homeomorphism of pairs $(A_n,E_n)\to(A,E)$. 
In particular, for all large $n$, $E_n$ is homeomorphic to $E$, 
\item\cite[1.2]{petrunin:extremal}\label{lim-intr-extr} $E_n\GHto E$ as length metric spaces (with the
intrinsic metrics induced from $A_n$ and $A$).
\end{enumerate}
The first property is a coproduct of the proof of Perelman's stability theorem.
The proof of the second is an application of quasigeodesics.

\item \cite[1.4]{petrunin:extremal}\label{1st-var}{\it The first variation formula.} Assume $A\in\Alex$ and
$E\i A$ is an extremal subset, let us denote by $|\!**|_E$ its intrinsic metric. 
Let $p,q \in E$ and $\alpha(t)$ be a curve in $E$ starting from $p$ in direction
$\alpha^+(0)\in \Sigma_p E$.
Then 
$$ |\alpha(t)\,q|_E=|p q|_E-\cos\phi\cdot t + o(t).$$
where $\phi$ is the minimal (intrinsic) distance in $\Sigma_p E$ between
$\alpha^+(0)$ and a direction of a shortest path in $E$ from $p$ to $q$ (if $\phi>\pi$, we assume $\cos\phi=-1$).

\item {\it Generalized Lieberman's Lemma.} Any minimizing geodesic for the
induced intrinsic metric on an extremal subset is a quasigeodesic in the ambient
space.

See \ref{lib-lem} for the proof and discussion.
\setcounter{extr-prop}{\value{enumi}}
\end{enumerate}

Let us denote by $\Ext(x)$ the minimal extremal subset which contains a point
$x\in A$. 
Extremal subsets which can be obtained this way will be called
\emph{primitive}. 
Set $$\Ext^\circ(x)=\{y\in A|\Ext(y)=\Ext(x)\};$$
the set $\Ext^\circ(x)$ is $\Ext(x)$ with its proper extremal subsets removed.
Let us call $\Ext^\circ(x)$ the \emph{main part} of $\Ext(x)$.
From the property~\ref{E-F} on page~\pageref{E-F},  $\Ext^\circ(x)$ is open and everywhere
dense in $\Ext(x)$.
Clearly the main parts of primitive extremal subsets form a disjoint
covering of $M$.

\begin{enumerate}
\setcounter{enumi}{\value{extr-prop}}
\item \cite[3.8]{perelman-petrunin:extremal} {\it Stratification.}\label{strata} The main part of a primitive extremal subset is a
topological manifold. In particular, the main parts of primitive extremal subsets
stratify Alexandrov's space into topological manifolds. 

This follows from theorem~\ref{thm:dist-extr} and the Morse lemma (property~\ref{morse} page~\pageref{morse}); see also example~\ref{ex:t-cone}, page~\pageref{ex:t-cone}.
\setcounter{extr-prop}{\value{enumi}}
\end{enumerate}

\subsection{Applications}

Extremal subsets make possible to give a more precise formulations of some known theorems.
Here is the simplest example, a version of the radius sphere theorem:

\begin{thm}{\bf Theorem.}\label{thm:extr-sph}
Let $A\in\Alex^m(1)$, $\diam A>\tfrac\pi2$ and $A$ have no extremal subsets. Then $A$
is homeomorphic to a sphere.
\end{thm}

From lemma~\ref{lemma:rad-big-U_p} and theorem~\ref{thm:dist-extr}, we have
$A\in\Alex(1)$, $\Rad A>\tfrac\pi2$ implies that $A$ has no extremal subsets. 
That is, 
this theorem does indeed generalize the radius sphere
theorem~\ref{cor:no-ext+rad-sph}(\ref{rad-sph}).

\Proof. Assume $p,q\in A$ realize the diameter of $A$.
Since  $A$ has no extremal subsets, from example~\ref{ex:t-cone}, page~\pageref{ex:t-cone}, it  follows that a small spherical neighborhood
of $p\in A$ is homeomorphic to $\RR^m$. 
From angle comparison, $\dist_p$ has only two critical points $p$ and $q$. 
Therefore, this theorem follows from the Morse lemma  (property~\ref{morse} page~\pageref{morse}) applied to $\dist_p$. \qeds

The main result of such type is the result in
\cite{perelman:collapsing}. 
It roughly states that a collapsing to a compact
space without proper extremal subsets carries a natural Serre bundle structure.

This theorem is analogous to  the following:
\begin{thm}{\bf Yamaguchi's fibration theorem
\cite{yamaguchi:bundle}.}
Let $A_n\in\Alex^m(\kappa)$ and $A_n\GHto M$, $M$ be a Riemannian manifold.

Then there is a sequence of locally trivial fiber bundles $\sigma_n\:A_n\to M$.
Moreover, $\sigma_n$ can be chosen to be \emph{almost submetries}%
\footnote{that is, a Lipshitz and co-Lipschitz with constants almost 1.} 
and the diameters of
its fibers converge to $0$.
\end{thm}
The conclusion in Perelman's theorem is weaker, but on the other hand it is just
as good for practical purposes. 

Here is a source of examples: 
of a collapse to the spaces with
extremal subsets which do not have the homotopy lifting property. 
Take a compact Riemannian manifold $M$ with an isometric and
non-free action by a compact connected Lie group $G$, then $(M\times \eps
{\cdot}G)/G\GHto M/G$ as $\eps\to0$ and since the curvature of $G$ is non-negative, by
O'Naill's formula, we get that the curvature of $(M\times \eps {\cdot}G)/G$ is uniformly bounded
below.

\begin{thm}
{\bf Homotopy lifting theorem.}
\label{thm:per-ser} Let $A_n\GHto A$, $A_n\in\Alex^m(\kappa)$, $A$ be compact
without proper extremal subsets and $K$ be a finite simplicial complex. 

Then, given a homotopy 
$$F_t\:K\to A,\ \ t\in [0,1]$$ 
and a sequence of maps 
$G_{0;n}\:K\to A_n$ such that $G_{0,n}\to F_0$ as $n\to\infty$ one can extend
$G_{0;n}$ by homotopies 
$$G_{t;n}\:K\to A$$
such that $G_{t;n}\to F_t$ as $n\to\infty$.
\end{thm}

An alternative proof is based on Lemma~\ref{lem:hom-approx}.

\begin{thm}{\bf Remark.}\label{rem:hom-seq} {\rm 
As a corollary of this theorem one obtains that for all large $n$
it is possible to write a homotopy exact sequence:

$$\cdots\pi_k(F_n)
\longrightarrow\pi_k(A_n) \longrightarrow\pi_k(A) 
\longrightarrow\pi_{k-1}(F_n)\cdots,$$
where the space $F_n$ can be obtained the following way:
Take a point $p\in A$, and fix $\eps>0$ so that $\dist_p\:A\to\RR$ has no critical values
in the interval $(0,2{\cdot}\eps)$. 
Consider a sequence of points $A_n\ni p_n\to p$ and take 
$F_n=B_\eps(p_n)\i A_n$.
In particular, if $p$ is a regular point then for large $n$, $F_n$ is homotopy
equivalent to a \emph{regular fiber over $p$}%
\footnote{\label{reg-fib} The regular fiber is constructed the
following way: take a distance chart $G\:B_{2{\cdot}\eps}(p)\to \RR^{k}$, $k=\dim A$ around $p\in A$
and lift it to $A_n$. 
It defines a map $G_n\:B_\eps(p_n)\to\RR^k$.
Then take
$F_n=G_n^{-1}\circ G(p)$ for large $n$. 
If $A_n$ are Riemannian then $F_n$ are manifolds and they do not depend on $p$ up to a homeomorphism.
Moreover, $F_n$ are almost non-negatively curved in a generalized sense; 
see \cite[definition 1.4]{KPT}.}. }
\end{thm}

Next we give two corollaries of the above
remark.
The last assertion of the following theorem was conjectured by Shioya in
\cite{shioya} and  was proved by Mendon\c{c}a in \cite{mendonca:shioya}.

\begin{thm}{\bf Theorem \cite[3.1]{perelman:collapsing}.} 
Let $M$ be a complete noncompact Riemannian manifold of
nonnegative sectional curvature. 
Assume that its asymptotic cone $\Cone_\infty(M)$
has no proper extremal subsets, then $M$ splits isometrically into the product
$L\times N$, where $L$ is a compact Riemannian manifold and $N$ is a non-compact
Riemannian manifold of the same dimension as $\Cone_\infty(M)$.

In particular, the same conclusion holds if radius of the ideal boundary of $M$ is
at least $\tfrac\pi2$.
\end{thm}

The proof is a direct application of theorem~\ref{thm:per-ser} and remark~\ref{rem:hom-seq}
for collapsing 
$$\eps{\cdot} M\GHto\Cone_\infty(M),\ \ \text{as}\ \ \eps\to0.$$

\begin{thm}{\bf Theorem \cite[3.2]{perelman:collapsing}.}
Let $A_n\in\Alex^m(1)$, $A_n\GHto A$ be a collapsing sequence 
(that is $m>\dim A$), 
then $\Cone(A)$ has proper
extremal subsets.
In particular, $\Rad A\le\tfrac\pi2$.
\end{thm}

The last assertion of this theorem (in a stronger form) has been proven 
by Grove and Petersen in \cite[3(3)]{grove-petersen:rad-sphere}.

The proof is a direct application of theorem~\ref{thm:per-ser} and remark~\ref{rem:hom-seq}
for collapsing of spherical suspensions
$$\Sigma(A_n)\GHto\Sigma(A),\ \ n\to\infty.$$

\section{Quasigeodesics}\label{QG}

The class of quasigeodesics
generalizes the class of geodesics to nonsmooth metric spaces.
It was first introduced by Alexandrov in \cite{alexandrov:qg} for
$2$-dimensional convex hypersurfaces in the Euclidean space, as the curves which
``turn'' right and left simultaneously.
They were studied further by Alexandrov, Burago, Pogorelov and Milka
in \cite{alexandrov-burago},  \cite{pogorelov:qg}
\cite{milka:qg}.
It was generalized to surfaces with bounded integral
curvature by Alexandrov, see \cite{alexandrov:int-qg}, 
and to multidimensional polyhedral spaces by Milka, see \cite{milka:poly1},
\cite{milka:poly2}.
For multi-dimensional Alexandrov's spaces they were introduced in the my master thesis, see also \cite{perelman-petrunin:qg} by Perelman and me.

It should be noted that the class of quasigeodesics
described here has nothing to do with the Gromov's quasigeodesics in
$\delta$-hyperbolic spaces. 

In Alexandrov's spaces, quasigeodesics behave more naturally than geodesics, mainly: 
\begin{enumerate}[$\diamond$]
\item There is a quasigeodesic starting in any direction from any point; 
\item The limit of quasigeodesics is a quasigeodesic.
\end{enumerate}

Quasigeodesics have beauty on their own, but also
due to the generalized Lieberman lemma (\ref{lib-lem}), they are very useful in the
study of intrinsic metric of extremal subsets, in particular the boundary of
Alexandrov's space. 

Since quasigeodesics behave almost as geodesics, they are often used instead of geodesics in
the situations when there is no geodesic in a given direction.
In most of these applications one can instead use the radial curves of gradient
exponent, see section~\ref{gexp}; 
a good example is the proof of
theorem~\ref{thm:dist-to-bry}, see footnote~\ref{qg-grad},
page~\pageref{qg-grad}.
In this type of argument, radial curves could be considered as a simpler and
superior tool since they can be defined for infinitely dimensional Alexandrov's spaces.

\subsection{Definition and properties}

In section~\ref{CF}, we defined $\lambda$-concave functions  as those locally
Lipschitz functions whose restriction to any unit-speed minimizing geodesic is
$\lambda$-concave. 
Now consider a curve $\gamma$ in an Alexandrov's space such that restriction of any
$\lambda$-concave function to $\gamma$ is $\lambda$-concave.
It is easy to see that for any Riemannian manifold $\gamma$ has to be a unit-speed
geodesic. 
In a general Alexandrov's space $\gamma$ should only be a
quasigeodesic.

\begin{thm}{\bf Definition.} A curve $\gamma$ in an Alexandrov's space is called
\emph{quasigeodesic} if for any $\lambda\in\RR$, given a $\lambda$-concave
function $f$, the real-to-real function $f\circ\gamma$ is $\lambda$-concave.
\end{thm}

Although this definition works for any metric space, it is only reasonable to
apply it for the spaces where we have $\lambda$-concave functions, and
Alexandrov's spaces seem to be the perfect choice.

The following is a list of corollaries from this definition:
\newcounter{qg-prop}
\begin{enumerate}
\item Quasigeodesics are unit-speed curves.
That is, if $\gamma(t)$ is a quasigeodesic then for any $t_0$ we have 
$$ \lim_{t\to t_0}\frac{|\gamma(t)\gamma(t_0)|}{|t-t_0|}=1.$$ 

To prove that quasigeodesic $\gamma$ is $1$-Lipschitz at some $t=t_0$,
it is enough to apply the definition for $f=\dist_{\gamma(t_0)}^2$ and use the fact
that in any Alexandrov's space $\dist_p^2$ is $(2+O(r^2))$-concave in
$B_r(p)$.
The lower bound is more complicated, see theorem~\ref{thm:unit-speed}.

\item For any quasigeodesic the right and left tangent vectors $\gamma^+$,
$\gamma^-$ are uniquely defined unit vectors. 

To prove, take a partial limits $\xi^\pm\in T_{\gamma(t_0)}$ for
$$\frac{\log_{\gamma(t_0)}\gamma(t_0\pm\tau)}{\tau},\ \ \text{as}\ \ 
\tau\to0+$$
It exists since quasigeodesics are 1-Lipschitz (see the previous property).
For any semiconcave function $f$,  $(f\circ\gamma)^\pm$ are well defined,
therefore 
$$(f\circ\gamma)^\pm(t_0)=d_{\gamma(t_0)}f(\xi^\pm).$$
Taking $f=\dist_q^2$ for different $q\in A$, one can see that $\xi^\pm$ is
defined uniquely by this identity, and therefore $\gamma^\pm(t_0)=\xi^\pm$.

\item {\it Generalized Lieberman's Lemma.} Assume $\gamma$ be a unit-speed geodesic in an extremal subset $E$ equipped with induced intrinsic metric.
Then $\gamma$ is a quasigeodesic ambient
Alexandrov's space
\setcounter{extr-prop}{\value{enumi}}

See \ref{lib-lem} for the proof and discussion.

\item\label{exist-qg} For any point $x\in A$, and any direction $\xi\in \Sigma_x$
there is a quasigeodesic $\gamma\:\RR\to A$ such that $\gamma(0)=x$ and
$\gamma^+(0)=\xi$.

Moreover, if $E\i A$ is an extremal subset and $x\in E$, $\xi\in \Sigma_x E$,
then $\gamma$ can be chosen to lie completely in $E$.

The proof is  given in the appendix~\ref{constr-qg}.

\setcounter{qg-prop}{\value{enumi}}
\end{enumerate}

Applying the definition locally, 
we get that if $f''+\kappa{\cdot}f\le 1$ then $f\circ\gamma$ is $(1-\kappa{\cdot}f\circ\gamma)$-concave (see section~\ref{variation-CF}). 
In particular, if $A$ is an Alexandrov's space with curvature $\ge\kappa$, $p\in
A$ and $h_p(t)=\rho_\kappa\circ\dist_p\circ\gamma(t)$%
\footnote{Function
$\rho_\kappa\:\RR\to\RR$ is defined on page~\pageref{rho_k}} 
then we have the
following inequality in the \emph{barrier sense}
$$h_p''\le 1- \kappa{\cdot} h_p.$$

This inequality can be reformulated in an equivalent way:
Let $A\in \Alex^m(\kappa)$, $p\in A$ and $\gamma$ be a quasigeodesic, then
function
$$t\mapsto \tilde\angle_\kappa(|\gamma(0)p|,|\gamma(t)p|,t)$$
is decreasing for any $t>0$ (if $\kappa>0$ then one has to assume
$t\le\pi/\sqrt\kappa$).

In particular, 
$$
\tilde\angle_\kappa(|\gamma(0)p|,|\gamma(t)p|,t)\le\angle(\uparrow_{\gamma(0)}^p,\gamma^+(0))$$ 
for any $t>0$ (if $\kappa>0$ then in addition $t\le\pi/\sqrt\kappa$).

Let us give a more geometric property using the notion of {\it developing} defined below. 

{\it Any quasigeodesic in an Alexandrov's space with curvature $\ge \kappa$, has a
convex $\kappa$-developing with respect to any point.}

\begin{thm}{\bf Definition of developing.} 
Fix a real $\kappa$. 

Let $X$ be a metric space, $\gamma\:[a,b]\to X$ be a 1-Lipschitz curve and $p\in
X\backslash\gamma$. 
If $\kappa>0$, assume in addition that $|p\gamma(t)| < \pi/ \sqrt{\kappa}$ for
all 
$t\in [a,b]$. 

Then there exists a unique (up to rotation) curve
$\tilde\gamma\: [a,b]\to \Lob_\kappa$, parametrized by the arclength, and such
that
$|o\tilde\gamma(t)|=|p\gamma(t)|$ for all $t$ and some fixed $o\in \Lob_\kappa$,
and the segment
$o\tilde\gamma(t)$ turns clockwise as $t$ increases
(this is easy to prove).
Such a curve $\tilde\gamma$ is called the \emph{$\kappa$-development of $\gamma$
with respect to $p$}. 

The development $\tilde\gamma$ is called \emph{convex} if for every $t\in
(a,b)$, for sufficiently small $\tau > 0$ the
curvilinear triangle, bounded by the segments $o\tilde\gamma(t\pm\tau)$ and the
arc $\tilde\gamma|_{t-\tau,t+\tau}$, is convex. 
\end{thm}

This definition is given by Alexandrov in \cite{alexandrov:devel};
it is based on earlier construction in \cite{liberman} by Liberman.
In \cite{milka:qg}, Milka shows that the developing of a quasigeodesic on a convex surface is convex.

\begin{enumerate}\setcounter{enumi}{\value{qg-prop}}

\item \label{thm:eq-def-qg}  Let $A\in\Alex^m(\kappa)$,
$m>1$%
\footnote{This condition is only needed to ensure that the set
$A\backslash\gamma$ is everywhere dense.}.
A curve $\gamma$ in $A$ is a quasigeodesic if and only if it is parametrized by
arc-length and one of the following properties is fulfilled:
\begin{enumerate}[(i)]
\item For any point $p\in A\backslash\gamma$ the $\kappa$-developing of $\gamma$
with respect to $p$ is convex.
\item \label{part:h}For any point $p\in A$, if
$h_p(t)=\rho_\kappa\circ\dist_p\circ\gamma(t)$, then we have the following
inequality in \emph{a barrier sense}
$$h_p''\le 1- \kappa{\cdot} h_p.$$
\item Function
$$t\mapsto \tilde\angle_\kappa(|\gamma(0)p|,|\gamma(t)p|,t)$$
is decreasing for $t>0$.
\item\label{comp-inq} The inequality
$$\angle(\uparrow_{\gamma(0)}^p,\gamma^+(0))\ge
\tilde\angle_\kappa(|\gamma(0)p|,|\gamma(t)p|,t)$$ 
holds for all small $t>0$.
\end{enumerate}
The ``only if'' part has already been proven above, and the ``if'' part follows
from corollary~\ref{cor:eq-qg}

\item\label{pr:limit-qg} A pointwise limit of quasigeodesics is a quasigeodesic.
More generally:

\label{thm:limit-qg}
\textit{Assume $A_n\GHto A$, $A_n\in\Alex^m(\kappa)$, $\dim A=m$ 
(that is, it is not a
collapse).\\
Let $\gamma_n\:[a,b]\to A_n$ be a sequence of quasigeodesics which converges
pointwise to a curve $\gamma\:[a,b]\to A$.
Then $\gamma$ is a quasigeodesic.}

As it follows from lemma~\ref{lem:lifting}, the statement in the
definition is correct for any $\lambda$-concave function $f$ which has
controlled convexity type $(\lambda,\kappa)$.
That is, $\gamma$ satisfies the property~\ref{propr:weak-qg}. 
In particular, the $\kappa$-developing of $\gamma$ with respect to any point $p\in A$ is convex, and as it is noted in remark~\ref{rmk:stronger-unit-speed}, $\gamma$ is a
unit-speed curve.
Therefore, from corollary~\ref{cor:eq-qg} we get that it is a quasigeodesic.

\setcounter{qg-prop}{\value{enumi}}
\end{enumerate}

\noi Here is a list of open problems on quasigeodesics:

\begin{enumerate}[(i)]
\item Is there an analog of the Liouvile theorem for ``quasigeodesic flow''?
\item Is it true that any finite quasigeodesic has bounded variation of turn?  

or

Is it possible to approximate any finite quasigeodesic by sequence of
broken lines with bounded variation of turn?

\item Is it true that in an Alexandrov's space without boundary there is an
infinitely long geodesic?
\end{enumerate}

As it was noted by A.~Lytchak, the first and last questions can be reduced to the following:
Assume $A$ is a compact Alexandrov's $m$-space without boundary. Let us set
$V(r)=\int_A\vol_m(B_r(x))$,
then $V(r)=\vol_m(A)\omega_m r^m+o(r^{m+1})$.
The technique of \emph{tight maps} makes it possible to prove only that
$V(r)\z=\vol_m(A)\omega_m r^m+O(r^{m+1})$.
Note that if $A$ is a Riemannian manifold with boundary then $V(r)=\vol_m(A)\omega_m r^m+\vol_{m-1}(\partial A)\omega'_m
r^{m+1}+o(r^{m+1})$.

\subsection{Applications.}

The quasigeodesics is the main technical tool in the questions linked to the intrinsic
metric of extremal subsets, in particular the boundary of Alexandrov's space. 
The main examples are
the proofs of convergence of intrinsic metric of extremal subsets and
the first variation formula (see properties~\ref{lim-intr-extr} and \ref{1st-var}, on page~\pageref{lim-intr-extr}).
 
Below  we give a couple of simpler examples:

\begin{thm}{\bf Lemma.}\label{lemma:rad-big-U_p}
Let $A\in\Alex^m(1)$  and  $\Rad A>\tfrac\pi2$. 
Then for any $p\in A$
the space of directions $\Sigma_p$ has radius $> \tfrac\pi2$.
\end{thm}

\Proof. Assume that $\Sigma_p$ has radius
$\leq\tfrac\pi2$, and let $\xi\in \Sigma_p$ be a direction, such that $\bar B_\xi
(\tfrac\pi2)=\Sigma_p$. 
Consider a quasigeodesic $\gamma$ starting at $p$ in direction $\xi$. 

Then for $q=\gamma(\tfrac\pi2)$ we have $\bar B_q(\tfrac\pi2)=A$. Indeed, for
any point $x\in A$ we have 
$\angle(\xi,\uparrow_p^x)\le\tfrac\pi2$. Therefore, by the comparison inequality
(property~\ref{comp-inq}, page~\pageref{comp-inq}),
$|x q|\le\tfrac\pi2$. 
This contradicts our assumption that $\Rad A> \tfrac\pi2$. \qeds

\begin{thm}{\bf Corollary.}\label{cor:no-ext+rad-sph}
Let $A\in\Alex^m(1)$ and  $\Rad A>\tfrac\pi2$ then 
\begin{enumerate}[(i)]
\item\label{no-ext} $A$ has no extremal subsets.
\item\label{rad-sph}(radius sphere theorem) $A$ is homeomorphic to an
$m$-sphere.
\end{enumerate}
\end{thm}

Part (\ref{rad-sph}) was proved by Grove and Petersen in \cite{grove-petersen:rad-sphere}.
Another proof follows immediately from \cite[1.2, 1.4.1]{perelman-petrunin:extremal}; 
theorem~\ref{thm:extr-sph} gives a slight generalization.

\Proof. Part (\ref{no-ext}) is obvious. 

Part (\ref{rad-sph}): From lemma \ref{lemma:rad-big-U_p}, $\Rad \Sigma_p>\tfrac\pi2$. 
Since $\dim \Sigma_p<m$, by the induction hypothesis we have $\Sigma_p\simeq S^{m-1}$.
Now the Morse lemma (see property~\ref{morse}, page~\pageref{morse}) for
$\dist_p\:A\to\RR$ gives that $A\simeq\Sigma(\Sigma_p)\simeq S^m$, here $\Sigma(\Sigma_p)$ denotes a spherical suspension over $\Sigma_p$.\qeds

\section{Simple functions}
\label{adm}
\addtocounter{subsection}{1}
\setcounter{thm}{0}

This is a short technical section.
Here we introduce \emph{simple functions}, a subclass of semiconcave
functions which on one hand includes all functions we need and in addition is liftable; 
that is, for any such function one can construct a nearby function on a
nearby space with ``similar'' properties.

Our definition of simple function is a modification of two different definitions of so called ``admissible functions'' given by Perelman in
\cite[3.2]{perelman:morse}, see also \cite[5.1]{kapovitch:stability}. 

\begin{thm}{\bf Definition} Let $A\in \Alex$, a function $f\:A\to \RR$ is called
\emph{simple} if there is a finite set of points $\{q_i\}_{i=1}^N$ and a
semiconcave function $\Theta\:\RR^N\to\RR$ which is non-decreasing in each
argument such that
$$f(x)= \Theta(\dist_{q_1}^2,\dist_{q_2}^2,\dots,\dist_{q_N}^2)$$
\end{thm}

It is straightforward to check that simple functions are semiconcave. Class of simple functions is
closed under summation, multiplication by a positive constant
and taking the minimum (as well as multiplication by positive simple functions).

In addition this class is liftable; 
that is, given a converging sequence of
Alexandrov's spaces $A_n\GHto A$ and a simple function $f\:A\to\RR$ there is a
way to construct a sequence of functions $f_n\:A_n\to\RR$ such that $f_n\to f$.
Namely, for each $q_i$ take a sequence $A_n\ni q_{i,n}\to q_i\in A$ and consider
function $f_n\:A_n\to\RR$ defined by
$$f_n=\Theta(\dist_{q_{1,n}}^2,\dist_{q_{2,n}}^2,\dots,\dist_{q_{N,n}}^2).$$

\subsection{Smoothing trick.}\label{smooth}
Here we present a trick which is very useful for doing local analysis in
Alexandrov's spaces, it was introduced by Otsu and Shioya in \cite[section 5]{otsu-shioya}.

Consider function
$$\widetilde\dist_p=\oint\limits_{B_\eps(p)}\dist_x\cdot \d x.$$
In this notation, we do not specify $\eps$ assuming it
to be very small.

It is easy to see that $\widetilde\dist_p$ is semiconcave.
Note also that 
$$\d_y\widetilde\dist_p
=
\oint\limits_{B_\eps(p)}d_y\dist_x\cdot \d x.$$
If $y\in A$ is regular, 
that is $T_y$ is isometric to Euclidean space, 
then for almost all $x\in B_\eps(p)$ the differential $d_y\dist_x\:T_y\to\RR$ is
a linear function.
Therefore $\widetilde\dist_p$ is differentiable at every regular point.
That is 
$$\d_y\widetilde\dist_p\:T_y\to\RR$$
is a linear function for any regular $y\in A$.

The same trick can be applied to any simple function
$$f(x)=\Theta(\dist_{q_1}^2,\dist_{q_2}^2,\dots,\dist_{q_N}^2).$$
This way we obtain function
$$\tilde f(x)=
\oint_{B_\eps(q_1)\times B_\eps(q_2)\times\cdots\times B_\eps(q_N)}
\Theta(\dist_{x_1}^2,\dist_{x_2}^2,\dots,\dist_{x_N}^2)
\cdot d x_1\cdot d x_2\cdots d x_N,$$
which is differentiable at every regular point;
that is, if $T_y$ is isometric to
the Euclidean space then 
$$\d_y\tilde f\:T_y\to\RR$$
is a linear function.

\section{Controlled concavity}
\addtocounter{subsection}{1}
\setcounter{thm}{0}

In this and the next sections we introduce a couple of techniques which use
comparison of $m$-dimensional Alexandrov's space with a model space of the same
dimension \label{lob-k-m}$\Lob_\kappa^m$ (that is $m$-dimensional simply connected Riemannian manifold with constant curvature $\kappa$).
These techniques were introduced by Perelman in \cite{perelman:morse}
and \cite{perelman:DC}.

We start with the local existence of a strictly concave function on an Alexandrov's space.

\begin{thm}{\bf Theorem~\cite[3.6]{perelman:morse}.}
\label{thm:strictly-concave}
Let $A\in\Alex$. 

For any point $p\in A$ there is a strictly concave function $f$ defined in an
open neighborhood of $p$.

Moreover, given $v\in T_p$, the differential, $d_p f(x)$, can be chosen
arbitrarily close to $x\mapsto -\<v,x\>$
\end{thm}

\begin{wrapfigure}{r}{35mm}
\begin{lpic}[t(-10mm),b(0mm),r(0mm),l(0mm)]{pics/strictly-concave(0.25)}
\lbl[tr]{125,05;$q$}
\lbl[br]{32,166;$\gamma(t)$}
\lbl[l]{44,159;$\alpha(t)$}
\end{lpic}
\end{wrapfigure}

\Proof. 
Consider the real-to-real function 
$$\phi_{r,c}(x)=(x-r)- c{(x-r)^2}/r,$$
so we have 
$$\phi_{r,c}(r)=0,\ \ \phi_{r,c}'(r)=1\ \ \phi_{r,c}''(r)=- {2c}/{r}.$$ 

Let $\gamma$ be a unit-speed geodesic, fix a point $q$ and set 
$$\alpha(t)=\angle(\gamma^+(t),\uparrow_{\gamma(t)}^{q}).$$
If $r>0$ is sufficiently small and $|q\gamma(t)|$ is sufficiently close to
$r$, then direct calculations show that
$$(\phi_{r,c}\circ\dist_q\circ\gamma)''(t)
\le 
\frac{3-c\cdot \cos^2\alpha(t)}{r}.$$

Now, assume $\{q_i\}$, $i=\{1,..,N\}$ is a finite set of points such that $|p q_i|=r$ for any $i$. 
For $x\in A$ and $\xi_x\in \Sigma_x$, set $\alpha_i(\xi_x)=\angle(\xi_x,\uparrow_p^{q_i})$. 
Assume we have a collection $\{q_i\}$ such
that for any $x\in B_\eps(p)$ and $\xi_x\in \Sigma_x$ 
we have  $\max_i\{|\alpha_i(\xi_x)-\tfrac\pi2|\}\ge\eps>0$. 
Then  taking in the above inequality $c>3N/\cos^2\eps$, we get that the function
$$f=\sum_i \phi_{r,c}\circ\dist_{q_i}$$
is strictly concave in $B_{\eps'}(p)$ for some positive $\eps'<\eps$.

To construct the needed collection $\{q_i\}$, note that for small $r>0$ one can
construct $N_\delta\ge \Const/\delta^{(m-1)}$ points $\{q_i\}$ such that $|p q_i|=r$
and $\tilde\angle_\kappa q_i p q_j>\delta$ (here $\Const=\Const(\Sigma_p)>0$).
On the other hand, the set of directions which is orthogonal to a given direction
is smaller than $S^{m-2}$ and therefore contains at most
$\Const(m)/\delta^{(m-2)}$ directions with angles at least $\delta$. 
Therefore, for small enough $\delta>0$, $\{q_i\}$ forms the needed collection.

If $r$ is small enough, points $q_i$ can be chosen so that all directions
$\uparrow_p^{q_i}$ will be $\eps$-close to a given direction $\xi$ and
therefore the second property follows.
\qeds

Note that in the theorem~\ref{thm:strictly-concave} (as well as in
theorem~\ref{exist-control}), the function $f$ can be chosen to have maximum value $0$ at $p$,
$f(p)=0$ and with $d_p f(x)$ arbitrary close to $-|x|$.
It can be constructed by taking the minimum of the functions in these theorems. Whence the claim below follows.

\begin{thm}{\bf Claim.}\label{cor:convex-nbhd}
For any point of an Alexandrov's space there is an arbitrary small closed convex
neighborhood.
\end{thm}

Applying rescaling and passing to the limit, 
one can estimate the size of the convex
hull in an Alexandrov's space in terms of the volume of a ball containing it:

\begin{thm}{\bf Lemma on strictly concave convex hulls
\cite[4.3]{perelman-petrunin:extremal}.}
For any $v>0$, $r>0$ and $\kappa\in \RR$, $m\in\NN$ there is $\eps>0$ such that,
if $A\in\Alex^m(\kappa)$ and $\vol B_{r}(p)\ge v$ then for any $\rho<\eps\cdot r$, 
$$\diam\Conv B_\rho(p)\le \rho/\eps.$$

In particular, for any compact Alexandrov's $A$ space there is $\Const\in \RR$ such that for
any subset $X\i A$
$$\diam \l(\Conv X \r)\le \Const\cdot\diam X.$$
\end{thm}

\subsection{General definition.} 
The above construction can be generalized and optimized in many ways to fit
particular needs. 
Here we introduce one such variation which is not the most general, but general
enough to work in most applications.

Let $A$ be an Alexandrov's space and $f\:A\to \RR$,
$$f=\Theta(\dist^2_{q_1},\dist^2_{q_2},\dots,\dist^2_{q_N})$$ be a \emph{simple
function} (see section~\ref{adm}).
If $A$ is $m$-dimensional, we say that such a function $f$ has \emph{controlled
concavity of type} $(\lambda,\kappa)$ at $p\in A$, if for any $\eps>0$ there is
$\delta>0$, such that for any collection of points $\{\tilde p,\tilde q_i\}$ in
the \emph{model $m$-space}%
\footnote{that is, a simply connected $m$-manifold with
constant curvature $\kappa$.} $\Lob_\kappa^m$  satisfying 
$$|\tilde q_i\tilde q_j|> |q_i  q_j|-\delta\ \ \text{and}\ \  \bigl||\tilde
p\tilde q_i|-|p q_i|\bigr|<\delta \ \ \text{for all}\ \ i,j,$$ 
we have that the function 
$\tilde f\:\Lob_\kappa^m\to \RR$ defined by 
$$\tilde f=\Theta(\dist^2_{\tilde q_1},\dist^2_{\tilde q_2},..,\dist^2_{\tilde
q_n})$$
is $(\lambda-\eps)$-concave in a small neighborhood of $\tilde p$.

The following lemma states that the controlled concavity is stronger than the usual
concavity.

\begin{thm}{\bf Lemma.} \label{contr-concave}
Let $A\in\Alex^m(\kappa)$.

If a simple function 
$$f=\Theta(\dist^2_{q_1},\dist^2_{q_2},..,\dist^2_{q_N}),\ \ f\:A\to\RR$$ 
has a controlled concavity type $(\lambda,\kappa)$ at each point $p\in \Omega$, then
$f''\le \lambda$ in $\Omega$. 
\end{thm}

The proof is just a direct calculation similar to that in the proof
of~\ref{thm:strictly-concave}. 
Note also, that the function constructed in the proof of
theorem~\ref{thm:strictly-concave} has controlled concavity.
In fact from the
same proof follows:

\begin{thm}{\bf Existence.} \label{exist-control}
Let $A\in\Alex$, $p\in A$, $\lambda,\kappa\in\RR$.
Then there is a function $f$ of controlled concavity $(\lambda,\kappa)$
at $p$.

Moreover, given $v\in T_p$, the function $f$ can be chosen so that its differential $d_p f(x)$ will be arbitrary close to $x\mapsto -\<v,x\>$.
\end{thm}

Since functions with a controlled concavity are simple they admit liftings, and
from the definition it is clear that these liftings also have controlled
concavity of the same type.
More precisely, we get the following.

\begin{thm}{\bf Concavity of lifting.} \label{lem:lifting}
Let $A\in\Alex^m$.

Assume a simple function 
$$f\:A\to \RR,\ \ f=\Theta(\dist^2_{q_1},\dist^2_{q_2},..,\dist^2_{q_N})$$
has controlled concavity type $(\lambda,\kappa)$ at $p$.

Let $A_n\in\Alex^m(\kappa)$,
$A_n\GHto A$ (so, no collapse) and $\{p_n\},\{q_{i,n}\}\in
A_n$ be sequences of points such that $p_n\to p\in A$ and  $q_{i,n}\to q_i\in A$
for each $i$.

Then for all large $n$, the liftings of $f$,
$$f_n\:A_n\to \RR,\ \
f_n=\Theta(\dist^2_{q_{1,n}},\dist^2_{q_{2,n}},..,\dist^2_{q_{N,n}})$$
have controlled concavity type $(\lambda,\kappa)$ at $p_n$. 
\end{thm}

\subsection{Applications} 
\label{app-con-con}

As was already noted, in the theorems~\ref{thm:strictly-concave} and
\ref{exist-control}, the function $f$ can be chosen to have a maximum value $0$ at $p$,
and with $d_p f(x)$ arbitrary close to $-|x|$. 
This observation was used by Kapovitch in \cite{kapovitch:regularity} to solve the second part of Petersen's problem 32 from \cite{petersen:list}:

\begin{thm}
{\bf Petersen's problem.}\label{smoothable}
Let $A$ be a smoothable Alexandrov's $m$-space; 
that is,
there is a sequence of Riemannian $m$-manifolds $M_n$ with curvature $\ge\kappa$
such that $M_n\GHto A$.

Prove that the space of directions $\Sigma_x A$ for any point $x\in A$ is
homeomorphic to the standard sphere.
\end{thm}

Note that Perelman's stability theorem (see  \cite{perelman:spaces2},
\cite{kapovitch:stability}) only gives that $\Sigma_x A$ has to be homotopically
equivalent to the standard sphere.

\bigskip
\noi\textit{Sketch of the proof:}
Fix a big negative $\lambda$ and construct a function $f\:A\to\RR$ with $d_p
f(x)\approx -|x|$ and controlled concavity of type $(\lambda,\kappa)$.
From  \ref{contr-concave}, the liftings  $f_n\:M_n\to\RR$ of  $f$ (see
\ref{lem:lifting}) are strictly concave for large $n$.
Let us slightly smooth the functions $f_n$ keeping them strictly concave.
Then the level sets $f^{-1}_n(a)$, for values of $a$, which are little below the maximum of $f_n$,
have strictly positive curvature and are diffeomorphic to the standard
sphere%
\footnote{Since $f$ has only one critical value above $a$ and it is a local maximum.}.

Let us denote by $p_n\in M_n$ a maximum point of $f_n$.
Then it is not hard to choose a sequence $\{a_n\}$ and a sequence of rescalings
$\{\lam_n\}$ so that $(\lam_n{\cdot}M_n,p_n)\GHto (T_p,o_p)$ and 
$\lam_n\cdot f^{-1}_n(a_n)\i \lam_n {\cdot}M_n$ converge to a
convex hypersurface $S$ close to $\Sigma_p\i T_p$.
Then, from Perelman's stability theorem, it follows that $S$ and therefore $\Sigma_p$
is homeomorphic to the standard sphere.
\qeds

\noi{\bf Remark.} From this proof it follows that $\Sigma_p$ is
itself smoothable. Moreover, there is a non-collapsing sequence of Riemannian metrics $g_n$ on $S^{m-1}$ such that $(S^{m-1},g_n)\GHto\Sigma_p$. 
This observation makes possible to proof a similar statement for iterated spaces of directions of smoothable Alexandrov space.

\bigskip

In the case of collapsing, the liftings $f_n$ of a function $f$ with controlled concavity
type do not have the same controlled concavity type.

Nevertheless, the liftings are semiconcave and moreover, as was noted by Kapovitch in
\cite{kapovitch:collapsing}, if
$M_n$ is a sequence of $(m+k)$-dimensional Riemannian manifolds with curvature $\ge
\kappa$, $M_n\GHto A$, $\dim A=m$, then one has a good control over the sum of
$k+1$ maximal eigenvalues of their Hessians. 
In particular, a construction as in the proof of theorem~\ref{thm:strictly-concave} gives a strictly concave function
on $A$ for which the liftings $f_n$ on $A_n$ have Morse index $\le k$.
It follows that one can retract an $\eps$-neighborhood of $p_n$ to a $k$-dimensional CW-complex%
\footnote{it is unknown whether it could be retracted to an $k$-submanifold. If true, it would give some interesting applications}%
, where $p_n\in A_n$ is a maximum point of $f_n$ and $\eps$ does not depend on $n$.
This observation gives a lower bound for the \emph{codimension of
a collapse}%
\footnote{in our case, it is $k$; the difference between the dimension of spaces from the collapsing sequence and the dimension of the limit space} 
to particular spaces. 
For example, for any lower curvature bound $\kappa$, the codimension of a collapse to $\Sigma(\mathbb H\mathrm P^m)$%
\footnote{that is, the spherical suspension over $\mathbb H\mathrm P^m$} 
is
at least 3, and for $\Sigma(C a\! \operatorname{P}^2)$  is at least 8 (it is expected to be $\infty$). 
In addition, it yields the following funny sphere theorem, 
it is funny since it does not assume positiveness of curvature.

\begin{thm}{\bf Funny sphere theorem.}
If a $4{\cdot}(m+1)$-dimensional Riemannian manifold $M$ with sectional curvature $\ge\kappa$ is
sufficiently close%
\footnote{that is, $\eps$-close for some $\eps=\eps(\kappa,m)$} 
to $\Sigma(\mathbb H\mathrm P^m)$, then it
is homeomorphic to a sphere.
\end{thm}

The controlled concavity also gives a short proof of the following result:

\begin{thm}{\bf Theorem.}\label{thm:unit-speed} Any quasigeodesic is a
unit-speed curve.
\end{thm}

\Proof.  To prove that a quasigeodesic $\gamma$ is $1$-Lipschitz at some $t=t_0$,
it is enough to apply the definition for $f=\dist_{\gamma(t_0)}^2$ and use the fact
that in any Alexandrov's space $\dist_p^2$ is $(2+O(r^2))$-concave in
$B_r(p)$.

Note that if $A_n,A\in\Alex^m(\kappa)$, $A_n\GHto A$ without collapse, and
$\gamma_n$ in $A_n$ is a sequence of quasigeodesics which converges to a curve
$\gamma$ in $A$, then $\gamma$ has the following property%
\footnote{from statement~\ref{thm:limit-qg}, page~\pageref{thm:limit-qg}%
, we that $\gamma$ is a quasigeodesic, but its proof is based on this theorem}:

\begin{thm}{\bf Property.}\label{propr:weak-qg}
For any function $f$ on $A$ with controlled concavity type $(\lambda,\kappa)$ we
have that $f\circ\gamma$ is $\lambda$-concave.
\end{thm}

If $\gamma$  is a quasigeodesic in $A$ with $\gamma(0)=p$, then the curves
$\gamma(t/\lam)$ are quasigeodesics in $\lam{\cdot} A$. 
Therefore, as $\lam\to\infty$, the limit curve
$$\gamma_\infty(t)=\left[\begin{matrix}
|t|\cdot\gamma^+(0)&\text{if}\ \  t\ge0\\
|t|\cdot\gamma^-(0)&\text{if}\ \  t<0\\
          \end{matrix}\right.$$
in $T_p$ has the above property.
By a construction similar%
\footnote{Setting $v=\gamma^\pm(0)\in T_p$ and
$w=2\gamma^\pm(0)$, this function can be presented as a sum
$$f=A(\phi_{r,c}\circ\dist_{o}+\phi_{r,c}\circ\dist_w)+B\sum_i\phi_{r',c'}
\circ\dist_{q_i},$$
for appropriately chosen positive reals 
$A,\ B,\ r,\ r',\ c,\ c'$ and a collection of points $q_i$ such that,
$\angle o p q_i=\tilde\angle_0 o p q_i=\tfrac\pi2$, $|p q_i|=r$ .} 
to theorem~\ref{thm:strictly-concave}, for any $\eps>0$ there is a function $f$ of controlled
concavity type $(-2+\eps,-\eps)$ on a neighborhood of $\gamma^\pm\in T_p$ such that
$$f(t\cdot\gamma^\pm)=-(t-1)^2+o((t-1)^2).$$
Applying the property above we get $|\gamma^\pm(0)|\ge 1$. 
\qeds

\begin{thm}{\bf Remark.}\label{rmk:stronger-unit-speed} {\rm Note that we have proven
a slightly stronger statement; namely, if a curve $\gamma$ satisfies the
property~\ref{propr:weak-qg} then it is a unit-speed curve}.
\end{thm}

\begin{thm}
{\bf Question.} Is it true that for any point $p\in A$ and any $\eps>0$, there is a
$(-2+\eps)$-concave function $f_p$ defined in a neighborhood of $p$, such that $f_p(p)=0$ and
$f_p\ge-\dist_p^2$?
\end{thm}

Existence of a such function would be a useful technical tool. 
In particular, it
would allow for an easier proof of the above theorem.

\section{Tight maps}
\setcounter{subsection}{1}
\setcounter{thm}{0}

The tight maps considered in this section give a more flexible version of
distance charts.

Similar maps (so called \emph{regular maps}) were used by Perelman in \cite{perelman:spaces2} and \cite{perelman:morse}.
Later he modified them to nearly this form in \cite{perelman:DC}. 
Recently Lytchak and Nagano used this technique  
for Alexandrov's spaces with upper curvature bound.

\begin{thm}{\bf Definition.}\label{def:tight} Let $A\in \Alex^m$ and $\Omega\i
A$ be an open subset.
A collection of semiconcave functions
$f_0,f_1,\dots,f_\ell$ on $A$
is called \emph{tight in $\Omega$} if 
$$\sup_{x\in \Omega,\,i\not=j}\{d_x f_i(\nabla_x f_j)\}<0.$$
In this case the map 
$$F\:\Omega\to \RR^{\ell+1},\ \ F\:x\mapsto (f_0(x),f_1(x),\dots,f_\ell(x))$$
is called \emph{tight}.

A point $x\in\Omega$ is called a \emph{critical point} of $F$ if $\min_i\{d_x f_i\}\le
0$, otherwise the point $x$ is called \emph{regular}.
\end{thm}

\begin{thm}{\bf Main example.} If $A\in\Alex^m(\kappa)$ and
$a_0,a_1,\dots,a_\ell,p\in A$ such that 
$$\tilde\angle_\kappa a_i p a_j>\tfrac\pi2\ \ \text{for all}\ \  i\not=j$$ 
then the map
$x\mapsto (|a_0x|,|a_1x|,\dots,|a_\ell x|)$ is tight in a neighborhood of
$p$.
\end{thm}

The inequality in the definition follows from inequality $(**)$ 
on page \pageref{**-polar-inq} and a subsequent to it example~(\ref{polar}).

This example can be made slightly more general.
Let $f_0,f_1,...,f_\ell$ be a collection of simple functions
$$f_i=\Theta_i(\dist_{a_{1,i}}^2,\dist_{a_{2,i}x}^2,\dots, \dist_{a_{n_i,i}x}^2)$$
and the sets of points $K_i=\{a_{k,i}\}$ satisfy the following inequality
$$\tilde\angle_\kappa x p y>\tfrac\pi2\ \ \text{for any}\ \  x\in K_i,\ \  y\in K_j,\ \ i\not=j.$$
Then the map
$x\mapsto(f_0(x),f_1(x),...,f_\ell(x))$
is tight in a neighborhood of $p$.
 We will call such a map a \emph{simple tight map}.

Yet further generalization is given in the property~\ref{lip+} below.

The maps described in this example have an important property, they are liftable
and their lifts are tight. 
Namely, given a converging sequence $A_n\GHto A$, $A_n\in\Alex^m(\kappa)$ and a
simple tight map $F\:A\to \RR^{\ell+1}$ around $p\in A$, the construction in
section~\ref{adm} gives simple tight maps $F_n\:A_n\to \RR^\ell$ for large $n$, $F_n\to F$.

I was unable to prove that tightness is a stable property in a sense
formulated in the question below.
It is not really important for the theory since all maps which appear naturally 
are simple (or, in the worst case they are as in the generalization and
as in the property \ref{lip+}).
However, for the beauty of the theory it would be nice to have an answer to
the following question.

\begin{thm}{\bf Question.} Assume $A_n\GHto A$, $A_n\in\Alex^m(\kappa)$,
$f,g\:A\to \RR$ is a tight collection around $p$ and $f_n,g_n\:A_n\to\RR$, $f_n\to
f$, $g_n\to g$ are two sequences of $\lambda$-concave
functions and $A_n\ni p_n\to p\in A$. 
Is it true that for all large $n$, the collection $f_n,g_n$ must be
tight around $p_n$?

If not, can one modify the definition of tightness so that
\begin{enumerate}[(i)]
\item it would be stable in the above sense,
\item the definition would make sense for all semiconcave functions
\item the maps described in the main example above are tight?
\end{enumerate}
\end{thm}

Let us list some properties of tight maps with sketches of proofs:
\begin{enumerate}
\item\label{lip+} 
Let $x\mapsto(f_0(x),f_1(x),...,f_\ell(x))$ be a tight map in an open subset $\Omega\i A$, then there is $\eps>0$ such that if $g_0,g_1,...,g_n$ is a collection of
$\epsilon$-Lipschitz semiconcave functions in $\Omega$ then the map 
$$x\mapsto (f_0(x)+g_0(x),f_1(x)+g_1(x),...,f_\ell(x)+g_\ell(x))$$
is also tight in $\Omega$.
\item \label{open-reg}The set of regular points of a tight map is open. 

Indeed, let $x\in\Omega$ be a regular point of tight map $F=(f_0,f_1,\dots,f_\ell)$. 
Take real $\lambda$ so that $f_i''\le \lambda$ for all $i$ in a fixed neighborhood of
$x$.
Take a point $p$ sufficiently close to $x$ such that $d_x f_i(\uparrow_x^p)>0$
and moreover $f_i(p)-f_i(x)>\tfrac\lambda2{\cdot}
|x p|^2$ for each $i$.
Then, from $\lambda$-concavity of $f_i$, there is a small neighborhood
$\Omega_x\ni x$ such that for any $y\in\Omega_x$ and $i$ we have $\d_y
f_i(\uparrow_y^p)\ge\eps$
for some fixed $\eps>0$.

\item \label{minus-funct}If one removes one function from a tight
collection (in $\Omega$) then (for the corresponding map) all points of $\Omega$ become
regular. 
In other words, the projection of a tight map $F$ to any coordinate
hyperplane is a tight map with all regular points (in $\Omega$).

This follows from the property \ref{grad-onto} on page \pageref{grad-onto}
applied to the flow for the removed function $f_i$.

\item \label{extra-funct} 
The converse also holds;
that is,
if $F$ is regular at $x$ then
one can find a semiconcave function $g$ such that map $z\mapsto(F(z),g(z))$ is
tight in a neighborhood of $x$.
Moreover, $g$ can be chosen to have an arbitrary controlled concavity type.

Indeed, one can take $g=\dist_p$, where $p$ as in the property~\ref{open-reg}.
Then we have 
$$d_x g(v)=-\max_{\xi\in\Uparrow_x^p}\{\<\xi,v\>\}$$
and therefore
$$d_x g(\nabla_x f_i)
=
-\max_{\xi\in\Uparrow_x^p}\{\<\xi,\nabla_x f_i\>\}
\le
-\max_{\xi\in\Uparrow_x^p} \{d_x f(\xi)\}
\le
-\eps.$$
On the other hand, from inequality $(**)$  on page \pageref{**-polar-inq} and
example (\ref{polar}) subsequent to it, we have
$$d_x f_i(\nabla_x g)+\min_{\xi\in\Uparrow_x^p}\{d_x f_i(\xi)\}\le 0.$$
The last statement follows from the construction in
theorem~\ref{thm:strictly-concave}. 

\item\label{co_Lipschitz} 
A tight map is open and even \emph{co-Lipschitz}%
\footnote{A map $F\: X\to Y$ between
metric spaces is called $L$-co-Lipschitz in $\Omega\i X$ if for any ball $B_r(x)\i
\Omega$ we have $F(B_r(x))\supset B_{r/L}(F(x))$ in $Y$} 
in a neighborhood of any regular point.

This follows from lemma~\ref{lem:tight-dir}.

\item Let $A\in\Alex$, $\Omega\i A$ be an open subset.
If $F\:\Omega\to\RR^{\ell+1}$ is tight then $\ell\le\dim A$.

Follows from the properties~\ref{minus-funct} and \ref{co_Lipschitz}.

\item{\it Morse lemma.}\label{morse} A tight map admits a local splitting in a neighborhood of its regular point, and a proper everywhere regular tight map is a locally trivial fiber bundle. Namely

\begin{enumerate}[(i)]
\item If $F\:\Omega\to\RR^{\ell+1}$ is a tight map and $p\in \Omega$ is a regular point, then there is a neighborhood $\Omega\supset\Omega_p\ni p$ and homeomorphism 
$$h\:\Upsilon\times F(\Omega_p)\to \Omega_p,$$ 
such that
$F\circ h$ coincides with the projection to the second coordinate $\Upsilon\times F(\Omega_p)\to F(\Omega_p)$.
\item If $F\:\Omega\to\Delta\i\RR^{\ell+1}$ is a proper tight map and all points in $\Delta\z\i \RR^{\ell+1}$ are regular values of $F$, then $F$ is a locally trivial fiber bundle.
\end{enumerate}
The proof is a backward induction on $\ell$, see 
\cite[1.4]{perelman:morse}, \cite[1.4.1]{perelman:spaces2} or
\cite[6.7]{kapovitch:stability}.
\end{enumerate}

The following lemma is an analog of lemmas 
\cite[2.3]{perelman:morse} and \cite[2.2]{perelman:DC}.

\begin{thm}{\bf Lemma.}\label{lem:tight-dir} Let $x$ be a regular point of a tight
map $$F\:x\mapsto(f_0(x),f_1(x),\dots,f_\ell(x)).$$
Then there is $\eps>0$ and a neighborhood $\Omega_x\ni x$ such that for any
$y\in \Omega_x$ and $i\in\{0,1,\dots,\ell\}$ there is a unit vector $w_i\in \Sigma_x$
such that $\d_x f_i(w_i)\ge\eps$ and $d_x f_j(w_i)=0$ for all $j\not=i$. 

Moreover, if $E\i A$ is an extremal subset and $y\in E$ then $w_i$ can be
chosen in $\Sigma_y E$.

\end{thm}

\Proof. Take $p$ as in the property \ref{open-reg} page~\pageref{open-reg}.
Then we can find a neighborhood $\Omega_x\ni x$ and $\eps>0$ so that for any
$y\in \Omega_x$
\begin{enumerate}[(i)]
\item $d_y f_i(\uparrow_y^p)>\eps$ for each $i$;
\item $-d_y f_i(\nabla_y f_j)>\eps.$ for all $i\not=j$.
\end{enumerate}

Note that if $\alpha(t)$ is an $f_i$-gradient curve in $\Omega_x$ then 
$$(f_i\circ\alpha)^+>0
\ \ \text{and}\ \ 
(f_j\circ\alpha)^+\le -\eps
\ \ \text{for any}\ \ 
j\not=i.$$
Applying  lemma~\ref{lem:stable-grad-curves} for
$(\lam{\cdot} A,y)\GHto T_y$, $\lam\cdot[f_i-f_i(y)]\to \d_y f_i$, we get the same inequalities for 
$\d_y f_i$-gradient curves on $T_y$.
That is, if $\beta(t)$ is an 
$\d_y f_i$-gradient curve in $T_y$ then 
$$(\d_y f_i\circ\beta)^+>0
\ \ \text{and}\ \ 
(\d_y f_j\circ\beta)^+\le-\eps
\ \
\text{for any}\ \ 
j\not=i.$$
Moreover, $\d_y f_i(v)>0$
implies $\<\nabla_v(\d_y f_i),\uparrow_{v}^o\><0$, 
therefore in this case $|\beta(t)|^+\z>0$.

Take $w_0\in T_y$ to be a maximum point for $d_y f_0$ on the set 
$$\{v\in T_y|f_i(v)\ge0, |v|\le 1\}.$$
Then
$$d_y f_0(w_0)\ge d_y f_0(\uparrow_y^p)>\eps.$$
Assume for some $j\not=0$ we have $f_j(w_0)>0$.
Then 
$$\min_{i\not=j} \{d_{w_0}d_y f_i,d_{w_0}\nu\}\le 0,$$
where the function $\nu$ is defined by $\nu\:v\mapsto -|v|$; 
this is a concave function
on $T_y$.
Therefore, if $\beta_j(t)$ is a $d_y f_j$-gradient curve with an end%
\footnote{it does exist by property \ref{grad-onto} on page \pageref{grad-onto}} 
point at
$w_0$, then moving along $\beta_j$ from $w_0$ backwards decreases only
$d_y f_j$,  and increases the other $d_y f_i$ and $\nu$ in the first order; this is a
contradiction.

To prove the last statement it is enough to show that $w_0\in T_y E$, which follows since $T_y E\i T_y$ is an extremal subset (see property~\ref{ext-tangent} on page~\pageref{ext-tangent}).  
\qeds

\begin{thm}{\bf Main theorem.}\label{thm:tight-map}
Let $A\in\Alex^m(\kappa)$, $\Omega\i A$ be the interior of a compact convex subset,
and 
$$F\:\Omega\to \RR^{\ell+1},\ \ F\:x\mapsto(f_0(x),f_1(x),\dots,f_\ell(x))$$ 
be a tight map. 
Assume all $f_i$ are strictly concave.
Then 
\begin{enumerate}[(i)]
\item the set of critical points of $F$ in $\Omega$ forms an $\ell$-submanifold $M$ 
\item  $F\:M\to\RR^{\ell+1}$ is an embedding. 
\item  $F(M)\i \RR^{\ell+1}$ is a convex hypersurface which lies in the boundary of
$F(\Omega)$%
\footnote{In fact $F(M)=\partial F(\Omega)\cap F(\Omega)$.}%
.
\end{enumerate}
\end{thm}

\begin{thm}{\bf Remark.} 
{\rm The condition that all $f_i$ are strictly concave
seems to be very restrictive, but that is not really so; if $x$ is a regular
point of a tight map $F$ then, using properties~\ref{lip+} and~\ref{extra-funct}
on page~\pageref{lip+}, one can find $\eps>0$ and $g$ such that 
$$F'\:y
\mapsto
(f_0(y)+\eps g(y),\dots,f_\ell(y)+\eps g(y),g(y))$$ 
is tight in a small neighborhood of $x$ and all its coordinate functions are
strictly concave. 
In particular, in a neighborhood of $x$ we have 
$$F=L\circ F'$$
where $L\:\RR^{\ell+2}\to \RR^{\ell+1}$ is linear.}
\end{thm}

\begin{thm}{\bf Corollary.}\label{cor:conv-chart}
In the assumptions of theorem~\ref{thm:tight-map}, if in addition $m=\ell$ then
$M=\Omega$, $F(\Omega)$ is a convex hypersurface in $\RR^{m+1}$ and $F\:\Omega\to
\RR^{m+1}$ is a locally bi-Lipschitz embedding. 
Moreover, each projection of $F$ to a coordinate hyperplane is a locally
bi-Lipschitz homeomorphism.
\end{thm}

\Proof\ {\it of theorem~\ref{thm:tight-map}}. 
Let $\gamma\:[0,s]\to A$ be a minimal unit-speed geodesic connecting $x,y\in \Omega$, so $s=|x y|$.
Consider a straight segment $\bar\gamma$ connecting $F(x)$ and $F(y)$: 
$$\bar\gamma\:[0,s]\to\RR^{\ell+1},
\ \  
\bar\gamma(t)=F(x)+\tfrac
t s\cdot\l[F(y)-F(x)\r].$$
Each function $f_i\circ\gamma$ is concave, therefore all coordinates of 
$$F\circ\gamma(t)-\bar\gamma(t)$$ 
are non-negative.
This implies that the Minkowski sum%
\footnote{equivalently 
$Q=\{(x_0,x_1,\dots,x_\ell)\in\RR^{\ell+1}|\exists(y_0,y_1,\dots,y_\ell)\in
F(\Omega)\forall i\  x_i\le y_i\}$.} 
$$Q=F(\Omega)+(\RR_-)^{\ell+1}$$ 
is a convex set.

Let $x_0\in \Omega$ be a critical point of $F$. 
Since $\min_i\{\d_{x_0}f_i\}\le 0$, at least one of coordinates
of $F(x)$ is smaller than the corresponding coordinate of $F(x_0)$ for any $x\in\Omega$.
In particular, $F$ sends its critical point to the boundary of $Q$.

Consider map 
$$G\:\RR^{\ell+1}\to A,
\ \ 
G\:(y_0,y_1,\dots,y_\ell)
\mapsto
\argmax\{\min_i\{f_i-y_i\}\}$$
where \label{argmax}$\argmax\{f\}$ denotes a maximum point of $f$.
The function $\min_i\{f_i-y_i\}$ is strictly concave; 
therefore
$\argmax\{\min_i\{f_i-y_i\}\}$ is uniquely defined and $G$ is continuous in the
domain of definition.%
\footnote{We do not need it, but clearly
$$G(y_0,y_1,\dots,y_\ell)
=
G(y_0+h,y_1+h,\dots,y_\ell+h)$$ 
for any $h\in \RR$.} 
The image of $G$ coincides with the set of critical points of $F$ and moreover $G\circ
F|_M=\id_M$. 
Therefore $F|_M$ is a homeomorphism%
\footnote{In general, $G$ is not Lipschitz
(even on $F(M)$); 
even in the case  $f_i''\le-1)$ for all $i$, 
it is only
possible to prove that $G$ is H\"older continuous of class $C^{0;\frac12}$. (In
fact the statement in \cite{perelman:spaces2}, page 20, lines 23--25 is wrong
but the proposition 3.5 is still OK.)}.
\qeds

\Proof\ {\it of corollary~\ref{cor:conv-chart}}. It only remains to show that
$F$ is locally bi-Lipschitz. 

Note that for any point $x\in \Omega$, one can find $\eps>0$ and a
neighborhood $\Omega_x\ni x$, so that for any direction $\xi\in \Sigma_y$, $y\in\Omega_x$ one can choose $f_i$, $i\z\in \{0,1,\dots,m\}$, such that $d_x f_i(\xi)\le -\eps$. 
Otherwise, by a slight perturbation%
\footnote{as in the property~\ref{lip+} on page~\pageref{lip+}} of collection $\{f_i\}$ we get a map $F\:A^m\to \RR^{m+1}$ regular at $y$, which contradicts property~\ref{co_Lipschitz}.

Therefore applying it for $\xi=\uparrow_z^y$ and $\uparrow_y^z$, $z,y\in\Omega$, we get two values $i,j$ such that 
$$f_i(y)-f_i(z)
\ge 
\eps{\cdot}|y z|
\ \ \text{and}\ \ 
f_j(z)-f_j(y)\ge \eps{\cdot}|y z|.
$$
Therefore $F$ is bi-Lipschitz.

Clearly $i\not=j$ and therefore at least one of them is not zero.
Hence the projection map $F'\:x\mapsto(f_1(x),\dots,f_m(x))$ is also locally
bi-Lipschitz.
\qeds

\subsection{Applications.}
\label{app-tight} 

One series of applications of tight maps is Morse theory 
for Alexandrov's spaces, it is based on the
main theorem~\ref{thm:tight-map}. 
It includes Morse lemma (property~\ref{morse} page~\pageref{morse}) and

\begin{enumerate}[$\diamond$]

\item{\it Local structure theorem \cite{perelman:morse}.}
Any small spherical neighborhood of a point in an Alexandrov's space is homeomorphic to a cone over its boundary.

\item{\it Stability theorem \cite{perelman:spaces2}.} 
For any compact $A\in \Alex^m(\kappa)$ there is $\eps>0$ such that if $A'\in \Alex^m(\kappa)$ is $\eps$-close to $A$ then $A$ and $A'$ are homeomorphic.

\end{enumerate}
The other series is the regularity results on an Alexandrov's space. 
These results obtained by Perelman in \cite{perelman:DC} 
are improvements of earlier results of Otsu and Shioya in \cite{otsu-shioya} and \cite{otsu:second-der}. 
It use mainly
the corollary~\ref{cor:conv-chart} and the smoothing trick; see subsection~\ref{smooth}. 
\begin{enumerate}[$\diamond$]
\item Components of metric tensor of an Alexandrov's space in a chart are continuous
at each regular point%
\footnote{that is, at each point with Euclidean tangent space}. 
Moreover they have bounded variation and are differentiable almost everywhere.
\item The Christoffel symbols  in a chart are well defined as signed Radon measures.
\item Hessian of a semiconcave function on an Alexandrov's space is defined almost
everywhere. 
That is, if $f\:\Omega\to\RR$ is a semiconcave function, then for almost
any $x_0\in \Omega$ there is a symmetric bi-linear form $\Hess_f$ such that
$$f(x)
=
f(x_0)+d_{x_0}f(v)+\Hess_f(v,v)+o(|v|^2),$$
where $v=\log_{x_0}x$. 
Moreover, $\Hess_f$ can be calculated using  standard formulas in the above
chart.
\end{enumerate}

Here is yet another, completely Riemannian application. This statement has been
proven by Perelman, a sketch of its proof is included in an appendix to
\cite{petrunin:PL}.
The proof is based on the following observation: if $\Omega$ is an open subset of a Riemannian manifold and  $F\:\Omega\to\RR^{\ell+1}$ is a tight map with strictly concave coordinate functions, then its level sets $F^{-1}(x)$ inherit the lower curvature bound. 
\begin{enumerate}[$\diamond$]
\item \emph{Continuity of the integral of scalar curvature.} 
Given  a compact Riemannian manifold $M$, let us define $\mathcal F(M)=\int_M\Sc$. 
Then $\mathcal F$ is continuous on the space of Riemannian
$m$-dimensional manifolds with uniform lower curvature and upper diameter
bounds.%
\footnote{In fact $\mathcal F$ is also bounded on the set of Riemannian
$m$-dimensional manifolds with uniform lower curvature, this is proved in
\cite{petrunin:scalar} by a similar method.}
\end{enumerate}

\section{Please deform an Alexandrov's space.}
\setcounter{subsection}{1}

In this section we discuss a number of related open problems. 
They seem to be very hard, but I think it is worth to
write them down in order to indicate the border between known and unknown things.

The main problem in Alexandrov's geometry is to find a way to vary Alexandrov's
space, or simply to find a nearby Alexandrov's space to a given Alexandrov's
space. 
Lack of such variation procedure makes it impossible to use Alexandrov's
geometry in the way it was designed to be used: 

For example, assume you want to solve the Hopf conjecture;
that is, you want to find out
if $S^2\times S^2$ carries a metric with positive sectional curvature. 
Assuming such metrics exist, 
there is a volume maximizing Alexandrov's metrics $d$
on $S^2\times S^2$ with curvature $\ge 1$.
(There is no reason to believe
that this metric $d$ is Riemannian, but from Gromov's compactness theorem such
Alexandrov's metric should exist.)
Provided we have a procedure to vary $d$ while keeping its curvature $\ge 1$,
we could find some special properties of $d$  and in ideal situation we could arrive to a contradiction.

At the moment, except for boring rescaling, there is no variation procedure available.
The following conjecture (if true) would give such a procedure.
Although it will not be sufficient to solve the Hopf conjecture, it will give some information about the critical Alexandrov's metric.

\begin{thm}{\bf Conjecture.}\label{conj:bry}
The boundary of an Alexandrov's space equipped with induced intrinsic metric is an
Alexandrov's space with the same lower curvature bound.
\end{thm}

This also can be reformulated as:

\bigskip

\noi{\bf \ref{conj:bry}$'\!\!\!$. Conjecture.}{\it\ Let $A$ be an Alexandrov space without boundary. Then a convex hypersurface in $A$ equipped with induced intrinsic metric is an Alexandrov's space with the same lower curvature bound.}

\bigskip

This conjecture, if true, would give a variation procedure. 
For example if $A$ is a non-negatively curved Alexandrov's space and $f\:A\to\RR$ is concave (so $A$ is necessarily open) then for any $t$ the graph
$$A_t
=
\{(x,t{\cdot} f(x))\in A\times\RR\}$$
with induced intrinsic metric would be an Alexandrov's space. 
Clearly $A_t\GHto A$ as $t\to0$. 
An analogous construction exists for semiconcave functions on closed manifolds, but
one has to take a \emph{parabolic cone} 
(see footnote~\ref{par-cone} on page~\pageref{par-cone}) 
instead of the product.

So far, even for a convex hypersurface in a Riemannian manifold, there is only one proof available (it is given by Alexander, Kapovitch and me in \cite{akp} 
and based on the result of Buyalo in \cite{buyalo:convex-surface}) 
which uses smoothing and the Gauss formula. 
There is one synthetic proof by Milka
(see \cite{milka-conv}) for a convex surface in the Euclidean space, but this proof heavily relies on Euclidean
structure and it seems impossible to generalize it even to the Riemannian case.

There is a chance of attacking this problem by proving a type of the Gauss formula for
Alexandrov's spaces. 
One has to start with defining a curvature tensor of Alexandrov's spaces (it
should be a measure-valued tensor field), then prove that the constructed tensor is really responsible for the geometry of the space. 
Such things were already done in the two-dimensional case 
and for the spaces with two-sided curvature bound, 
see \cite{reshetnyak:curvature} and
\cite{nikolaev:curvature} respectively.
So far the best results in this direction are given by Perelman in \cite{perelman:DC}, see also section~\ref{app-tight} for more details.
This approach, if works, would give something really new in the area.

Almost everything that is known so far about the intrinsic metric of a boundary is also known for the intrinsic metric of a general extremal subset.
In \cite{perelman-petrunin:extremal}, it was conjectured  that an analog of conjecture~\ref{conj:bry} is true for any \emph{primitive extremal subset}, but it turned out to be wrong; a simple example was constructed in \cite{petrunin:extremal}. 
All such examples appear when codimension of extremal subset is $\ge 3$.

\begin{thm}{\bf Conjecture.}\label{conj:codim=2}
Let $A\in\Alex(\kappa)$, $E\i A$ be a primitive extremal subset and $\codim
E=2$ then $E$ equipped with induced intrinsic metric belongs to  $\Alex(\kappa)$
\end{thm}

The following question is closely related to conjecture~\ref{conj:bry}.

\begin{thm}{\bf Question.}\label{qst:lift-conc}
Assume $A_n\GHto A$, $A_n\in \Alex^m(\kappa)$, $\dim A=m$ (that is, 
it is not a collapse).

Let $f$ be a $\lambda$-concave function of an Alexandrov's space $A$. 
Is it always possible to find a sequence of $\lambda$-concave functions
$f_n\:A_n\to \RR$ which converges to $f\:A\to\RR$? 
\end{thm}

Here is an equivalent formulation:

\bigskip

\noi{\bf \ref{qst:lift-conc}.$'$}{\bf \,Question.} {\it Assume $A_n\GHto A$, $A_n\in \Alex^m(\kappa)$, $\dim
A=m$ (that is, it is not a collapse) and $\partial A=\emptyset$.

Let $S\i A$ be a convex hypersurface. 
Is it always possible to find a sequence of convex hypersurfaces $S_n\i A_n$
which converges to $S$?
}

\bigskip

If true, this would give a proof of  conjecture~\ref{conj:bry} for the
case of a \emph{smoothable Alexandrov's space} (see page~\pageref{smoothable}).

In most of (possible) applications, Alexandrov's spaces appear as
limits of Riemannian manifolds of the same dimension.
Therefore, even in this reduced generality, 
an answer is valuable.

The question of whether an Alexandrov space is smoothable is also far from being solved.
From Perelman’s stability theorem, if an Alexandrov's space has topological singularities then it is not smoothable.
Moreover, by a theorem of Kapovitch in \cite{kapovitch:regularity},
one has that any space of
directions of a smoothable Alexandrov's space is homeomorphic to the sphere. 
Except for the 2-dimensional case, it is only known that any polyhedral metric of non-negative curvature on a 3-manifold is smoothable (see \cite{matveev:smooth}).
There is yet no procedure of smoothing an Alexandrov's space even in a neighborhood of a regular point.

Maybe a more interesting question is whether smoothing is unique up to a diffeomorphism. 
If the answer is positive it would imply in particular that any Riemannian
manifold with curvature $\ge 1$ and $\diam>\tfrac\pi2$ is diffeomorphic(!) to the
standard sphere, see \cite{grove:grove's-question} for details.
Again, from Perelman's stability theorem (\cite{perelman:spaces2}), it follows that any two smoothings must be homeomorphic. 
In fact it seems likely that any two smoothings are PL-homeomorphic; see
\cite[question 1.3]{kapovitch:stability} and discussion right before it.
It seems that today there is no technique which might approach the general uniqueness 
problem (so maybe one should try to construct a counterexample).

One may also ask similar questions in the collapsing case.
In \cite{spher-susp-W},
Petersen and Wilhelm constructed Alexandrov's spaces with curvature $\ge1$ which can not be presented as a limit of an (even collapsing) sequence of Riemannian manifolds with curvature
$\ge \kappa>\tfrac14$. 
In \cite{kapovitch:collapsing}, Kapovitch found some lower bounds for codimension of collapse with arbitrary lower curvature bound to some special Alexandrov's spaces, see section~\ref{app-con-con} for more discussion.
It is expected that the same spaces (for example, the spherical suspension over the Cayley plane) can not be approximated by sequence of
Riemannian manifolds of any fixed dimension and any fixed lower curvature bound, but so far this question remains open.

\appendix\section{Existence of quasigeodesics}
\label{constr-qg}
\setcounter{thm}{0}

This appendix is devoted to the proof of property~\ref{exist-qg} on
page~\pageref{exist-qg}. 
Namely we prove the following.

\begin{thm}\label{thm:exist-qg}{\bf Existence theorem.}
Let $A\in\Alex^m$, then for any point $x\in A$, and any direction $\xi\in \Sigma_x$
there is a quasigeodesic $\gamma\:\RR\to A$ such that $\gamma(0)=x$ and
$\gamma^+(0)=\xi$.

Moreover if $E\i A$ is an extremal subset and $x\in E$, $\xi\in \Sigma_x E$
then $\gamma$ can be chosen to lie completely in $E$.
\end{thm}
The proof is quite long; it was obtained by Perelman around 1992;
here we present a simplified proof similar to \cite{perelman-petrunin:qg} which is
based on the gradient flow technique. 
We include a complete proof here, since otherwise it would never be published.

Quasigeodesics will be constructed in three big steps. 
\begin{description}
\item [\ref{step1}] Monotonic curves $\longrightarrow$ convex curves.
\item [\ref{step2}] Convex curves $\longrightarrow$ pre-quasigeodesics.
\item [\ref{step3}] Pre-quasigeodesics $\longrightarrow$ quasigeodesics.
\end{description}

In each step, we construct a better type of curves from a given type of curves
by an extending-and-chopping procedure and then passing to a limit.
The last part is most complicated.

The second part of the theorem \ref{thm:exist-qg} 
is
proved in the subsection~\ref{qg-extrim}.

\setcounter{subsection}{-1}
\subsection{Step 0: Monotonic curves} 

As a starting point we use radial curves, which do exist for any initial data
(see section~\ref{gexp}), and by lemma~\ref{lem:monotonic} 
they are monotonic in the
sense of the following definition:

\begin{thm}{\bf Definition.}
A curve $\alpha(t)$ in an Alexandrov's space $A$ is called \emph{monotonic} with
respect to a parameter value $t_0$ if for any $\lambda$-concave function $f$,
$\lambda\ge 0$, we have that function
$$t
\mapsto
\frac{f\circ\alpha(t+t_0)-f\circ\alpha(t_0)-\tfrac\lambda2{\cdot}t^2}t$$ 
is non-increasing for $t>0$. 
\end{thm}

The following construction makes a new monotonic curve out of two. 
It will be used in the next section to construct \emph{convex curves}.

\begin{thm}{\bf Extension.}\label{ext-mono}
Let $A\in\Alex$, $\alpha_1\:[a,\infty)\to A$ and $\alpha_2\:[b,\infty)\to A$ be two
monotonic curves with respect to $a$ and $b$ respectively. 

Assume 
$$a\le b,
\ \ \alpha_1(b)=\alpha_2(b)
\ \ \text{and}\ \
\alpha^+_1(b)=\alpha^+_2(b).$$ 
Then its joint
$$\beta\:[a,\infty)\to A,
\ \ 
\beta(t)=\l[
\begin{matrix}
\alpha_1(t)&\text{if}&t< b\\
\alpha_2(t)&\text{if}&t\ge b
\end{matrix}\right.
$$
is monotonic with respect to $a$ and $b$.
\end{thm}

\Proof.
It is enough to show that 
$$t 
\mapsto 
\frac{f\circ\alpha_2(t+a)-f\circ\alpha_1(a)-\tfrac\lambda2{\cdot}t^2}t
$$
is non-increasing for $t\ge b-a$.
By simple algebra, it follows from the following two facts:
\begin{enumerate}[$\diamond$]
\item $\alpha_2$ is monotonic and therefore
$$t 
\mapsto 
\frac{f\circ\alpha_2(t+b)-f\circ\alpha_2(b)-\tfrac\lambda2{\cdot}t^2}t$$
is non-increasing for $t>0$.
\item From monotonicity of $\alpha_1$,
\begin{align*}
(f\circ\alpha_2)^+(b)
&=
d_{\alpha_1(b)}
f(\alpha_1^+(b))
=
\\
&=
(f\circ\alpha_1)^+(b)
\le
\\
&\le\frac{f\circ\alpha_1(b)+f\circ\alpha_1(a)-\tfrac\lambda2 (b-a)^2}{b-a}.
\end{align*}

\end{enumerate}
\qedsf

\subsection{Step 1: Convex curves.} \label{step1}

In this step we construct \emph{convex curves} with arbitrary initial data.

\begin{thm}{\bf Definition.}
A curve $\beta\:[0,\infty)\to A$ is called \emph{convex} if for any
$\lambda$-concave function $f$, $\lambda\ge 0$, we have that function $$t\mapsto
f\circ\beta(t)-\tfrac\lambda2{\cdot}t^2$$
is concave. 
\end{thm}

\noi{\bf Properties of convex curves.} 
Convex curves have the following properties; the proofs are either trivial or the
same as for quasigeodesics:

\begin{enumerate}
\item\label{conv-mono} A curve is convex if and only if it is monotonic with
respect to any value of parameter.
\item\label{conv-lip} Convex curves are $1$-Lipschitz.
\item\label{conv-tang} Convex curves have uniquely defined right and left
tangent vectors.
\item\label{limit-convex} A limit of convex curves is convex and the natural parameter converges to the natural parameter of the limit curves (the proof the last statement is based on the same idea as theorem~\ref{thm:unit-speed}).
\end{enumerate}

The next is a construction similar to \ref{ext-mono} which gives a new convex
curve out of two. 
It will be used in the next section to construct \emph{pre-quasigeodesics}.

\begin{thm}{\bf Extension.}\label{ext-conv}
Let $A\in\Alex$, $\beta_1\:[a,\infty)\to A$ and $\beta_2\:[b,\infty)\to A$ be two
convex curves. 
Assume $$a\le b,\ \ \beta_1(b)=\beta_2(b)\ \ \text{and}\ \ 
\beta^+_1(b)=\beta^+_2(b)$$ 
then its joint
$$\gamma\:[a,\infty)\to A,\ \ \gamma(t)=\l[
\begin{matrix}
\beta_1(t)&\text{if}&t\le b\\
\beta_2(t)&\text{if}&t\ge b
\end{matrix}\right.$$
is a convex curve.
\end{thm}

\Proof. Follows immediately from \ref{ext-mono} and property~\ref{conv-mono} above.

\begin{thm}{\bf Existence.} 
Let $A\in\Alex$, $x\in A$ and $\xi\in \Sigma_x$. 
Then there is a convex curve $\beta_\xi\:[0,\infty)\to A$ such that
$\beta_\xi(0)=x$ and $\beta_\xi^+(0)=\xi$.
\end{thm}

\Proof.
For $v\in T_x A$, consider the radial curve
$$\alpha_v(t)=\gexp_x(t v)$$
According to lemma~\ref{lem:monotonic} if $|v|=1$ then $\alpha_v$ is
$1$-Lipschitz and monotonic. 
Moreover, straightforward calculations show that the same is true for $|v|\le 1$.

Fix $\eps>0$. Given a direction $\xi\in \Sigma_x$, let us consider the following recursively
defined sequence of radial curves $\alpha_{v_n}(t)$ such that $v_0=\xi$ and
$v_{n}=\alpha^+_{v_{n-1}}(\eps)$.
Then consider their joint
$$\beta_{\xi,\eps}(t)=\alpha_{v_{\lfloor t/\eps\rfloor}}(t-\eps\lfloor
t/\eps\rfloor).$$
Applying an extension procedure~\ref{ext-mono} we get that
$\beta_{\xi,\eps}\:[0,\infty)\to A$ is monotonic with respect to any $t=n{\cdot}\eps$.

By property~\ref{conv-mono} on page~\pageref{conv-mono}, passing to a partial
limit $\beta_{\xi,\eps}\to\beta_{\xi}$ as $\eps\to 0$ we get a convex curve
$\beta_{\xi}\:[0,\infty)\to A$.

It remains to show that $\beta_{\xi}^+(0)=\xi$.

Since $\beta_\xi$ is convex, its right tangent vector is well defined and
$|\beta^+_\xi(0)|\z\le 1$%
\footnote{see properties~\ref{conv-tang}
and~\ref{conv-lip}, page~\pageref{conv-tang}}.
On the other hand, since $\beta_{\xi,\eps}$ are monotonic with respect to $0$, for any semiconcave function $f$
we have 
\begin{align*}
d_x f(\beta^+_{\xi}(0))
&=
(f\circ\beta_{\xi})^+(0)\le
\\
&\le
\lim_{\eps_i\to0}(f\circ\beta_{\xi,\eps})^+(0)=
\\
&=d_x f(\xi). \end{align*}
Substituting in this inequality $f=\dist_{y}$ with
$\angle(\uparrow_x^{y},\xi)<\eps$, we get $$\<\beta^+_\xi(0),\uparrow_x^{y}\>  >
 1-\eps$$
for any $\eps>0$. 
Together with $|\beta^+_\xi(0)|\le 1$ (property \ref{conv-lip} on page
\pageref{conv-lip}), it implies that 
$$\beta^+(0)=\xi.$$
\qedsf

\subsection{Step 2: Pre-quasigeodesics} \label{step2}

In this step we construct a
\emph{pre-quasigeodesic} with arbitrary initial data.

\begin{thm}{\bf Definition.} 
A convex curve $\gamma\:[a,b) \rightarrow A$ is called a pre-quasigeodesic if for
any $s\in [a,b)$ such that ${|\gamma^+(s)|}>0$, the curve $\gamma^s$ defined by
$$\gamma^s(t)=\gamma\left(s+\frac{t}{|\gamma^+(s)|}\right)$$
is convex for $t\ge0$, and if ${|\gamma^+(s)|}=0$ then $\gamma(t)=\gamma(s)$ for
all $t\ge s$.
\end{thm}

Let us first define entropy of pre-quasigeodesic, which measures ``how far'' a given
pre-quasigeodesic is from being a quasigeodesic.

\begin{thm}{\bf Definition.}\label{def:entropy}
Let $\gamma$ be a pre-quasigeodesic in an Alexandrov's space.

The \emph{entropy} of $\gamma$, $\mu_\gamma$ is the measure on the set of parameters
defined by
$$ \mu_\gamma ((a,b))=\ln |\gamma^+(a)|-\ln |\gamma^-(b)|.$$
\end{thm}
\noi Here are its main properties:
\begin{enumerate}
\item The entropy of a pre-quasigeodesic $\gamma$ is zero if and only if $\gamma$
is a quasigeodesic. 
\item\label{lim-entropy} For a converging sequence of pre-quasigeodesics
$\gamma_n\to \gamma$, the entropy of the limit is a weak limit of entropies,
$\mu_{\gamma_n}\rightharpoonup\mu_\gamma$. 
The later follows from property \ref{limit-convex} on page \pageref{limit-convex}. 
\end{enumerate}

The next statement is similar to \ref{ext-mono} and \ref{ext-conv}; 
it makes a new pre-quasigeodesic out of two. 
It will be used in the next section to construct \emph{quasigeodesics}.

\begin{thm}{\bf Extension.}\label{ext-pqg}
Let $A\in\Alex$, $\gamma_1\:[a,\infty)\to A$ and $\gamma_2\:[b,\infty)\to A$ be
two pre-quasigeodesics. 
Assume $$a\le b,\ \ \gamma_1(b)=\gamma_2(b),\ \  \gamma^-_1(b)\ \
\text{is polar to}\ \ \gamma^+_2(b)\ \  \text{and}\ \ |\gamma^+_2(b)|\le|\gamma^-_1(b)|$$ 
then its joint
$$\gamma\:[a,\infty)\to A,\ \ \gamma(t)=\l[
\begin{matrix}
\gamma_1(t)&\text{if}&t\le b\\
\gamma_2(t)&\text{if}&t\ge b
\end{matrix}\right.$$
is a pre-quasigeodesic.
Moreover, its entropy is defined by $$\mu_\gamma|_{(a,b)}=\mu_{\gamma_1},\ \
\mu_\gamma|_{(b,c)}=\mu_{\gamma_2}\ \ \text{and}\ \
\mu_\gamma(\{b\})=\ln|\gamma^+(b)|-\ln|\gamma^-(b)|.$$
\end{thm}

\Proof. The same as for \ref{ext-mono}.
\qeds

\begin{thm}{\bf Existence.} 
Let $A\in\Alex$, $x\in A$ and $\xi\in \Sigma_x$. 
Then there is a pre-quasigeodesic $\gamma\:[0,\infty)\to A$ such that
$\gamma(0)=x$ and $\gamma^+(0)=\xi$.
\end{thm}

\Proof. 
Let us choose for each point $x\in A$ and each direction $\xi\in \Sigma_x$ a convex
curve $\beta_\xi\:[0,\infty)\to A$ such that $\beta_\xi(0)=x$,
$\beta_\xi^+(0)=\xi$. 
If $v=r\xi$, then set
$$\beta_v(t)=\beta_\xi(r t).$$
Clearly $\beta_v$ is convex if $0\le r\le 1$.

Let us construct a convex curve $\gamma_\varepsilon \: [0,\infty)\rightarrow M$
such that there is a representation of $[0,\infty)$ as  a countable union of disjoint
half-open intervals $[a_i,\bar a_i)$, such that $|\bar a_i-a_i|\le \eps$ and for
any $t\in [a_i,\bar a_i)$ we have
$$|\gamma_\eps^+(a_i)|\ge |\gamma_\eps^+(t)| \ge
(1-\eps)\cdot|\gamma_\eps^+(a_i)|.\eqno(*)$$
Moreover, for each $i$, the curve $\gamma_\varepsilon^{a_i}\:[0,\infty)\to A$,
$$\gamma_\varepsilon^{a_i}(t)=\gamma_\varepsilon\left({a_i}+ \frac{t}{|\gamma_\varepsilon^+({a_i})|}\right)$$
is also convex.

Assume we already can construct $\gamma_\eps$ in the interval $[0,t_{\max})$, and cannot do it any further. 
Since $\gamma_\eps$ is 1-Lipschitz, we can extend it continuously to
$[0,t_{\max}]$.
Use lemma~\ref{lem:polar} to
construct a vector $v^*$ polar to $\gamma^-_\eps(t_{\max})$ with $|v^*|\le
|\gamma^-_\eps(t_{\max})|$.
Consider the joint of $\gamma_\eps$ with a short half-open segment of $\beta_v$, a longer curve with the desired property. This is a contradiction.

Let $\gamma$ be a partial limit of $\gamma_{\eps}$ as $\eps\to0$.
From property~\ref{limit-convex} on page~\pageref{limit-convex}, we get that for
almost all $t$ we have $|\gamma^+(t)|=\lim|\gamma_{\eps_n}^+(t)|$.
Combining this with inequality $(*)$ shows that for any $a\ge0$
$$\gamma^{a}(t)=\gamma\left({a}+\frac{t}{|\gamma^+({a})|}\right)$$
is convex.
\qeds

\subsection{Step 3: Quasigeodesics}\label{step3}

We will construct quasigeodesics in an $m$-dimensional Alexandrov's space, assuming
we already have such a construction in all dimensions $<m$. 
This construction is much easier for the case of an Alexandrov's space with only
\textit{$\delta$-strained points}; in this case we construct a sequence of special pre-quasigeodesics only by
extending/chopping procedures (see below) and then pass to the limit.
In a general Alexandrov's space we argue by contradiction, we assume that $\Omega$ is a maximal
open set such that for any initial data one can construct an
$\Omega$-quasigeodesic 
(that is, a pre-quasigeodesic with zero {\it entropy} on
$\Omega$, see \ref{def:entropy}), 
and arrive at a contradiction with the assumption
$\Omega\not=A$.

The following extension and chopping procedures are essential in the construction:

\begin{thm}{\bf Extension procedure.}\label{extension} Given a pre-quasigeodesic
$\gamma\:[0,t_{\max})\to A$ we can extend it as a pre-quasigeodesic
$\gamma\:[0,\infty)\to A$ so that 
$$\mu_\gamma(\{t_{\max}\})=0.$$
\end{thm}
\Proof. 
Let us set $\gamma(t_{\max})$ to be the limit of $\gamma(t)$ as $t\to t_{\max}$
(it exists since pre-quasigeodesics are Lipschitz).

From Milka's lemma~\ref{lem:milka}, we can construct a vector 
$\gamma^+(t_{\max})$ which is polar to $\gamma^-(t_{\max})$
and such that
$|\gamma^+(t_{\max})|=|\gamma^-(t_{\max})|$.
Then extend $\gamma$ by a pre-quasigeodesic in the direction $\gamma^+(t_{\max})$. 
By \ref{ext-pqg}, we get
$$\mu_\gamma\{t_{\max}\}=\ln|\gamma^+(t_{\max})|-\ln|\gamma^-(t_{\max}
)|=0.$$
\qedsf

\begin{thm}{\bf Milka's lemma} {\it (existence of the polar direction)}.
\label{lem:milka}
For any unit vector $\xi\in \Sigma_p$ there is a polar unit vector $\xi^*$;
that is, there is
$\xi^*\in \Sigma_p$ such that 
$$\<\xi,v\>+\<\xi^*,v\>\ge 0$$
for any $v\in T_p$.
\end{thm}

The proof is taken from \cite{milka:poly1}. That is the only instance where we use
existence of quasigeodesics in lower dimensional spaces.

\Proof. Since $\Sigma_p$ is an Alexandrov's $(m-1)$-space with curvature $\ge 1$, given
$\xi\in \Sigma_p$ we can construct a quasigeodesic in $\Sigma_p$ of length
$\pi$, starting at $\xi$; the comparison inequality (theorem
\ref{thm:eq-def-qg}(\ref{comp-inq})) implies that if $\xi^*$ is the other end of
this quasigeodesic then
$$|\xi\,\eta|_{\Sigma_q}+|\eta\,\xi^*|_{\Sigma_q}=\angle(\xi,\eta)+\angle(\eta,
\xi^*)\leq\pi\ \ \text{for all}\ \  \eta\in \Sigma_p.$$ 
The later is
equivalent to the statement that $\xi$ and $\xi^*$ are polar in $T_p$.
\qeds

\begin{thm}{\bf Chopping procedure.}\label{chopping} Given a pre-quasigeodesic
$\gamma\:[0,\infty)\to A$, for any $t\ge 0$ and $\eps>0$ there is $\bar t>t$ such
that 
$$\mu_\gamma\l((t,\bar t)\r) <\eps[\theta+\bar t-t],\ \ \bar t-t<\eps,\ \ \theta<\eps,$$
where 
$$\vartheta=\vartheta(t,\bar t) =
\angle\l(\gamma^+(t),\uparrow_{\gamma(t)}^{\gamma(\bar t)}\r).$$
\end{thm}
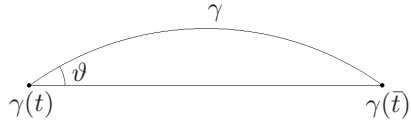
\begin{lpic}[t(-0mm),b(0mm),r(0mm),l(0mm)]{pics/angle(0.2)}
\lbl[t]{8,8;$\gamma(t)$}
\lbl[t]{245,8;$\gamma(\bar t)$}
\lbl[l]{33,21;$\vartheta$}
\lbl[b]{130,55;$\gamma$}
\end{lpic}

\Proof. For all sufficiently small $\tau>0$ we have $$\theta(t,t+\tau)<\eps$$
and from convexity of $\gamma^t$ it follows that 
$$\mu\l((t,t+\tau/3)\r) <C{\cdot}\theta^2(t,t+\tau).$$
The following exercise completes the proof.
\qeds

\begin{thm}{\bf Exercise.}
Let the functions $h,g\:\RR_+\to \RR_+$ be such that for any sufficiently small
$s$,
$$h(s/3)\le g^2(s),\ s\le g(s)\ \hbox{and}\ \lim_{s\to0} g(s)=0.$$
Show that for any $\varepsilon>0$ there is  $s>0$ such that
$$h(s)< 10{\cdot}g^2(s)\ \hbox{and}\ g(s)\le\varepsilon.$$
\end{thm}

\noi{\bf Construction in the $\mathbf\delta$-strained case.}
From the extension procedure, it is sufficient to construct a quasigeodesic
$\gamma\:[0,T)\to A$ with any given initial data $\gamma^+(0)=\xi\in \Sigma_p$ for some
positive $T=T(p)$. 

\textit{The plan:} Given $\eps>0$, we first construct a pre-quasigeodesic
$$\gamma_\eps\:[0,T)\to A,\ \ \  \gamma_\eps^+(0)=\xi$$ 
such that one can present $[0,T)$ as  a countable union of disjoint
half-open intervals
$[a_i,\bar a_i)$ with the following property ($\theta$ is defined in the chopping procedure~\ref{chopping}):
$$\mu\l([a_i,\bar a_i)\r) <\eps{\cdot}\theta(a_i,\bar a_i),\ \ \bar a_i-a_i<\eps,\ \
\theta(a_i,\bar a_i)<\eps.\eqno(\star)$$
Then we show that the entropies $\mu_{\gamma_\eps}([0,T))\to 0$ as $\eps\to0$
and passing to a partial limit of $\gamma_\eps$ as $\eps\to0$ we get a
quasigeodesic.

\bigskip
\noi\textit{Existence of $\gamma_\eps$:} Assume that we already can construct
$\gamma_\eps$ on an interval $[0,t_{\max})$, $t_{\max}<T$ and cannot construct it
any further, then applying the extension procedure~\ref{extension} for
$\gamma_\eps\:[0,t_{\max})\to A$ and then chopping it (\ref{chopping}) starting
from $t_{\max}$, we get a longer curve with the desired property; that is a contradiction.

\bigskip
\noi\textit{Vanishing entropy:} From $(\star)$ we have that 
$$\mu_{\gamma_\eps}([0,T))<\eps\cdot\l[T+\sum_i\theta(a_i,\bar a_i)\r].$$
Therefore, to show that $\mu_{\gamma_\eps}([0,T))\to 0$, it only remains to show
that the sum
\[\sum_i\theta(a_i,\bar a_i)\le \Const\] 
where $\Const$ is independent
of $\eps$.

That will be the only instance, where we apply that  $p$ is $\delta$-strained for a small
enough $\delta$.

Note that there is $\eps=\eps(\delta)\to0$ as $\delta\to0$ and
$T=T(p)>0$ such that there is a finite collection of points $\{q_k\}$ which
satisfy the following property: for any $x\in B_T(p)$ and $\xi\in \Sigma_x$ there is
$q_k$ such that
$\angle (\xi,\uparrow_x^{q_k})<\eps$.
Moreover, we can assume $\dist_{q_k}$ is $\lambda$-concave in $B_T(p)$ for some
$\lambda>0$.

Since the distance functions $\dist_{q_k}$ are 1-Lipschitz and $\lambda$-concave in $B_T(p)$,
for any convex curve $\gamma\:[0,T)\to B_T(p)\i A$, the measures
$\chi_k$ on $[0,T)$, defined by
$$\chi_k((a,b))=(\dist_{q_k}\circ\gamma)^-(b)-(\dist_{q_k}
\circ\gamma)^+(a)+\lambda{\cdot}(b-a),$$
are positive and their total mass is bounded by $\lambda T+2$.

Let $x\in B_T(p)$, and $\delta$ be small enough. 
Then for any two directions
$\xi,\nu\z\in \Sigma_x$ there is $q_k$ which satisfies the following property:
$$\tfrac1{10}{\cdot}\angle_x (\xi,\nu)\le \d_x\dist_{q_k}(\xi)-\d_x\dist_{q_k}(\nu)\ \ \
\text{and}\ \ \ \d_x\dist_{q_k}(\nu)\ge0. \eqno(*)$$
\label{inq:*}
Substituting in this inequality 
$$\xi=\gamma^+(a_i)/|\gamma^+(a_i)|,\ \ \
\nu=\uparrow_{\gamma(a_i)}^{\gamma(\bar a_i)},$$ 
and applying lemma~\ref{lem:angle-d}, we get 
$$\theta(a_i,\bar a_i)
=\angle(\xi,\nu)
\le 
10{\cdot}\sum_n\chi_k([a_i,\bar a_i)).$$
Therefore 
$$\sum_i\theta(a_i,\bar a_i)
\le 10{\cdot}N{\cdot}(\lambda T+2),$$
where $N$ is the number of points in the collection $\{q_k\}$.
\qeds

\begin{thm}{\bf Lemma.} \label{lem:angle-d} Let $A\in\Alex$, $\gamma\:[0,t]\to A$ be a convex curve $|\gamma^+(0)|=1$ 
and $f''\le\lambda$ for $\lambda\ge0$.
Set $p=\gamma(0)$, $q=\gamma(t)$, \
$\xi=(\gamma)^+(0)$ and $\nu=\uparrow_p^q$. 
Then
$$\d_p f(\xi)-\d_p f(\nu)\le
(f\circ\gamma)^+(0)-(f\circ\gamma)^-(t)+\lambda{\cdot}t,$$ 
provided that $\d_p f(\nu)\ge 0$.
\end{thm}

\begin{lpic}[t(2mm),b(2mm),r(0mm),l(0mm)]{pics/a35(0.4)}
\lbl[tr]{0,0;$p$}
\lbl[t]{121,-1;$q$}
\lbl[rb]{21,16;$\xi$}
\lbl[t]{25,-1;$\nu$}
\lbl[b]{70,22;$\gamma$}
\end{lpic}

\Proof. Clearly, 
$$f(q)\leq f(p)+\d_p f(\nu){\cdot}|p q|+\tfrac\lambda2{\cdot}
|p q|^2\le
f(p)+\d_p f(\nu)t+\tfrac\lambda2{\cdot}t^2.$$ 
On the other hand, 
$$f(p)\le f(q)-(f\circ\gamma)^-(t){\cdot}t+\tfrac\lambda2{\cdot}t^2.$$ 
Clearly, $\d_p f(\xi)= (f\circ\gamma)^+(0)$, whence the result.
\qeds

\noi{\bf What to do now?}
We have just finished the proof for the case, where all points of $A$ are $\delta$-strained. 
From this proof it follows that if we denote by $\Omega_\delta$ the subset of all
$\delta$-strained points of $A$ (which is an open everywhere dense set, see
\cite[5.9]{BGP}), then for any initial data one can construct a pre-quasigeodesic
$\gamma$ such that $\mu_\gamma(\gamma^{-1}(\Omega_\delta))=0$.
Assume $A$ has no boundary; set $\mathfrak C=A\backslash
\Omega_\delta$. 
In this case it seems unlikely that we hit $\mathfrak C$ by shooting a pre-quasigeodesic in a generic direction.
If we could prove that it almost never happens, then we obtain existence of
quasigeodesics in all directions as the limits of quasigeodesics in generic directions
(see property~\ref{pr:limit-qg} on page~\pageref{pr:limit-qg}) and passing to doubling in case $\partial A\not=\emptyset$.
Unfortunately, we do not have any tools so far to prove such a thing%
(It could be possible if we would have an analog of the Liouvile theorem for
``pre-quasigeodesic flow''.)

Instead we generalize inequality $(*)$ on page \pageref{inq:*}.

\begin{thm}{\bf The $\mathbf{(*)}$ inequality.} \label{inq:di-inq}
Let $A\in\Alex^m(\kappa)$ and $\mathfrak C\i A$ be a closed subset.
Let $p\in \mathfrak C$ be a point with  $\delta$-maximal $\vol_{m-1} \Sigma_p$;
that is
$$\vol_{m-1} \Sigma_p+\delta
>
\inf_{x\in\mathfrak C}\{\vol_{m-1} \Sigma_p\}.$$ 
Then, if $\delta$ is small enough, there is a finite set of points $\{q_i\}$ and
$\eps>0$, such that for any $x\in\mathfrak C\cap \bar B_\eps(p)$ and any pair of directions $\xi\in \Sigma_x\mathfrak C$%
\footnote{$\Sigma_x\mathfrak C$ is defined on page~\pageref{U_pX}.} 
and $\nu\in \Sigma_x$ we can choose $q_i$ so that 
$$\tfrac{1}{10}{\cdot}\angle_x(\xi,\nu)
\le 
\d_x\dist_{q_i}(\xi)-\d_x\dist_{q_i}(\nu)\ \ \
\text{and}\ \ \ \d_x\dist_{q_k}(\nu)\ge0.$$
\end{thm}

\Proof. 
We can choose $\eps>0$ so small that for any $x\in \bar B_\eps(p)$, 
$\Sigma_x$ is almost bigger than $\Sigma_p$.%
\footnote{that is for small $\delta>0$ there is a
map $f\:\Sigma_p\to \Sigma_x$ such that $|f(x)f(y)|>|x y|-\delta$.} 
Since $\vol_{m-1} \Sigma_p$ is almost maximal we get that for any 
$x\in  \mathfrak C\cap \bar B_\eps(p)$, $\Sigma_x$ is almost isometric to $\Sigma_p$.
In particular, if one takes a set $\{q_i\}$ so that directions $\uparrow_p^{q_i}$
form a sufficiently dense set and $\angle q_i p q_j\approx\tilde\angle_\kappa q_i
p q_j$, then directions $\uparrow_x^{q_i}$ will form a sufficiently dense set in
$\Sigma_x$ for all $x\in\mathfrak C\cap \bar B_\eps(p)$.

Note that for any $x\in  \mathfrak C\cap \bar B_\eps(p)$ and $\xi\in
\Sigma_x\mathfrak C$,  there is an almost isometry $\Sigma_x\to \Sigma(\Sigma_\xi \Sigma_x)$ such
that $\xi$ goes to north pole of the spherical suspension $\Sigma(\Sigma_\xi \Sigma_x)=\Sigma_\xi
T_x$.(Otherwise, taking a point $y\in \mathfrak C$, close to $x$ in direction
$\xi$ we would get that $\vol_{m-1}\Sigma_y$ is essentially bigger than
$\vol_{m-1}\Sigma_x$, which is impossible since both are almost equal to
$\vol_{m-1}\Sigma_p$.)

Using these two properties, we can find $q_i$ so that 
$\uparrow^\nu_\xi\approx\uparrow_\xi^{\uparrow_x^{q_i}}$ 
in $\Sigma_\nu (\Sigma_x A)$ and
$\angle(\xi,\uparrow_x^{q_i})>\tfrac\pi2$, hence the statement follows.
\qeds

Now we are ready to finish construction in the general case. 
Let us define a subtype
of pre-quasigeodesics:

\begin{thm}{\bf Definition.}
Let $A\in \Alex$ and $\Omega\i A$ be an open subset. 
A pre-quasigeodesic $\gamma\:[0,T)\to A$ is called $\Omega$-quasigeodesic if  its
entropy vanishes on $\Omega$; 
that is, 
$$\mu_\gamma(\gamma^{-1}(\Omega))=0$$
\end{thm}

From property~\ref{lim-entropy} on page~\pageref{lim-entropy}, it follows that
the limit of $\Omega$-quasigeodesics is a $\Omega$-quasigeodesic. 
Moreover, if for any initial data we can construct an
$\Omega$-quasigeodesic and an $\Omega'$-quasigeodesic, then it is possible to construct an
$\Omega\cup \Omega'$-quasigeodesic for any initial data; for $\Upsilon\Subset\Omega\cup \Omega'$, $\Upsilon$-quasigeodesic can be constructed by joining together
pieces of $\Omega$ and $\Omega'$-quasigeodesics and $\Omega\cup \Omega'$-quasigeodesic can be constructed as a limit of $\Upsilon_n$-quasigeodesics as $\Upsilon_n\to \Omega\cup \Omega'$.

Let us denote by $\Omega$ the maximal open set such that for any initial data
one can construct an $\Omega$-quasigeodesic. 
We have to show then that $\Omega=A$. 

Let $\mathfrak C=A\backslash \Omega$, and let $p\in \mathfrak C$ be the point
with almost maximal $\vol_{m-1} \Sigma_p$. 
We will arrive to a contradiction by constructing  a $B_\eps(p)\cup
\Omega$-quasigeodesic for any initial data. 

Choose a finite set of points $q_i$ as in \ref{inq:di-inq}.
Given $\eps>0$, it is enough to construct an $\Omega$-quasigeodesic
$\gamma_\eps\:[0,T)\to A$, for some fixed $T>0$ with the given initial data $x\in
\bar B_\eps(p)$, $\xi\in \Sigma_x$, such that the entropies
$\mu_{\gamma_\eps}((0,T))\to 0$ as $\eps\to0$.
 
The $\Omega$-quasigeodesic $\gamma_\eps$ which we are going to construct will
have the following property: one can present $[0,T)$ as a countable union of disjoint
half-open intervals $[a_i,\bar a_i)$ such that 
$$\text{if}\ \ \
\frac{\gamma^+(a_i)}{|\gamma^+(a_i)|}\in \Sigma_{\gamma(a_i)}\mathfrak C \ \ \
\text{then}\ \ \   \mu_\gamma([a_i,\bar a_i))\le \eps{\cdot}\theta(a_i,\bar a_i)$$ 
and
$$\text{if}\ \ \ \frac{\gamma^+(a_i)}{|\gamma^+(a_i)|}\not\in
\Sigma_{\gamma(a_i)}\mathfrak C\ \ \ \text{then}\ \ \ \mu_\gamma([a_i,\bar a_i))=0$$

Existence of $\gamma_\eps$ is being proved the same way as in the $\delta$-strained case,
with the use of one additional observation:
if $$\frac{\gamma^+(t_{\max})}{|\gamma^+(t_{\max})|}\not\in
\Sigma_{\gamma(a_i)}\mathfrak C$$ then any $\Omega$-quasigeodesic in this direction
has zero entropy for a short time.

Then, just as in the $\delta$-strained case, applying
inequality~\ref{inq:di-inq} we get that $\mu_{\gamma_{\eps}}(0,T)\to0$ as
$\eps\to 0$. 
Therefore, passing to a partial limit $\gamma_\eps\to\gamma$ gives a
$B_\eps(p)\cup\Omega$-quasigeodesic $\gamma\:[0,T)\to A$ for any initial data in
$B_\eps(p)$.
\qeds

\subsection{Quasigeodesics in extremal subsets.}\label{qg-extrim}

The second part of theorem~\ref{thm:exist-qg} follows from the above construction, but we have to
modify Milka's lemma~\ref{lem:milka}:

\begin{thm}{\bf Extremal Milka's lemma.} Let $E\i T_p$ be an extremal subset of a tangent cone then for any vector $v\in E$ there is a polar vector $v^*\in E$
such that $|v|=|v^*|$.
\end{thm}

\Proof. Set $X=E\cap \Sigma_p$. If $\Sigma_\xi X\not=\emptyset$ then the proof is the
same as for the standard Milka's lemma; it is enough to choose a direction in
$\Sigma_\xi X$ and shoot a quasigeodesic $\gamma$ of length $\pi$ in this direction such that $\gamma\i X$ ($\gamma$ exists from the induction hypothesis). 

Since $E$ is extremal,
if $X=\{\xi\}$ then $B_{\pi/2}(\xi)=\Sigma_p$. 
Therefore $\xi$ is polar to itself.

Otherwise, if $\Sigma_\xi X=\emptyset$ and $X$ contains at least two points, choose $\xi^*$ to be closest point in $X\backslash\xi$ from $\xi$. 
Since $X\i \Sigma_p$ is extremal we have that for any $\eta\in \Sigma_p$ $\angle_{\Sigma_p}\eta\xi^*\xi\le\tfrac\pi2$ and since $\Sigma_\xi X=\emptyset$ we have $\angle_{\Sigma_p}\eta\xi\xi^*\le\tfrac\pi2$.
Therefore, from triangle comparison we have
$$|\xi\eta|_{\Sigma_p}+|\eta\xi^*|_{\Sigma_p}
=
\angle(\xi,\eta)+\angle(\eta,\xi^*)\le\pi$$
\qedsf

\end{document}